\documentclass[reqno]{amsart}
\usepackage{amsthm,amsmath,amssymb,amscd, mathrsfs,graphics, hyperref}
\usepackage{diagrams}
\newlength{\marg}
\usepackage[lmargin=\marg,tmargin=\marg,bmargin=\marg,rmargin=\marg]{geometry}

\newarrow {Mapsto} |--->
\newarrow {Iso} ---->
\newcommand{\Ad}{\operatorname{Ad}}
\newcommand{\Aut}{\operatorname{Aut}} 
\newcommand{\BG}{\B G}
\newcommand{\B}{\cb}

\newcommand{\Ch}{\mathrm{ch}} 
\newcommand{\DD}{\Delta} 

\newcommand{\Dt}{\widetilde{D}}
\newcommand{\F}{\mathbf{F}} 
\newcommand{\Gac}{\Gammac} 
\newcommand{\Gammac}{{\Gamma_\mathrm{cut}}}

\newcommand{\Gatc}{{\Gat_{\mathrm{cut}}}} 
\newcommand{\Gat}{{\widetilde{\Gamma}}} 
\newcommand{\Ga}{\Gamma} 
\newcommand{\Gb}{\overline{G}}
\newcommand{\Hom}{{\mathscr{H}\!\mathrm{om}}}
\newcommand{\IF}{\mathrm{C}} 
\newcommand{\IIF}{\widetilde{\IF}} 

\newcommand{\Ind}{\operatorname{Ind}} 
\newcommand{\In}{\operatorname{In}} 
\newcommand{\Irrep}{\operatorname{Irrep}} 

\newcommand{\I}{\mathbf{I}} 
\newcommand{\Kmgnloc}[1]{K(\MM_{g,n}#1)_{\loc}} 
\newcommand{\Kgatcloc}{K(\MM_{\Gatc})_{\loc}} 
\newcommand{\LG}{P} 

\newcommand{\MM}{\M^G}

\newcommand{\M}{\overline{\cM}}
\newcommand{\N}{\overline{\mathscr{N}}}
\newcommand{\Out}{\operatorname{Out}} 

\newcommand{\Rep}{\operatorname{Rep}} 
\newcommand{\Res}{\operatorname{Res}} 
\newcommand{\Tdch}{\mathrm{Td}^{\vee}}  
\newcommand{\Td}{\mathrm{Td}}  

\newcommand{\V}{\mathbf{V}} 
 
\newcommand{\bCh}{\mathbf{ch}} 
\newcommand{\bGacgn}{{\bGa_{cut, g, n}}}
\newcommand{\bGac}{{\bGa_{cut}}} 
\newcommand{\bGatcgn}{{\bGat_{cut,g,n}}}
\newcommand{\bGatc}{{\bGat_{cut}}} 
\newcommand{\bGat}{{\widetilde{\boldsymbol{\Gamma}}}} 
\newcommand{\bGa}{{\boldsymbol{\Gamma}}} 

\newcommand{\bmb}{{\overline{\mathbf{m}}}}

\newcommand{\bm}{\mathbf{m}}
\newcommand{\bone}{\mathbf{1}} 
\newcommand{\brk}{\mathbf{rk}} 

\newcommand{\cF}{{\cf}}
\newcommand{\cM }{\mathscr{M}}

\newcommand{\cb}{{\mathscr B}}
 
\newcommand{\cc}{{\mathscr C}}

\newcommand{\ce}{{\mathscr E}}
\newcommand{\cfb}{\overline{\cf}}
\newcommand{\cf}{{\mathscr F}}

\newcommand{\cg}{\mathscr{G}} 

\newcommand{\cltm}{\widetilde{\cl}_-}
\newcommand{\cltp}{\widetilde{\cl}_+}
\newcommand{\clt}{\widetilde{\cl}} 
\newcommand{\cl}{{\mathscr L}} 
 
\newcommand{\conor}{{\mathfrak{C}}}
\newcommand{\co}{{\mathscr O}}

\newcommand{\crrt}{\widetilde{\crr}} 
\newcommand{\crr}{{\mathscr R}}
\newcommand{\csh}{{\mathbb{S}}}  
\newcommand{\cst}{\mathfrak{S}} 
\newcommand{\cs}{{\mathscr S}} 

\newcommand{\ct}{{\mathscr T}}

\newcommand{\cx}{\mathscr{X}} 
\newcommand{\cX}{\cx}
\newcommand{\cy}{\mathscr{Y}} 
\newcommand{\cY}{\cy}
\newcommand{\cz}{\mathscr{Z}} 
 
\newcommand{\dsand}{\quad \text{ and } \quad} 

\newcommand{\ee}{\varepsilon} 

\newcommand{\ga}{\gamma}
\newcommand{\gb}{\overline{\gamma}}
\newcommand{\gen}[1]{\langle{#1}\rangle} 
\newcommand{\gt}{{\tilde{g}}}  
 
\newcommand{\hodt}{\widetilde{\hod}} 
\newcommand{\hod}{\mathbb{E}} 

\newcommand{\id}{\operatorname{id}} 
\newcommand{\ie}{{i.e.,}} 
\newcommand{\ih}{\widehat{i}}

\newcommand{\itt}{\widetilde{i}}

\newcommand{\jjm}{{j_-}}
\newcommand{\jjpm}{{j_\pm}}
\newcommand{\jjp}{{j_+}}

\newcommand{\jt}{\widetilde{j}}
\newcommand{\kat}{\tilde{\kappa}} 
\newcommand{\km}{{k_-}}
\newcommand{\kpm}{{k_\pm}}
\newcommand{\kp}{{k_+}}

\newcommand{\loc}{\mathfrak {H}} 
\newcommand{\mb}{{\overline{m}}} 

\newcommand{\mmm}{{m_-}} 

\newcommand{\mpp}{{m_+}} 

\newcommand{\mut}{\widetilde{\mu}}
\newcommand{\nc}{\mathbb{C}} 
\newcommand{\nor}{{\mathfrak{N}}}  
\newcommand{\nq}{\mathbb{Q}}
 
\newcommand{\nz}{\mathbb{Z}}

\newcommand{\one}{1} 

\newcommand{\pibar}{\overline{\pi}}
\newcommand{\pic}{\operatorname{Pic}}
\newcommand{\pih}{\widehat{\pi}}

\newcommand{\pr}{\operatorname{pr}}
\newcommand{\psim}{\psi_-}
\newcommand{\psip}{\psi_+}
\newcommand{\psitm}{\psit_-}
\newcommand{\psitp}{\psit_+}
\newcommand{\psit}{\widetilde{\psi}}
\newcommand{\pt}{\widetilde{p}}

\newcommand{\rhot}{\widetilde{\rho}}

\newcommand{\rk}{\mathrm{rk}} 

\newcommand{\rpp}{r_+}

\newcommand{\sh}{\widehat{s}}
\newcommand{\sigbar}{\overline{\sigma}}
\newcommand{\sigb}{\boldsymbol{\sigma}} 
\newcommand{\sigm}{{\sig_-}}
\newcommand{\sigpm}{\sigma_{\pm}}
\newcommand{\sigp}{{\sig_+}}
\newcommand{\sig}{\sigma} 

\newcommand{\st}{\operatorname{st}}

\newcommand{\tauh}{\widehat{\tau}}
\newcommand{\taut}{\widetilde{\tau}}

\newcommand{\td}{\mathbf{td}}  %

\newcommand{\tr}{\operatorname{Tr}}

\newcommand{\vac}{\mathbf{1}}

\newcommand{\xm}{x_-}
\newcommand{\xp}{x_+}

\renewcommand{\L}{{\mathbb{L}}} 

\renewcommand{\th}{\widehat{t}}
\renewcommand{\tt}{\widetilde{t}}
 
\newtheorem{result}{Main Results} 

\newtheorem{thm}{Theorem}[section] 
\newtheorem{lm}[thm]{Lemma}
\newtheorem{prop}[thm]{Proposition} 
\newtheorem{crl}[thm]{Corollary}

\theoremstyle{definition}
\newtheorem{rem}[thm]{Remark} 

\newtheorem{df}[thm]{Definition} 

\newtheorem{df-pr}[thm]{Definition-Proposition}

\theoremstyle{remark} 
 \renewcommand{\thenota}{\kern-1ex}

\begin{document}
\title[Riemann-Hurwitz theorem] {A representation-valued relative Riemann-Hurwitz theorem and the Hurwitz-Hodge bundle}


\author [T. J. Jarvis]{Tyler J. Jarvis} \address {Department of
Mathematics, Brigham Young University, Provo, UT 84602, USA}
\email{jarvis@math.byu.edu} 

\author [T. Kimura]{Takashi Kimura} \address {Department of
Mathematics and Statistics; 111 Cummington Street, Boston University;
Boston, MA 02215, USA } \email{kimura@math.bu.edu} 
\thanks{Research of both authors was partially supported by NSF grant DMS-0605155} 

\date{\today}

\begin{abstract}

We provide a formula describing the $G$-module structure of the Hurwitz-Hodge
bundle for  admissible $G$-covers in terms of the Hodge bundle of the base
curve, and more generally, for describing the $G$-module structure of the
push-forward to the base of any sheaf on a family of admissible $G$-covers.
This formula can be interpreted as a representation-ring-valued
relative Riemann-Hurwitz formula for families of admissible $G$-covers. 
\end{abstract}

\maketitle \setcounter{tocdepth}{1}

\tableofcontents

\section{Introduction}

Much of the motivation for this paper arises from the classical Riemann-Hurwitz
Theorem, which we now briefly recall. 

 Let $\pih: E\rTo C$
be a surjective morphism between smooth projective curves (or compact Riemann
surfaces), where $E$ has genus $\gt$ and $C$ has genus $g$. Of particular
interest to us is the case where $E$ has 
$N$ connected components and
$C$ is connected. The \emph{Riemann-Hurwitz formula} is
\begin{equation}\label{eq:OriginalRH}
\gt = N + (g - 1) \deg(\pih) + S^{\pih},
\end{equation}
where $\deg(\pih)$ is the degree of the map $\pih$ and $S^{\pih}$ is defined by
\begin{equation*}
S^{\pih} := \sum_{q\in E} \frac{r(q)-1}{2},
\end{equation*}
where $r(q)$ is the ramification index of $q$. This formula is
important because it relates the global quantities $\gt$, $g$, 
$N$,
and
$\deg(\pih)$ to each other through $S^{\pih}$, which is the sum of local
quantities. 

The Riemann-Hurwitz formula has a particularly simple form when the smooth
complex curve $E$ has an 
action of a finite group $G$ such that every nontrivial element of $G$ has a finite fixed point set and $E\rTo^{\pih} C$ is the quotient map 
where $C$ is a smooth, connected curve of genus $g$.

Let $\{\,p_1,\ldots,p_n\,\} \subset C$ be the branch points of $\pih$.  The restriction of the curve $E$ to
$C-\{\,p_1,\ldots,p_n\,\}$ has a free $G$ action, \ie\ the 
restriction is a principal $G$-bundle over $C-\{\,p_1,\ldots,p_n\,\}$. 
It will be convenient to choose, for all $j\in \{1,\dots,n\}$, points $\pt_j$ in
$\pih^{-1}(p_j)$. This gives a pair $(E;\pt_1,\ldots,\pt_n)\rTo^{\pih} 
(C; p_1,\ldots,p_n)$, called a smooth, \emph{pointed, admissible $G$-cover} of genus
$g$ with $n$ marked points. The isotropy group of the point $\pt_j$ is a
cyclic group $\gen{m_j}$ of order $r_j$, where the generator $m_j$ acts on the tangent space $T_{\pt_j} E$ by multiplication by $\exp(2\pi i/r_j)$. Since $\pih^{-1}(p_j)$ is isomorphic to $G/\gen{m_i}$ as a $G$-set, we have
$S^{\pih} = \sum_{i=1}^n \frac{|G|}{r_i} \frac{r_i-1}{2}$. 
Furthermore, if $G_0$ is a subgroup of $G$ fixing a connected component of
$E$, then 
$N$ is the order of the set of cosets $G/G_0$. Therefore, the
Riemann-Hurwitz formula in this case can be written
\begin{equation}\label{eq:OriginalGRH}
\gt = \frac{|G|}{|G_0|}  + (g - 1) |G| + \sum_{i=1}^n S_{m_i},
\end{equation}
where, for all $i=1,\ldots,n$,
\begin{equation}
S_{m_i} := \frac{|G|}{2}
\left(1-\frac{1}{r_i} \right).
\end{equation}

There should be a generalization of the Riemann-Hurwitz formula, Equation
(\ref{eq:OriginalRH}), to families of curves
$\ce\rTo^{\pih}\cc\rTo^{\pibar}T$, where each fiber $\ce_t$ has genus $\gt$ and 
$N$ connected components,
and the fiber $\cc_t$ is a connected, genus-$g$ curve for all $t$ in
$T$. That is, we would like an equation in the K-theory of $T$
which, after taking the (virtual) rank, yields Equation
(\ref{eq:OriginalRH}). 
Although the answer is not known in general, 
we provide a complete answer in this paper for the case where 
$E$ forms an admissible $G$-cover for a finite group $G$.

We prove a generalization of Equation
(\ref{eq:OriginalGRH}) for flat, projective families of 
pointed  admissible $G$-covers $\ce\rTo^{\pih}\cc\rTo^{\pibar}T$,
which may also include curves and covers with nodal singularities. Our
generalization of the Riemann-Hurwitz formula (given in
Equation~(\ref{eq:BF}))  takes values in $\Rep(G)$, the
representation ring of  $G$. Since any such family of curves is the pull-back
of the universal family of $G$-curves  from the moduli space $\MM_{g,n}$ of
pointed admissible $G$-covers, we need only prove the $\Rep(G)$-valued
Riemann-Hurwitz  formula for the universal family of $G$-covers when the base
is  $T= \MM_{g,n}$. 

Every term in Equation (\ref{eq:OriginalGRH})
has a counterpart in our $\Rep(G)$-valued Riemann-Hurwitz formula. The term
$\gt$ corresponds to $R^1\pi_*\co_\ce$ (the dual Hurwitz-Hodge 
bundle), while the term $g$ corresponds to $R^1\pibar_*\co_\ce$ (the dual
Hodge bundle), and  for all $i\in\{1,\ldots,n\}$ the $S_{m_i}$ term corresponds to an element $\cs_{m_i}$ in K-theory, constructed from the tautological line bundles 
associated to the $i$-th marked point of $E$. In other
words, the $\Rep(G)$-valued Riemann-Hurwitz formula expresses the dual
Hurwitz-Hodge bundle in terms of other tautological K-theory classes on
$\MM_{g,n}$.

Our generalized Riemann-Hurwitz formula (Equation~(\ref{eq:BF})) has several
interesting and useful consequences.  Among other things, it allows one to
explicitly describe, without reference to the universal $G$-cover, the virtual class for
orbifold Gromov-Witten invariants for global quotients in degree zero.  It also yields a 
new differential equation which computes arbitrary descendant Hurwitz-Hodge integrals.  
These applications will be explored elsewhere \cite{JK2}.  

We note that Tseng \cite{Ts} has also obtained a differential equation which can be used
to calculate descendant potential functions of Hurwitz-Hodge integrals. However, he works directly in the Chow ring, while our differential equation is a consequence of our $\Rep(G)$-valued relative Riemann-Hurwitz formula, which lives in equivariant K-theory.  
It would be very interesting to make explicit the relationship between his results and ours.

\subsection*{Relation to the results of \cite{JKK2}}

The paper \cite{JKK2} gives a simple formula for the obstruction bundle of stringy cohomology and K-theory of a complex manifold $X$ with an action of a finite group $G$.  That formula allows one to completely describe the stringy cohomology of $X$ 
in a manner which does not involve any complex curves or $G$-covers.  It also yields a simple formula for the obstruction bundle of Chen-Ruan orbifold cohomology of an arbitrary Deligne-Mumford stack.

One interpretation of the present result is as a generalization of the results of \cite{JKK2} to stringy (or orbifold) Gromov-Witten theory.   That is, the results of \cite{JKK2} are essentially limited to the case of genus zero with three marked points, while the results of this paper are for general families of arbitrary genus and arbitrary numbers of marked points.  Applications and details of how this can be applied to orbifold Gromov-Witten invariants are developed in \cite{JK2}.

\subsection*{Overview of background and results}
We now describe the background and results of this paper in more detail. 

The moduli space $\M_{g,n}$ of stable curves of genus
$g$ with $n$ marked points is endowed with tautological vector bundles
arising from the universal curve $\cc\rTo^\pibar \M_{g,n}$ and its universal
sections $\M_{g,n}\rTo^{\sigbar_i}\cc$, where $i=1,\ldots,n$.  Pulling back
the relative dualizing sheaf $\omega$ of $\pibar$ by a tautological section
yields a line bundle $\cl_i = \sigbar_i^*\omega$ on $\M_{g,n}$ whose fiber
over a stable curve $(C; p_1,\ldots,p_n)$ is the cotangent line $T^*_{p_i}
C$. Its first Chern class is the tautological class $\psi_i := c_1(\cl_i)$
for all $i=1,\ldots,n$. The Hodge bundle $\hod := \pibar_*\omega$ is a
rank-$g$ vector bundle over $\M_{g,n}$ whose dual bundle is isomorphic to
$\crr := R^1\pibar_*\co_\cc$ by Serre duality. The Chern classes of the Hodge
bundle  give tautological classes $\lambda_j := c_j(\hod)$ for $j=0,\ldots,g$. 

For any finite group $G$, there is a natural generalization of the moduli
space of stable curves $\M_{g,n}$, namely the space $\MM_{g,n}$  of pointed
admissible  $G$-covers of genus $g$ with  $n$ marked points.  A pointed
admissible $G$-cover is a morphism $(E; \pt_1,\ldots,\pt_n)\rTo (C;
p_1,\ldots,p_n)$, where $(C; p_1,\ldots,p_n)$ is an $n$-pointed stable curve
of genus $g$ together with an admissible $G$-cover $E$, and for each $i \in
\{1,\dots,n\}$ $\pt_i$ is a choice of one point in the fiber over $p_i$. Since
the pointed curve $(C; p_1,\ldots,p_n)$ can be recovered from $(E;
\pt_1,\ldots,\pt_n)$ by taking the quotient by $G$, we will sometimes denote
the pointed $G$-cover by just $(E;\pt_1,\ldots,\pt_n)$.  
Basic properties of 
the stack
$\MM_{g,n}$ of pointed $G$-covers of genus $g$ with $n$ marked points
are described in \cite{JKK}.

An (unpointed) admissible $G$-cover consists of the same data except that one
forgets the marked points $\pt_1,\ldots,\pt_n$  on ${E}$ while retaining the
marked  points $p_1,\ldots,p_n$ on $C$.   The stack of (unpointed) admissible
$G$-covers of genus-$g$ curves with $n$ marked points is equivalent
\cite{ACV} to the stack  $\M_{g,n}(\BG)$ of stable  maps into the stack
$\BG$, where $\BG = [\mathrm{pt}/G]$ is the quotient stack of a single point
modulo the trivial action of $G$. For any $G$, there are forgetful morphisms
$\MM_{g,n}\rTo^t \M_{g,n}(\BG)\rTo^s \M_{g,n}$, where the first morphism forgets the
marked points $\pt_1,\ldots,\pt_n$ on the $G$-cover ${E}$, while the second
forgets the $G$-cover ${E}$. When $G$ is the trivial group, then both of 
these morphisms are isomorphisms.  

Consider the universal $G$-cover $\ce\rTo^\pi\MM_{g,n}$, whose fiber over a
point $[E,\pt_1,\dots,\pt_n]\in\MM_{g,n}$ is the $G$-cover $E$ itself. The projection
$\pi$ is $G$-equivariant, where $\ce$ inherits a tautological $G$-action, and
where the action on $\MM_{g,n}$ is trivial. The universal section $\sigma_i$
associated to the $i$-th marked point on the $G$-cover yields a tautological
line bundle $\clt_i := \sigma^*\omega_{\pi}$ for all $i=1,\ldots,n$, where
$\omega_{\pi}$ is the relative dualizing sheaf of $\pi$. That is, $\clt_i$ is
the line bundle on $\MM_{g,n}$ whose fiber over $(E; \pt_1,\ldots,\pt_n)$ is
the cotangent line $T^*_{\pt_i}E$. Similarly, let $\crrt$ be the
$G$-equivariant vector bundle $R^1\pi_* \co_\ce$ on $\MM_{g,n}$ whose fiber
over a point $(E; \pt_1,\ldots,\pt_n)$ is the $G$-module $H^1(E,\co_E)$,
where $\co_E$ is the structure sheaf of $E$. $\crrt$ has rank equal to the
genus
$\gt$ 
of the (possibly disconnected) curve $E$. We will call  the
$G$-equivariant dual bundle $\hodt:=\pi_* \omega_\pi$ on $\MM_{g,n}$ the
\emph{Hurwitz-Hodge bundle}. When $G$ is the trivial group, then $\clt_i$, $\crrt$, and $\hodt$ reduce to $\cl_i$, $\crr$, and $\hod$, respectively. 

The tautological classes are obtained by taking Chern classes of these
bundles. The main point of this paper is that the $G$-equivariant K-theory
class of the Hurwitz-Hodge bundle $\crrt$ (or $\hodt$) admits an explicit
description in terms of 
the usual Hodge bundle $\crr$ (or $\hod$) obtained from the universal curve
(rather than the universal $G$-cover) and the tautological line bundles
$\clt_i$ for all $i\in\{1,\ldots,n\}$.
This formula admits an interpretation as a $\Rep(G)$-valued Riemann-Hurwitz
formula for families. Because the Chern class of a vector bundle only depends
upon its class in K-theory, this means that all the tautological classes  can
be written explicitly  as pullbacks, by the forgetful map $\MM_{g,n}\rTo
\M_{g,n}$, of  the tautological classes $\psi_i$ and $\lambda_i$ from the
moduli space of curves $\M_{g,n}$. In other words, one can \emph{remove all
references to the universal $G$-cover from the computation of the Chern
character of $\crrt$.}  

More precisely, endowing $\MM_{g,n}$ with a trivial 
$G$-action,\footnote{There is another $G$-action present on $\MM_{g,n}$ as well, namely that arising from the action on the sections of $\ce$.  However, the trivial action will be more useful for our purposes and will be the action we use throughout this paper unless otherwise specified.} 
we can ``factorize''
the $G$-equivariant K-theory $K_G(\MM_{g,n})$ of  coherent sheaves 
on $\MM_{g,n}$ into a geometric part and a
representation-theoretic part, i.e.,  there is an algebra
isomorphism $$\Phi:K_G(\MM_{g,n})\rIso 
K(\MM_{g,n})\otimes\Rep(G),$$ where $K(\MM_{g,n})$ is the Grothendieck group of
coherent sheaves 
on $\MM_{g,n}$ and $\Rep(G)$ is the representation ring of
$G$.  This is because a $G$-equivariant 
sheaf over a base with a trivial $G$-action has a fiber-preserving $G$-action and, hence, an associated representation of $G$.

In the special case of $g=0$ and $n=3$, the stack $\MM_{0,3}$ is zero
dimensional, and so the $G$-equivariant $K$-theory of each connected
component is isomorphic to $\Rep(G)$.  In \cite[Lm 8.5]{JKK2} and
(independently) in \cite{Kan} there is a formula which gives a simple
expression for the representation $\Phi(\crrt)$ solely in terms of the
monodromies at the marked points $\pt_i$. The main result of this paper is a
generalization of that formula to all genera and to any number of marked points.

\begin{result}
We provide an explicit formula describing the $G$-equivariant structure of the K-theoretic push-forward of an arbitrary $G$-equivariant sheaf from a family of pointed, admissible $G$-covers to the base of the family.   That
is, for all finite groups $G$, for all genera $g\ge 0$, and for all $n\ge 0$
such that $2g-2+n > 0$, we give an explicit formula (Equation~(\ref{eq:FullBF}))  for the push-forward $\Phi(R\pi_* \cf) \in K(\MM_{g,n})\otimes\Rep(G)$ of any $G$-equivariant sheaf $\cf \in K_G(\ce)$.  Here, as above, $\ce\rTo^{\pi} \MM_{g,n}$ is the universal pointed admissible $G$-cover.  

As a corollary, we obtain a 
\emph{universal $\Rep(G)$-valued relative
Riemann-Hurwitz theorem}, that is, a simple, explicit formula
(Equation~(\ref{eq:BF})) for the dual Hurwitz-Hodge bundle $\Phi(\crrt) \in
K(\MM_{g,n})\otimes\Rep(G)$ in terms of 
the dual Hodge bundle $\crr$ pulled back from $\M_{g,n}$ and the tautological
line bundles $\clt_i$ for all $i\in \{1,\ldots,n\}$. 

The universal $\Rep(G)$-valued relative Riemann-Hurwitz formula allows us to
\begin{itemize}
\item Write a $\Rep(G)$-valued relative Riemann-Hurwitz formula for any pointed, flat, projective family of admissible $G$-covers.
\item Write down an explicit action of the automorphism group of $G$ on the
  $\Rep(G)$-valued relative Riemann-Hurwitz formula.
\item Write the tautological classes on $\MM_{g,n}$ in terms of tautological classes on $\M_{g,n}$. Moreover, we can do this in a way that also tracks the action of the group $G$.
\item Write $\Rep(G)$-valued generalizations of Mumford's identity.
\end{itemize}
\end{result}

Taking the rank of the $\Rep(G)$-valued Riemann-Hurwitz formula, Equation
(\ref{eq:BF}), yields the original Riemann-Hurwitz formula
(Equation~(\ref{eq:OriginalGRH})) 
for admissible $G$-covers.

The idea behind the proof of the Main Theorem (\ref{thm:FullBF}) is to apply
the Lefschetz-Riemann-Roch Theorem, a localization theorem in $G$-equivariant
$K$-theory, to the universal $G$-cover. This expresses the push-forward of
any $G$-equivariant sheaf on the universal $G$-cover over $\MM_{g,n}$ in
terms of the normal bundles to the fixed point loci which, for the universal
$G$-cover, can only occur at its punctures and nodes.

\subsection*{Structure of the paper}
The structure of this paper is as follows. 
In Section~\ref{sec:groups}, we recall some facts from representation theory,
equivariant K-theory and cohomology, and then prove a version of the
Lefschetz-Riemann-Roch theorem for stacks that we will need.
In the third section, we
review properties of the moduli spaces $\MM_{g,n}$ and $\M_{g,n}(\BG)$,
tautological bundles associated to them, and the universal $G$-curve and the
gluing morphisms on $\MM_{g,n}$. We also establish some basic properties of
the dual Hurwitz-Hodge bundle $\crrt$ with respect to gluing and forgetting
tails. We also introduce its associated tautological classes. In the fourth
section, we prove the Main Theorem and the 
$\Rep(G)$-valued Riemann-Hurwitz Theorem for families. As a corollary, we
prove a $\Rep(G)$-valued version of Mumford's Identity for the Hurwitz-Hodge
bundle.
In the fifth section, we calculate the $\Rep(G)$-valued Chern character of
the dual Hurwitz-Hodge bundle $\crrt$ on $\MM_{g,n}$. We then show that 
our $\Rep(G)$-valued Chern character of $\crrt$ reduces, in a special case,
to a formula for the ordinary Chern character of the Hurwitz-Hodge bundle
obtained by generalizing   Mumford's calculation of the ordinary Chern
character of the Hodge bundle on $\M_{g,n}$. 

\subsection*{Future directions for research}

The $\Rep(G)$-valued Riemann-Hurwitz formula allows one to use the structure
of the representation ring to obtain additional information about the
structure of its associated potential functions.  We hope to treat these in a
subsequent paper.

In this paper we have focused on developing formulas for global quotients.
It is natural to try extend these to general Deligne-Mumford stacks.    We
hope to treat this in a later paper.  Related to this, it would be very
interesting to generalize our results to the case where $G$ is a Lie group. 

Finally, the $\Rep(G)$-valued Riemann-Hurwitz Theorem, in nice cases,
can be regarded as a formula for a $G$-index theorem of the Dolbeault
operator $\overline{\partial}$ for families of Riemann surfaces, \ie\ as a
relative holomorphic Lefschetz theorem. It would be interesting to generalize
these results to other differential operators.

\subsection*{Acknowledgments}
We would like to thank the National Science Foundation (Grant DMS-0605155)
for their financial support and also the \emph{Recent Progress on the Moduli
Space of Curves} Workshop at the Banff International Research Station
for providing an opportunity to meet together to work on this project and to share our results. 
We would especially like to thank D.~Edidin for many helpful discussions and
suggestions. 

The second author would like to thank D.~Rohrlich for useful conversations
and to thank friends and colleagues for putting up with him during the
preparation of this paper. 
He would also like to thank the \emph{Aspects of Moduli} Workshop at the De
Giorgi Center at the Scuola Normale Superiore Pisa and the \emph{Geometry and
Physics} Program at the Hausdorff Research Institute for Mathematics where
much progress on an earlier version of this paper was made and the Institut
des Hautes \'Etudes Scientifiques where the current paper was completed. He
is grateful to them for their hospitality and financial support and for
providing a stimulating atmosphere.

\section{Groups, representations, equivariant K-theory, and 
  Lefschetz-Riemann-Roch
  }\label{sec:groups}

\subsection{Groups and Representations}
In this section, we introduce some notation and terminology associated to
groups and their representations.

Throughout this paper $G$ will denote a finite group.  The number of elements in a set $A$ will be denoted by $|A|$. The \emph{order} $|m|$ of an element $m$ in $G$
is defined to be the order $|\gen{m}|$ of the subgroup $\gen{m}$ generated by $m$. For every $m$ in $G$, the \emph{centralizer
$Z_G(m)$ 
of $m$ in $G$} is the subgroup of elements in $G$ which
commute with $m$. $Z_G(m)$ contains the cyclic subgroup $\gen{m}$ as a normal
subgroup. The \emph{set of conjugacy classes of $G$} is denoted by $\Gb$ and
for any element $m\in G$ we let $\mb \in \Gb$ denote the conjugacy class
containing $m$. For all $m$ in $G$, we have the useful identity 
\begin{equation} 
|\mb|
= |G|/|Z_G(m)|.
\end{equation}
Clearly, $|m|$ and $|Z_G(m)|$ are independent of the choice of representative
$m$ in $\mb$.
 
Let $\Rep(G;\nz)$ denote the \emph{(virtual) representation ring of $G$},
i.e., the Grothendieck group of finite-dimensional complex representations of
$G$. Since every representation of $G$ is 
uniquely decomposable into the direct sum of irreducible representations,
$\Rep(G;\nz)$ is 
a free $\nz$-module with a basis given by the set $\Irrep(G)$ of irreducible
representations of $G$. We denote the trivial irreducuble representation of
$G$ by $\bone$. 
The ring $\Rep(G;\nz)$ also has a \emph{metric} (a nondegenerate, symmetric
pairing) $\eta$, where $\eta(V,W) := \dim_\nc \textrm{Hom}_G(V,W)$. Furthermore, $\Irrep(G)$
is an orthonormal basis with respect to $\eta$. The \emph{regular
representation} $\nc[G]$ of $G$  is the group ring acted upon by left
multiplication and it satisfies the useful identity
\begin{equation}\label{eq:GroupRing}
\nc[G] = \sum_\alpha V_\alpha^{\oplus \dim V_\alpha},
\end{equation}
where the sum runs over all 
$V_\alpha \in \Irrep(G)$. 
In addition, $\Rep(G;\nz)$ has an involution called \emph{dualization},  corresponding to taking any representation
$W$ of $G$ and replacing it by its dual $W^*$. Dualization 
preserves the product and metric.  Since dualization also preserves
the set $\Irrep(G)$, Equation (\ref{eq:GroupRing}) implies that $\nc[G]^* =
\nc[G]$.  We will later need to work with $\nq$ coefficients, so we let
$\Rep(G) := \Rep(G;\nz)\otimes_\nz\nq$ and extend the product and metric
$\nq$-linearly. 

A representation of $G$ is determined, up to isomorphism, by its character,
i.e., if $\phi:G\rTo\Aut(V)$ is a $G$-module, then its character $\chi(V)$ is
the function defined by the trace $\chi_\ga(V) := \tr(\phi(\ga))$ for all
$\ga$ in $G$. The character $\chi(V)$ is a class function, i.e., it is a function on $\Gb$.
Evaluating at the identity element, we obtain $\chi_\one(V) = \dim
V$. In particular,
the character of the regular representation is 
\begin{equation}\label{eq:GroupRingChar}
\chi_{\ga}(\nc[G]) = |G| \delta_\ga^1
\end{equation}
for all $\ga$ in $G$, since left multiplication by $\ga$ acts without fixed
points unless $\ga$ is the identity. The character $\chi$ can be linearly
extended to $\Rep(G)$ to obtain a (virtual) character and every element in
$\Rep(G)$ is determined by its character. 

If $\nu:L\rTo G$ is a group homomorphism, then there is an associated ring
homomorphism $\nu^*:\Rep(G)\rTo \Rep(L)$ given by pulling back by $\nu$. 
An important special case is when $L$ is a subgroup of $G$ and $\nu$ is an
inclusion. The pullback of a $G$-module $V$ via the inclusion map yields an
$L$-module called the \emph{restriction} to $L$ and is denoted by $\Res^L_G V$. We
denote the induced ring homomorphism by $\Res^L_G :\Rep(G)\rTo
\Rep(L)$. Conversely, an $L$-module $W$ yields a $G$-module $\Ind^G_L W$,
called the \emph{induced module}, which is the tensor product of
$L$-modules $\nc[G]\otimes_{\nc[L]} W,$ where $L$ acts on $\nc[G]$ by right
multiplication. Induction yields a linear map
$\Ind^G_L:\Rep(L)\rTo\Rep(G)$. Restriction and induction are adjoint via
Frobenius reciprocity. 
Restriction obviously commutes with dualization; therefore, by Frobenius
reciprocity, induction commutes with dualization, as well.

An important special case for our purposes is the cyclic group $\gen{m}$
generated by an element $m$ of order $r := |m|$. Let $\V_m$ denote the
irreducible representation of $\gen{m}$, where $m$ acts as multiplication
by 
$\zeta_r :=\exp(- 2\pi i/r).$
For all integers $k$, let
\begin{equation}
\V_m^k := 
\begin{cases}
 \V_m^{\otimes k} & \text{if $k \geq 1$} \\
 \bone & \text{if $k = 0$} \\
(\V_m^{\otimes -k})^* & 
\text{if $k \leq -1$}.
\end{cases}
\end{equation}
In other words, $\V_m^k$ is the irreducible representation of $\gen{m}$ where
$m$ acts by multiplication by $\zeta_r^k$. Thus, we have the identity
\begin{equation}
\Irrep(\gen{m}) = \{\,\V_m^k \,\}_{k=0}^{r-1}
\end{equation}
and the isomorphism of algebras
\begin{equation}
\Rep(\gen{m}) \rTo^{\phi_m} \nq[\V_m] / \left< \V_m^r - \bone \right>, 
\end{equation}
where $\nq[\V_m]$ is the polynomial ring on the variable $\V_m$. Furthermore,
Equation (\ref{eq:GroupRing}) in this case reduces to
\begin{equation}\label{eq:GroupRingCyclic}
\nc[\gen{m}] = \sum_{k=0}^{r-1} \V_m^k.
\end{equation}

The canonical involution of groups $\sigb:\gen{m}\rTo\gen{m}$, taking a group
element to its inverse, induces a ring automorphism
$\sigb^*:\Rep(\gen{m})\rTo\Rep(\gen{m})$. It takes
\begin{equation}
\sigb^*\V_m^k := \V_{m^{-1}}^k = 
\begin{cases}
 \bone & \text{if $k=0$} \\
 \V_m^{r-k} & \text{if $ 1\leq k\leq r-1$}.
\end{cases}
\end{equation}
Therefore, we have
\[
\sigb^*\circ\phi_m = \phi_{m^{-1}}.
\]

A useful identity intertwining the cyclic group and conjugation is 
\begin{equation}\label{eq:IndCyclicConjInv}
\Ind^G_{\gen{m}} \V_m^k = \Ind^G_{\gen{\ga m\ga^{-1}}} \V_{\ga m\ga^{-1}}^k
\end{equation}
for all $\ga,m$ in $G$ and $k=0,\ldots,|m|-1$.

There are several useful dimensions associated to elements of $\Rep(G)$. For
all $W$ in $\Rep(G)$, the dimension $d(W)$ is
\begin{equation}\label{eq:etaCG}
d(W) := \chi_1(W) = \eta(W,\nc[G]) ,
\end{equation}
where  the second equality follows from  Equation (\ref{eq:GroupRing}).
Given, in addition, a group
element $m$ in $G$, we define $d_{m k}(W)$ in $\nq$ via
\begin{equation}\label{eq:dmk}
\Res^G_{\gen{m}} W =: \sum_{k=0}^{|m|-1} d_{m k}(W) \V_m^k.
\end{equation}
Frobenius reciprocity implies that for all $k=0,\ldots,|m|-1$,
\begin{equation}\label{eq:dFrobRec}
d_{m k}(W) = \widetilde{\eta}(\V_m^k,\Res^G_{\gen{m}} W) 
= \eta(\Ind^G_{\gen{m}}\V_m^k,W),
\end{equation}
where $\widetilde{\eta}$ is the metric on $\Rep(\gen{m})$.

Finally, we will need some maps associated to taking $G$-invariants. There is a $\nq$-linear map $\Rep(G)\rTo\Rep(G)$
which  associates to a $G$-module $W$, the submodule $W^G$  of $G$-invariant
vectors. This map commutes with induction, i.e., 
\begin{equation}\label{eq:InvOfInd}
(\Ind^G_L W)^G = W^L,
\end{equation}
and with dualization.
The regular representation is respected by induction, i.e., for any subgroup
$L$ of $G$,
\begin{equation}\label{eq:IndRegRep}
\Ind^G_L \nc[L] = \nc[G].
\end{equation}
This leads to a useful $\nq$-linear map $\I^G:\Rep(G)\rTo\Rep(G)$; namely, 
\begin{equation}\label{eq:IsupG}
\I^G(W) := W - \nc[G]\otimes W^G.
\end{equation}
Note that for any $W\in
\Rep(G)$ the $G$-invariants of $\I^G (W)$ vanish, as does $\I^G(\nc[G])$: 
\begin{equation}\label{eq:Ivanish}
\left(\I^G(W)\right)^G = 0 = \I^G(\nc[G]),
\end{equation}
and for any $\ga \in G$, its trace  is given by 
\begin{equation}
\chi_\ga (\I^G(W)) = \begin{cases} 
\chi_\ga(W) & \text{ if $\ga\neq 1$}\\
\chi_1(W) - |G|\otimes W^G  & \text{ if $\ga = 1$}.
\end{cases}
\end{equation}
Equations (\ref{eq:InvOfInd}) and (\ref{eq:IndRegRep})  imply that $\I$ also
commutes with induction, i.e.,   for all $W$ in $\Rep(L)$, 
\begin{equation}\label{eq:IandInduction}
\I^G(\Ind^G_L W) = \Ind^G_L (\I^L W).
\end{equation}
Finally, $\I$ also commutes with dualization, since $\nc[G]$ is self-dual
and taking $G$-invariants commutes with dualization.

\subsection{Equivariant K-theory and cohomology}

Let $\cX$ be a variety or orbifold (DM-stack) with the action of a finite group
$G$. Let $K(\cX)$ and $K_G(\cX)$ denote the Grothendieck group (with $\nq$ coefficients) of coherent sheaves and $G$-equivariant 
coherent sheaves on $\cX$, respectively. 

Suppose that $G$ acts trivially on $\cX$. In this
case, there is an algebra isomorphism 
\begin{equation}
\Phi:K_G(\cX)\rTo K(\cX)\otimes\Rep(G),
\end{equation}
 since the $G$-action on an
equivariant vector bundle acts linearly on each fiber. 
\begin{df}
Let $\chi_{\gamma}:\Rep(G) \rTo \nc$ denote the homomorphism which sends each
representation to its character at $\gamma$.  We will abuse notation and also
write $\chi_{\gamma}$ to denote the map $$\id\otimes\chi_{\gamma}:K(\cX)\otimes
\Rep(G) \rTo K(\cX)\otimes \nc .$$ 
\end{df}
\begin{rem}
The composition $\chi_{\gamma} \circ \Phi:K_G(\cX) \rTo K(\cX)\otimes\nc$ plays an important role in the Lefschetz-Riemann-Roch Theorem, which is a fundamental tool in this paper.

Note that $\chi_{\gamma} \circ \Phi $ takes a $G$-equivariant bundle $E$ to the
eigenbundle decomposition $$E \rMapsto^{\chi_{\gamma} \circ \Phi  }\sum_{\zeta}
E_{\gamma,\zeta} \otimes \zeta,$$  where the sum runs over all eigenvalues
$\zeta$ of $\gamma$ and $E_{\gamma,\zeta}$ is the eigenbundle of $E$ where
$\gamma$ acts with eigenvalue $\zeta$.

In particular, $\chi_{1} \circ \Phi:K_G(\cX) \rTo K(\cX)
$ is just the map which
forgets the $G$-equivariant structure of the elements of $K_G(\cX)$. 
\end{rem}

\begin{df}
If $\gamma$  is an element of a finite group $G$, and if $\cX$ 
has  a $G$-action, 
then we can define a morphism
$$\ell_{\gamma}:=\chi_{\gamma} \circ \Phi \circ {\operatorname{res}}:K_{\gen{\gamma}}(\cX)\rTo K(\cX^{\gamma})\otimes \nc$$
as the composition of $\chi_{\gamma} \circ \Phi $ with the obvious restriction $\operatorname{res}:K_{\gen{\gamma}}(\cX) \rTo K_{\gen{\gamma}}(\cX^{\gamma})$.

We can also compose $\ell_{\gamma}$ with the morphism that forgets the
$G$-equivariant structure and only retains the $\gen{\gamma}$-equivariant
structure, that is, $\ell_\ga\circ \Res^G_{\gen{\ga}}$.  When there is no
danger of confusion, we  denote this composition by $\ell_{\gamma}$ as well: 
$$\ell_{\gamma} : K_G(\cX) \rTo K(\cX^{\gamma}) \otimes \nc.$$
\end{df}
\begin{df} 
Let $\Ch:K(\cX)\rTo A^\bullet(\cX)$ denote the \emph{Chern character}
homomorphism.
Let $\bCh:K_G(\cX)\rTo A^\bullet(\cX)\otimes\Rep(G)$ denote the
\emph{$\Rep(G)$-valued Chern character} (algebra) 
homomorphism
$(\Ch\otimes\id_{\Rep(G)})\circ\Phi$. 
For all $j\geq 0$, let $\Ch_j:K(\cX)\rTo A^j(\cX)$ denote the projection of $\Ch$
onto $A^j(\cX)$ and, similarly, for $\bCh_j:K_G(\cX)\rTo A^j(\cX)\otimes\Rep(G)$.
\end{df}

\begin{rem}The map $\Phi$ commutes with push-forwards, in the sense that
$\Phi\circ f_* = (f_* \otimes \id_{\Rep(G)})\circ \Phi$ for any morphism
$f:\cX\rTo \cY$ between 
orbifolds with a trivial $G$ action.

Given $\cg\in K_G(\cX)$, if we
write $\Phi(\cg) = \sum_{\beta\in \Irrep(G)}\cg_\beta \otimes \beta$, then  we
see that for any \emph{representable} morphism $f:\cX\rTo\cY$
the $G$-equivariant form of the Grothendieck-Riemann-Roch Theorem is  
$$\bCh(f_*\cg) = \sum_\beta f_*(\Ch(\cg_\beta)\Td(T_f)) \otimes \beta =  (f_*
\otimes \id_{\Rep(G)}) \left(\bCh(\cg)(\Td(T_f)\otimes \vac)\right).$$ 
Since $T_f$ is $G$-invariant, we have $\Phi(T_f) = T_f \otimes \vac$, so if we
write $f_\star:=(f_* \otimes \id_{\Rep(G)})$ and 
$\td(T_f) :=(\Td \otimes \id_{\Rep(G)})\circ \Phi(T_f)$, 
then the $G$-equivariant Grothendieck-Riemann-Roch
Theorem can be written as
\begin{equation}
\bCh(f_*\cg)= f_\star (\bCh(\cg)\td(T_f)).
\end{equation}
Of course, if $f$ is not representable, there will be additional correction terms to the $GRR$ formula, as described, for example, in \cite{To}.
\end{rem}

\begin{df}
Let $\cX$ be connected. The \emph{$\Rep(G)$-valued rank}
$\brk:K_G(\cX)\rTo\Rep(G)$ is 
equal to $\bCh_0$, 
the projection of $\bCh$ onto
$A^0(\cX)\otimes\Rep(G)=\Rep(G)$. The \emph{virtual rank} is $\chi_1$ of
$\bCh$ and is denoted by $\rk := (\id_K \otimes \chi_1)\circ 
\bCh_0:K_G(\cX)\rTo
\nq$. Of course, if $\cX$ is not connected, the virtual rank of any $\cg\in K_G(\cX)$ is a locally constant 
$\nq$-valued function on $\cX$.
\end{df}

The Lefschetz-Riemann-Roch is a fundamental tool for this paper.   We need the theorem in the case of a family of semi-stable curves $\cx\rTo^f \cz$ over a base $\cz$ which is a Deligne-Mumford \emph{quotient stack}, that is, a DM stack $\cz$ for which there exists an algebraic space $Y$ and a linear algebraic group scheme $\LG$ such that $\cz=[Y/\LG]$.  

In order to make sense of the theorem in this case we will need to invert certain elements of the K-theory of $\cz$.  Specifically, we need to invert the element ${\chi_{\gamma} \circ \Phi \left(\lambda_{-1}(\conor_{\cx^{\gamma}/\cx})\right)}$  in $K(\cx^\ga)\otimes \nc$.  If $\cz$ (and hence $\cx$) were a projective variety, then 
${\chi_{\gamma} \circ \Phi \left(\lambda_{-1}(\conor_{\cx^{\gamma}/\cx})\right)}$ would already be invertible (see, for example, \cite[VI \S9]{FuLa}), but when $\cz$ is a more general DM stack, that may no longer be true, so we will need to localize, as follows.  
\begin{df} \label{df:localK}
Writing $\cz =[Y/\LG]$ as a global quotient, we have $K(\cz)= K_{\LG}(Y)$. 
Let $\loc$ denote the augmentation ideal of the representation ring
ring $\Rep(\LG)$ of $\LG$.   
\end{df}
\begin{prop}\label{prop:KisA}
If $\cz=[Y/\LG]$ is a DM quotient stack of finite dimension,
the localization $K(\cz)_\loc$ is isomorphic to the rational Chow ring $A^{\bullet}(\cz)$.

\end{prop}
\begin{proof}
By \cite[Thms 4.1, 5.1]{EG} the localization $K(\cz)_\loc$ of the ring
$K(\cz) = K_{\LG}(Y)$ at the ideal $\loc$ is isomorphic to the
completion $\widehat{K(\cz)}$ of ${K(\cz)}$ along the
ideal $\loc$.  Edidin and Graham \cite{EG} prove that $\widehat{K(\cz)}$ is isomorphic to  $\prod_{i=0}^{\infty}A_{\LG}^i(Y) = \prod_{i=0}^{\infty}A^i(\cz)$.

Moreover, since $\cz$ is a DM stack of finite dimension, the Chow groups $A^{i}(\cz)$ vanish for $i>>0$.  Therefore we have $$K(\cz)_\loc \cong \widehat{K(\cz)} \cong A^{\bullet}(\cz).$$ 

\end{proof}

\begin{thm}[Lefschetz-Riemann-Roch (LRR) for quotient stacks]\label{thm:DM-LRR}
Let $\cx\rTo^f \cz$ be a flat, projective family of semi-stable curves over a smooth DM quotient stack $\cz = [Y/\LG]$ with quasi-projective coarse moduli space.  Assume that a finite cyclic group $\gen{\gamma}$ acts on $\cx$ and acts trivially on $\cz$ such that $f$ is $\gen{\gamma}$-equivariant.  Denote by $f^{\gamma}$ the restriction of $f$ to the fixed-point locus $\cx^{\gamma}$.

Let $\conor_{\cx^{\gamma}/\cx}$ be the conormal sheaf of $\cx^{\gamma}$ in $\cx$, and let
$L_\cx: K_{\gen{\gamma}}(\cx) \rTo K(\cx^{\gamma})$ 
be defined as 
\begin{equation}\label{eq:L}
L_\cx (\cf)  := \frac{\ell_{\gamma}(\cf)}{\chi_{\gamma} \circ \Phi \left(\lambda_{-1}(\conor_{\cx^{\gamma}/\cx})\right)}.
\end{equation}
  Denote the K-theoretic push-forward along $f$ by $Rf_*(\cf) := \sum_{i=0}^\infty (-1)^i R^if_* \cf$, and let $\loc$ denote the augmentation ideal of the representation ring $\Rep(\LG)$ of $\LG$.   

The following diagram commutes:
\begin{equation}\label{eq:LRR}
\begin{diagram}
K_{\gen{\gamma}}(\cx)_{\loc} & \rTo^{L_\cx} & K(\cx^{\gamma})_{\loc} \otimes \nc\\
\dTo^{Rf_*}	& & \dTo^{Rf^{\gamma}_*}\\
K_{\gen{\gamma}}(\cz)_{\loc} & \rTo^{\ell_{\gamma}} & K(\cz)_{\loc}\otimes \nc 
\end{diagram}.
\end{equation}
\end{thm}
For a smooth scheme, this theorem is proved \cite[VI\S9]{FuLa}, in the singular  case in \cite{BFQ,Quart}, and in the non-projective case in \cite{Ta}.  We know of no published proof of the theorem that would apply to the case we need, so we give the proof here.
\begin{proof}

Since $\cz$ is a smooth DM quotient stack with quasi-projective coarse moduli, we can use the theorem of Kresch and Vistoli \cite[Thm 1]{KV}, which shows there is a smooth, quasi-projective scheme $Z$, and a finite, flat, surjective, local complete intersection (lci) morphism $g:Z\rTo \cz$.  We will endow $Z$ with the trivial $\gen{\ga}$-action.  The fiber product $X:=\cx \times_\cz Z$ is a semi-stable curve over $Z$, and it naturally inherits a $\gen{\ga}$-action from $\cx$. Denote the first and second projections from the fiber product by $g':X\rTo \cx$ and $f':X\rTo Z$, respectively.

Because the morphism $g$ is finite and surjective, the pullback $g^*:A^{\bullet}(\cz) \rTo A^{\bullet}(Z)$ is injective.  By Proposition~\ref{prop:KisA} the pullback is also injective on localized K-theory:
$$ g^*:K(\cz)_\loc \rInto K(Z)_\loc 
$$
Let $L_X:K_{\gen{\ga}}(X)_{\loc} \rTo K(X^{\ga})_{\loc} \otimes \nc$ be the homomorphism $L_X (\cf)  := \frac{\ell_{\gamma}(\cf)}{\chi_{\gamma} \circ \Phi \left(\lambda_{-1}(\conor_{X^{\gamma}/X})\right)},$
where $\conor_{X^{\gamma}/X}$ is the conormal bundle of $X^\ga$ in $X$.

By the usual Lefschetz-Riemann-Roch theorem for quasi-projective schemes \cite{BFQ}, the front face of the following diagram (Diagram~(\ref{diag:LRRproof})) commutes.  And since $g$ is flat, the two sides commute.  It is straightforward to check  that the bottom face commutes.

Note that $X^\ga$ is the fiber product $X^\ga = X\times_\cx \cx^{\ga}$, and the inclusions $\cx^\ga\rInto \cx$ and $X^\ga \rInto X$ are regular embeddings.  Therefore the excess conormal sheaf $E:= \left({g'}^{\ga *}\conor_{\cx^{\gamma}/\cx}\right)/\conor_{X^{\gamma}/X}$ is locally free \cite[\S6.3]{Fu}.  Moreover,  $g'$ is finite and flat, so the rank of $E$ is zero---that is, 
$${g'}^{\ga *}\conor_{\cx^{\gamma}/\cx} = \conor_{X^{\gamma}/X}.$$  This shows that the top face of the diagram commutes.

\begin{equation}\label{diag:LRRproof}
\begin{diagram}
&&K_{\gen{\ga}}(\cx)_{\loc} && \rTo^{L_\cx}&  & & K(\cx^{\ga})_{\loc} \otimes \nc \\
&\ldTo^{{g'}^*}&\vLine&&&&\ldTo^{{g'}^{\ga *}}\\
K_{\gen{\ga}}(X)_{\loc}&&\HonV&\rTo^{L_X}&&K(X^{\ga})_{\loc} \otimes \nc \\
\dTo^{Rf'_*}&&\dTo^{Rf_*}&&&\dTo^{R{f'}^{\ga}_*}&&\dTo^{{Rf}^{\ga}_*} \\
&&K_{\gen{\ga}}(\cz)_{\loc}&&\hLine^{\ell_\ga}&\VonH&\rTo&K(\cz)_{\loc}\otimes \nc\\
&\ldTo^{g^*}&&&&&\ldInto^{g^*}&&&&\\
 K_{\gen{\ga}}(Z)_{\loc}&&&\rTo^{\ell_\ga}&&K(Z)_{\loc}\otimes \nc 
\end{diagram}
\end{equation}

Since all the faces commute, except possibly the back, then for any $\cf\in K_{\gen{\ga}}(\cx)_{\loc}$, we have $$g^*Rf^{\ga}_*(L_\cx(\cf)) = g^* \ell_\ga R{f}_* \cf.$$
Since $g^*$ is injective, this gives 
$$Rf^{\ga}_*(L_\cx(\cf)) = \ell_\ga R{f}_* \cf,$$
as desired.\end{proof}

\section{Tautological bundles on the moduli space of (admissible) $G$-covers} 

\subsection{Moduli space of admissible $G$-covers}
We recall properties of the moduli spaces $\M_{g,n}(\BG)$ and $\MM_{g,n}$
that we will need throughout this paper.

For any finite group $G$, let $\M_{g,n}(\BG)$ denote the moduli space of
admissible $G$-covers and let $\MM_{g,n}$ denote the moduli space of pointed
admissible $G$-covers. We adopt the notation and definitions from \cite{JKK}.
We will denote by $\M$, $\M(\BG)$, and $\MM$ the 
disjoint
union 
over
all $g$ and $n$ of $\M_{g,n}$, $\M_{g,n}(\BG)$, and
$\MM_{g,n}$, respectively.
Note that our definition of $\M_{g,n}(\BG)$ differs slightly from that of
Abramovich-Corti-Vistoli \cite{ACV} in that the source curves possess
honest sections $\sigma_i$ instead of just gerbe markings.  

The map $\MM_{g,n}\rTo \M_{g,n}(\BG)$ is representable. And it is known that  $\M_{g,n}(\BG)$ is a quotient stack \cite{AGOT}.  Thus we have the following proposition.
\begin{prop}\label{thm:MMQuot}
The stack $\MM_{g,n}$ is a  quotient stack.
\end{prop} 

In particular, Lefschetz-Riemann-Roch (Theorem~\ref{thm:DM-LRR}) holds for the universal $G$-cover $\ce \rTo \MM_{g,n}$.

\subsubsection{Graphs}

Let $\bGa_{g,n}$ denote the set of (connected)
stable graphs of genus $g$ with $n$ tails and \emph{exactly one edge}. Elements in
$\bGa_{g,n}$ are either trees or loops. If $\Ga$ belongs to $\bGa_{g,n}$, then we let $\bGac$ denote the set of (possibly disconnected) stable graphs $\Gac$ obtained by cutting $\Ga$
along its edge together \emph{with a choice of ordering of the two new
edges}.  We denote the $n+1$st edge by $+$ and the $n+2$nd edge by $-$. For
example, if $\Ga$ is a tree and $n > 0$ or the genera of the two vertices are
different, then $\bGac$ contains two disconnected stable graphs. If $\Ga$ is a loop, then there is a single graph $\Gac$ in $\bGac$. \emph{Notice that the gluing morphisms on $\M$ are naturally
indexed by $\Gac$, rather than $\Ga$, since an ordering on the half edges must be specified.}

To describe operations
on $\MM_{g,n}(\bm)$, where $\bm$ belongs to $G^n$, we must decorate the graphs.   Let $\bGat_{g,n}(\bm)$ denote
the set of (connected) stable graphs of genus $g$ with exactly one edge and
with $n$ tails, such that the $i$-th tail is decorated by  $m_i$ for all
$i=1,\ldots,n$.   If $\Gat$ is a such decorated graph, then let
$|\Gat| \in \bGa_{g,n}$ indicate the same graph but without decorations.

Choose $\mb$ in
$\Gb$, and define $\bGat_{g,n}(\bm,\mb)$ to be the set of equivalence classes of the
set of (connected) stable graphs of genus $g$ with $n$ tails whose $i$-th
tail is decorated with $m_i$ for all $i=1,\ldots,n$, one half edge is
decorated with some group element $m_+=m$ in $\mb$, and the other half edge is
decorated by the group element $m_-=m^{-1}$. 

Two such graphs $\Gat, \Gat'\in \bGat_{g,n}(\bm,\mb)$ will be considered 
equivalent if there exists an isomorphism $\alpha$ of the underlying stable graphs 
such that the induced decorated graph $\alpha\Gat'$ differs from $\Gat$  only in
the group elements $(m_+, m_-)$ or $(m'_+, m'_-)$ associated to its two half
edges and if the pairs $(m_+, m_-)$ and $(m'_+,
m'_-)$  satisfy $(\ga \mpp\ga^{-1},\ga\mmm\ga^{-1}) = (\mpp',\mmm')$ for some $\ga$
in $G$.
We have
\[
\bGat_{g,n}(\bm) := \coprod_{\mb\in\Gb} \bGat_{g,n}(\bm,\mb).
\]

If $\Gat$ belongs to $\bGa_{g,n}(\bm,\mb)$, then 
$\bGatc(\Gat)$  denotes the set of decorated graphs $\Gatc$  obtained by cutting the edge of $\Gat$ following by
\emph{choosing an ordering for the newly created tails (as before) and a
group element $\mpp$ in $\mb$ associated to the ``$+$" tail and the group element
$\mmm := \mpp^{-1}$ associated to the ``$-$" tail.}  
Similarly,  
denote 
$$\bGatcgn:=\coprod_{\Gat \in \bGat_{g,n}} \bGatc(\Gat) \dsand\bGatcgn(\bm):=\coprod_{\Gat \in \bGat_{g,n}(\bm)} \bGatc(\Gat).$$
And let $\bGatcgn(\bm,m_+,m_-)\subseteq \bGatcgn$ denote the set of decorated cut graphs whose $i$th tail is decorated with $m_i$ and whose $+$ and $-$ tails are decorated with $m_+$ and $m_-$, respectively.

For any $\Ga$ in $\bGa_{g,n}$ we denote the closure in $\MM_{g,n}$ of the
substack of pointed admissible covers with dual graph $\Ga$ by $\MM_\Ga$.  And similarly, for any $\Gat$ in $\bGat_{g,n}(\bm,\mb)$ we denote the closure in $\MM_{g,n}$ of the substack of pointed admissible covers with decorated dual graph $\Gat$ by $\MM_\Gat$. We have
\[
\MM_\Ga = \coprod_{\bm \in G^n}\bigcup_{\mb\in\Gb}
\bigcup_{\substack{\Gat\in\bGat_{g,n}(\bm,\mb)\\ |\Gat| = \Ga}} \MM_\Gat.
\]

\subsubsection{Notation}

Throughout this paper we will use the following notation and refer to the following diagrams.

To reduce clutter, and when there is little chance of confusion, we will
adopt the following convention: when a map (e.g., $\pibar$) or an object
(e.g., $\cc$) is the pullback, in a Cartesian square, of a universal map
or object, we will use the same name for both maps---thus $\pibar$ is the
label for several parallel maps  in the following diagrams.   

We let $\pibar:\cc\rTo\M_{g,n}$ denote the universal curve over the stack $\M_{g,n}$, and let $s:\M_{g,n}(\BG) \rTo \M_{g,n}$ denote the morphism induced by forgetting the $G$-cover.  Similarly, let $t:\MM_{g,n}\rTo\M_{g,n}(\BG)$ denote the morphism induced by forgetting the sections of the universal $G$-cover.  

Let $\pih:\ce\rTo\cc$ denote the universal $G$-cover of $\pibar:\cc\rTo\M_{g,n}(\BG)$. For any $\Gat\in \bGat_{g,n}$ with $|\Gat| = \Ga$, we have corresponding inclusions $i_\Ga:\M_\Ga\rTo\M_{g,n}$ and $i_\Gatc:\M_\Gatc(\BG) \rTo \M_{g,n}(\BG)$ and $i_\Gatc:\MM_\Gatc\rTo\MM_{g,n}$, and inclusions $\ih_\Gatc:\cc_\Gatc \rTo \cc$ of the universal curves and inclusions $\itt_{\Gatc}:\ce_\Gatc \rTo \ce$ and inclusions of the universal $G$-covers.

All these morphisms, stacks, and universal objects fit together in the following diagram.

\begin{equation}\label{diag:MM-MBG-M}
\begin{diagram}
&&\ce &&&\rTo^\tt & \ce \\
&\ruTo^{\itt_\Gat}&\vLine&&&\ruTo^{\itt_\Gat}\\
\ce_\Gat&&\HonV&\rTo^{\tt_\Gat}&\ce_\Gat\\
\dTo^{\pih_\Gat}&&\dTo_{\pih}&&\dTo^{\pih_\Gat}&&\dTo_{\pih}\\
&&\cc&\hLine&\VonH&\rTo^\th&\cc&&&\rTo^\sh&\cc \\
&\ruTo^{\ih_\Gat}&\vLine&&&\ruTo^{\ih_\Gat}&\vLine&&&\ruTo^{\ih_\Ga}\\
\cc_\Gat&&\HonV&\rTo^{\th_{\Gat}}&\cc_\Gat &&\HonV&\rTo^{\sh_\Gat}&\cc_\Ga\\
\dTo^{\pibar_\Ga}&&\dTo^{\pibar}&&\dTo^{\pibar_\Gat}&&\dTo^\pibar &&\dTo^{\pibar_{\Ga}}&&\dTo^\pibar\\
&&\MM_{g,n}&\hLine&\VonH&\rTo^t&\M_{g,n}(\BG)&\hLine&\VonH&\rTo^s&\M_{g,n} \\
&\ruTo^{i_{\Gat}}&&&&\ruTo^{i_\Gat}&&&&\ruTo^{i_\Ga}\\
\MM_{\Gat}&&\rTo^{t_\Gat}&&\M_{\Gat}(\BG)&&&\rTo^{s_\Gat}&\M_{\Ga}
\end{diagram}
\end{equation}

We define $st:=s\circ t$ and 
$\pi = \pibar\circ\pih,$ and  let $\sigma_i:\MM_{g,n} \rTo \ce$ be the $i$th section of the pointed admissible $G$-cover $\pi$. Let $\sigbar_i:=\pih \circ \sigma_i:\M_{g,n} \rTo \cc$ be the $i$th section of the universal curve $\pibar$.

We will also find it useful to decompose the gluing morphisms.  If $\Gatc\in \bGatcgn$ glues to give the graph $\Gat\in \bGat_{g,n}$ and has an underlying graph $|\Gatc| = \Gac \in \bGacgn$, then the corresponding gluing morphisms $\rho_\Gatc:\MM_\Gatc \rTo \MM_{g,n}$ and $\rho_\Gac:\M_\Gac\rTo \M_{g,n}$ decompose into a composition $\rho_\Gatc = i_\Gat \circ \mu_\Gatc$, where $\mu$ is the obvious map from the stack of curves or $G$-covers with cut graph to the stack of curves or $G$-covers with uncut graph.

These fit together in the following diagram, to which we will refer throughout the paper.
\begin{equation}
\begin{diagram}\label{diag:cutting}
&&\ce_{\Gatc}&\rTo^{\mut_\Gatc}&\ce_{\Gat}& \rTo^{\itt_\Gat}&\ce\\
&&\dTo^{\pi_\Gatc}&&\dTo^{\pi_\Gat}&&\dTo^\pi\\
&&\MM_\Gatc &\rTo^{\mu_\Gatc} & \MM_\Gat & \rTo^{i_\Gat}&\MM_{g,n}\\
&&\dTo^{\st_\Gatc}&&\dTo^{\st_\Gat}&&\dTo^\st\\
&&\M_\Gac &\rTo^{\mu_\Gac} & \M_\Ga & \rTo^{i_\Ga}&\M_{g,n},
\end{diagram}
\end{equation}
\begin{rem}
Note that while the top two squares in Diagram~(\ref{diag:cutting}) are both Cartesian, usually neither of the bottom two squares 
is Cartesian.  

Moreover, the morphisms $\mu_\Gac$ and $\mu_\Gatc$ are generally not \'etale.  The morphism $\mu_\Gac$ fails to be \'etale precisely on the locus where $\M_\Ga$ has a self-intersection, that is, on the locus $\M_{\Ga'}$ of two-edged graphs $\Ga'$ which have the property that contracting \emph{either} of the two edges of the graph $\Ga'$ gives the one-edged graph $\Ga$.

Similarly, the morphism $\mu_\Gat$ fails to be \'etale on the locus $\MM_{\Gat'}$ of decorated two-edged graphs $\Gat'$ which have the property that contracting either of the two edges of $\Gat'$ gives $\Gat$.
\end{rem}

We let 
$\tau:\MM_{g,n+1}(\bm,1)\rTo\MM_{g,n}(\bm)$ denote the forgetting tails map.
Associated to the this map we have the following diagram:
\begin{equation}\label{diag:tails}
\begin{diagram}
&&\ce_{g,n+1} &&&\rTo^\tt & \ce_{g,n+1} \\
&\ldTo^{\taut}&\vLine&&&\ldTo^{\taut}\\
\ce_{g,n}&&\HonV&\rTo^{\tt}&\ce_{g,n}\\
\dTo^{\pih}&&\dTo_{\pih}&&\dTo^{\pih}&&\dTo_{\pih}\\
&&\cc_{g,n+1}&\hLine&\VonH&\rTo^\th&\cc_{g,n+1}&&&\rTo^\sh&\cc_{g,n+1} \\
&\ldTo^{\tauh}&\vLine&&&\ldTo^{\tauh}&\vLine&&&\ldTo^{\tauh}\\
\cc_{g,n}&&\HonV&\rTo^{\th}&\cc_{g,n} &&\HonV&\rTo^{\sh}&\cc_{g,n}\\
\dTo^{\pibar}&&\dTo^{\pibar}&&\dTo^{\pibar}&&\dTo^\pibar &&\dTo^{\pibar}&&\dTo^\pibar\\
&&\MM_{g,n+1}(\bm,1)&\hLine&\VonH&\rTo^t&\M_{g,n+1}(\BG)(\bmb,1)&\hLine&\VonH&\rTo^s&\M_{g,n+1} \\
&\ldTo^{\tau}&&&&\ldTo^{\tau}&&&&\ldTo^{\tau}\\
\MM_{g,n}(\bm)&&\rTo^{t}&&\M_{g,n}(\BG)(\bmb)&&&\rTo^{s}&\M_{g,n}
\end{diagram}
\end{equation}

\subsubsection{Basic Properties of the Morphisms and Stacks}

\begin{prop}\label{prop:FundClasses}
For all $\Ga \in \bGa_{g,n}$ and for all 
$\alpha$ in $A^\bullet(\M_\Ga)$, we have
\begin{equation}\label{eq:RamifCorrection}
\sum_{\substack{|\Gat|=\Ga \\ \Gat \in \bGat_{g,n}}} r_\Gat \, i_{\Gat *} \st_\Gat^*\alpha =  \st^*i_{\Ga *}\alpha,
\end{equation}
where 
$r_\Gat  = |m|$ 
if $\Gat \in \bGat_{g,n}(\bm,\mb)$.

In particular, the fundamental classes are related by 
the equality 
\begin{equation}\label{eq:FundClasses}
\sum_{|\Gat| = \Ga} 
r_\Gat [\MM_\Gat] = \st^*[\M_\Ga]
\end{equation}
in $A^\bullet(\MM_\Gat)$.
In addition, for any $\Gatc \in \bGatcgn(\bm, \mpp,\mmm)$, let $\Gat$ be the
graph obtained by  gluing the cut edge, let $\Gac := |\Gatc|$, and let $\Ga
:=|\Gat|$. For all $\alpha$ in $A^\bullet(\M_\Gac)$ we have the equality 
\begin{equation}\label{eq:muPushPull}
\mu_{\Gatc *}\st_\Gatc^*\alpha = 
\frac{|\Aut(\Gat)||Z_G(\mpp)|}{|\Aut(\Ga)| 
|\mpp|}
\st_\Gat^* \mu_{\Gac *}\alpha.
\end{equation}
For any class $\beta \in A^\bullet(
\bigcup_{\Gac\in\bGacgn}
\M_\Gac)$ and for any $\bm \in G^n$, 
we have
\begin{equation}\label{eq:StRhoFullPushPull}
\frac12 \sum_{\Gac\in\bGacgn} st^*\rho_{\Gac *} \beta = 
\frac12 \sum_{\Gatc\in \bGatcgn(\bm)}\frac{r_\Gatc^2}{|Z_G(\mpp)|}\rho_{\Gatc *}st_{\Gatc}^* \beta.
\end{equation}
\end{prop}
\begin{proof}
Let $F$ be the fibered product 
$F_\Gat:= \M_{\Gamma} \times_{\M_{g,n}}\MM_{g,n}.$
It is straightforward to see that the stack 
$\bigcup_{|\Gat|=\Ga} 
\MM_{\Gat}$ is the 
reduced induced closed substack underlying $F$ (i.e., the result of annihilating nilpotents in the structure sheaf).
So the fibered product $F$ also breaks up into a 
union of pieces indexed by $\Gat$ 
$$F = \bigcup_{|\Gat|=\Ga} F_{\Gat},$$
and the reduced induced closed substack underlying $F_\Gat$ is $\MM_\Gat$.   We have the following commutative diagram: 
\begin{equation}
\begin{diagram}
\bigcup_{|\Gat|=\Ga} \MM_{\Gat}\\
& \rdTo~{\bigcup j_{\Gat}}\rdTo(2,4)_{\bigcup \st_{\Gat}}\rdTo(4,2)^{{\bigcup i_{\Gat}}}\\
&& \bigcup_{|\Gat|=\Ga} F_\Gat & \rTo_{pr_2} & \MM_{g,n}
\\
&&\dTo_{pr_1} && \dTo_{\st}\\
&&\M_{\Gamma} & \rTo_{i_{\Gamma}} & \M_{g,n}.
\end{diagram}
\end{equation}
We also have that $j_{\Gat *}: A^{\bullet}(\MM_{\Gat}) \rTo A^{\bullet}(F_{\Gat})$ is an isomorphism, and $pr_1^*\alpha = r_\Gat j_{\Gat *} st_\Gat^*\alpha$, where $r_\Gat$ is the 
degree of
ramification  of $st$ along $\Gat$. That is,
$r_\Gat$ is the number of non-isomorphic pointed admissible $G$-covers over a
generic point of $\MM_{g,n}$ which degenerate to the same isomorphism class
of pointed admissible $G$-covers over a generic point of $\MM_{\Gat}$.  The
degeneration to $\MM_\Gat$ comes from contracting a cycle in the underlying
curve (with holonomy $m$) to a single point.  

Thus, after accounting for
automorphisms of the smooth versus the nodal $G$-covers, the number of
pointed admissible $G$-bundles over the smooth curve that contract to give
the same nodal $G$-cover is $r_\Gat =|m|$.  And we have 
$$ \sum_{|\Gat|=\Ga}r_\Gat i_{\Gat *} \st_\Gat^* \alpha = \sum_{|\Gat|=\Ga}  pr_{2 *} pr_1^* \alpha ={\st^* i_{\Ga *}\alpha},$$
where the last equality follows from the fact that $\st$ is flat.
This proves that Equation~(\ref{eq:RamifCorrection}) holds, and Equation~(\ref{eq:FundClasses}) is a special case.

For any $\Gatc$ denote the canonical map $\MM_\Gatc \rTo F':=\M_\Gac\times_{\M_{\Ga}} \MM_{\Gat}$ by $q$.  
It is easy to see that $q$ is finite and surjective; 
indeed, the product $F'$ consists of all triples $(\cc_{cut},(\ce_{glued}\to \cc_{glued},\pt_1,\dots,\pt_n), \alpha)$, 
where $\cc_{cut}$ is a curve in $\M_{\Gac}$, 
and $(\ce_{glued}\to \cc_{glued},\pt_1,\dots,\pt_n)$ is a pointed $G$-cover in $\MM_{\Gat}$, 
and $\alpha$ is an isomorphism between $\cc_{glued}$ and the curve obtained
by gluing $\cc_{cut}$. 
Normalizing $\ce$ at the node gives a new $G$-cover $\ce_{cut}$, and any
choice of $\pt_+$ with monodromy $\mpp$ will give an element of $\MM_\Gatc$
which maps to the original triple. 

Moreover, the degree of $q$ is 
\begin{equation}
\deg(q)= \frac{|\Aut(\Gat)||Z_G(\mpp)|}{|\Aut(\Ga)| |\mpp|},
\end{equation}
as can
be seen from the fact that 
\[
\deg(\mu_\Gatc) = \frac{|\Aut(\Gat)| |Z_G(\mpp)|}{|\mpp|},
\]
and the
degree of the second projection $pr_2:F' 
\rTo \MM_{\Gat}$ is the same as the degree of $\mu_\Gac$, namely
$\deg(\mu_\Gac) = {|\Aut(\Ga)|}$.  We now have
$$
\mu_{\Gatc *}\st_\Gatc^*\alpha = pr_{2 *}q_*q^* pr_1^* \alpha = \deg(q) \pr_{2 *}pr_1^* \alpha = 
\frac{|\Aut(\Gat)||Z_G(\mpp)|}{|\Aut(\Ga)| 
|\mpp|}
\st_\Gat^* \mu_{\Gac *}\alpha,
$$
which gives us Equation~(\ref{eq:muPushPull}). 

Finally, to prove Equation~(\ref{eq:StRhoFullPushPull}) we observe that
\begin{eqnarray}
\frac12 \sum_{\Gac\in\bGacgn} \st^*\rho_{\Gac *} \beta 
&=& \frac12 \sum_{\Gac\in\bGacgn} \sum_{\substack{|\Gat|=\Ga \\ \Gat \in \bGat_{g,n}}} r_\Gat \, i_{\Gat *} \st_\Gat^* \mu_{\Gac *} \beta\\
&=& \frac12 \sum_{\Gatc\in\bGatcgn} \frac{r_\Gat|\Aut(\Ga)|}{|\Aut(\Gat)|}  \, i_{\Gat *} \st_\Gat^* \mu_{\Gac *} \beta\\
&=& \frac12 \sum_{\Gatc\in \bGatcgn(\bm)}\frac{r_\Gatc^2}{|Z_G(\mpp)|}\rho_{\Gatc *}st_{\Gatc}^* \beta,
\end{eqnarray}
where the first equality follows from 
Equation~(\ref{eq:RamifCorrection}), the second from counting the number of graphs $\Gatc$ with $|\Gatc|=\Gac$, and the third from Equation~(\ref{eq:muPushPull}).
\end{proof}

\subsection{Tautological bundles and cohomology classes associated to the universal $G$-curve and their properties} 
\begin{df}
We define the bundle $\crrt$ on $\M_{g,n}(\BG)$ to be the 
push-forward  
\begin{equation}
\crrt := R^1\pi_*\co_\ce.
\end{equation}
Since the map $\MM_{g,n}\rTo^t \M_{g,n}(\BG)$ is flat, and the
universal admissible $G$-cover over $\MM_{g,n}$ is the pullback of the
admissible universal $G$-cover over $\M_{g,n}(\BG)$, the push-forward
$R^1\pi_*\co_\ce$ on $\MM_{g,n}$ is the pullback of $\crrt$ from
$\M_{g,n}(\BG)$.  We will abuse notation and also use $\crrt$ to denote this bundle on $\MM_{g,n}$. 
\end{df}

\begin{df}
Let $\clt_i$ be the line bundle given by pulling back the relative dualizing sheaf $\omega_{\pi}$ along the sections $\sigma_i$, and let $\psit_i$ be the first Chern class of $\clt_i$:
\begin{equation}
\clt_i := \sigma_i^*(\omega_\pi) \qquad
\psit_i := c_1(\clt_i).
\end{equation}
Similarly, on $\MM_\Gatc$, we have additional sections $\sigma_+$ and $\sigma_-$ and corresponding line bundles:
\begin{align}
\clt_+ &:= \sigma_+^* (\omega_\pi) \qquad &\psit_+ := c_1(\clt_+)\\
\clt_- &:= \sigma_-^* (\omega_\pi) \qquad &\psit_- := c_1(\clt_-).
\end{align}
\end{df}

\begin{prop}\label{prop:PsitPsi}
The bundles $\clt_i$ and $\cl_i$ in $\pic(\MM_{g,n})$ and the classes  $\psit_i := c_1(\clt_i)$ and $\psi_i := c_1(\cl_i)$ in
$A^1(\MM_{g,n}(\bm))$ 
are related by 
\begin{align}
\clt_i^{\otimes |m_i|} = \cl_i \qquad \text{ and }\qquad &
|m_i| \psit_i = \psi_i
\end{align}
for all $i=1,\ldots,n$.
Similarly, for all $\Gat$ in $\bGa(\bm)$, the classes $\psit_\pm := c_1(\clt_\pm)$ and $\psi_\pm
:= c_1(\cl_\pm)$ are related in 
$A^1(\MM_\Gatc)$ 
by  
\begin{equation}\label{eq:psit}
\psit_\pm = \frac{1}{r} \psi_\pm,
\end{equation}
where $r = |\mpp| = |\mmm|$.
\end{prop}

\begin{proof}
If $z$ is a local coordinate on $\ce$ near $\pt_i$, then $x:=z^{|m_i|}$ is a local coordinate on $\cc$ near $p_i$.  Locally on $\ce$ near $\pt_i$, the relative dualizing sheaf $\omega_\pi$ is generated by the one-form $dz$ and the pullback $\pih^*(\omega_{\pibar})$ is generated by $dx = |m_i| z^{|m_i|-1} dz$.  This shows that 
$$\pih^*(\omega_{\pibar}) = \omega_\pi \otimes \co(-|m_i|\Dt_i),$$ where $\Dt_i$ is the divisor in $\cc$ corresponding to the image of $\sigma_i$.
It is well-known (See, for example \cite[Lm 2.3]{JKV1}) that 
$$\sigma_i^*(\co(-\Dt_j)) = \begin{cases}\clt_i &\text{if $i=j$}\\
\co & \text{if $i\neq j$}
\end{cases}.
$$ 
Combining this with the fact that  
$\sigbar_i = \pih\circ\sigma_i$,
we have  
\begin{equation}
\cl_i = \sigma_i^*(\pih^*(\omega_{\pibar})) = \sigma_i^*\left(\omega_\pi \otimes \co\left(-(|m_i|-1)\Dt_i\right)\right) = \clt_i^{\otimes |m_i|}.  
\end{equation}
Taking first Chern classes completes the proof.
\end{proof}

We also need to define the analogue of Arbarello-Cornalba's kappa classes.
\begin{df}
For the universal curve $\pibar:\cc\rTo\MM_{g,n}$, let $$\omega_{\pibar,\mathrm{log}} := \omega_{\pibar}(\sum_{i=1}^n \bar{D}_i),$$
where $\bar{D}_i$ is the image of the $i$th section $\sigbar_i:\MM_{g,n} \rTo \cc$.

Similarly, for the universal $G$-cover $\pi:\ce \rTo\MM_{g,n}$, define 
\begin{equation}
\omega_{\pi,\mathrm{log}} := \omega_\pi (\sum_{i=1}^n \sum_{g \in G/\gen{m_i}} gD_i),
\end{equation}
where $D_i$ denotes the image of the section $\sigma_i$, and $\sum_{g\in
G/m_i} gD_i$ is the sum of all the translates of $D_i$. 
Now let
$$\kappa_a :=\pibar_*(c_1(\omega_{\pibar,\mathrm{log}})^{a+1}) \dsand 
\kat_a:=\pi_*(c_1(\omega_{\pi,\mathrm{log}})^{a+1}).$$
\end{df}

It is immediate to check that 
\begin{equation}
\omega_{\pi,\mathrm{log}}  = \pih^* \omega_{\pibar,\mathrm{log}},
\end{equation}
where $\pih:\ce\rTo \cc$ is the covering map.
Combining this with the fact that $\pih$ is finite of degree $|G|$, we have the following proposition. 
\begin{prop}
The classes $\kat_a$ and $\kappa_a$ on $\MM_{g,n}$ are related as follows:
\begin{equation}\label{eq:kappat}
\kat_a = |G|\kappa_a.
\end{equation}
\end{prop}
 
\begin{rem}
It is important to note that in addition to the definition of kappa classes given here, there is another common definition of kappa classes due to Mumford:
$$\kappa'_a:=\pibar_*(c_1(\omega_{\pibar})^{a+1}),
$$
and the obvious analogue for admissible $G$-covers would be
$$\kat'_a:=\pi_*(c_1(\omega_{\pi})^{a+1}).$$
Mumford's kappa classes don't behave as well as Arbarello-Cornalba's kappa classes, but the different definitions are related, as follows. 
\end{rem}
 \begin{prop}[{\cite[Eq~(1.5)]{ACComb}}]\label{prp:MumKappaAC}
\begin{equation}\label{eq:MumKappaAC}
\kappa_a = \kappa'_a + \sum_{i=1}^n\psi_i^{a}\dsand 
\kat_a= \kat'_a + 
\sum_{i=1}^n \frac{|G|}{|m_i|} \psit_i^{a}
\end{equation}
\end{prop}

\begin{figure}
\includegraphics{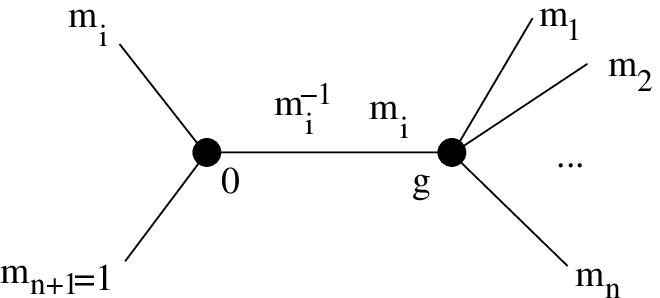}
\caption{\label{fig:Gati} The decorated graph $\Gat_i$ of Definition~\ref{df:Gati}.}
\end{figure}

Next we make a definition that will play an important role in the forgetting-tails morphism.
\begin{df}\label{df:Gati}
Let $\Gat_i$ be the one-edged tree with $n+1$ tails, and with one of its two vertices having only two tails, $p_i$ and $p_{n+1}$, decorated by $1$ and $m_i$, respectively, and the other vertex having $n-1$ tails $p_1,\dots,\hat{p}_i, \dots p_n$ decorated with $m_1, \dots, \hat{m}_i, \dots, m_n$, respectively (See Figure~\ref{fig:Gati}).
\end{df} 

\begin{thm} 
The following relations hold on $\MM$.  Those relations which do not involve
the gluing map $\rho$ also hold on $\M(\BG)$. 
\begin{enumerate}
\item The bundle $\crrt$ is preserved by the forgetting tails morphism $\tau$: 
\begin{equation}\label{eq:crrt-tails}
\tau^*\crrt = \crrt.
\end{equation}
The cotangent line bundles $\clt_i$ satisfy a relation like the Puncture Equation; namely, if $\Gat_i$ is the graph in Defintion~\ref{df:Gati}, then $\MM_{\Gat_i}$ is a divisor in $\MM_{g,n+1}$, and we have
\begin{equation}\label{eq:clt-tails}
\tau^*\clt_i = \clt_i\otimes \co(-\MM_{\Gat_i})
\end{equation}
in $\pic(\MM_{g,n+1}) \otimes \nq$ and
\begin{equation}\label{eq:psit-tails}
\tau^*\psit_i = \psit_i - \left[\MM_{\Gat_i}\right]
\end{equation}
in 
$A^1(\MM_{g,n+1},\nq)$.
\item Consider $\Gat$ in $\bGat_{g,n}(\bm,\mb)$ and choose  $\Gatc\in \bGatc(\Gat)$ of type $\mpp \in \mb$.
We have 
\begin{equation}\label{eq:rhoclt}
\rho_{\Gatc}^*\clt_i =  \clt_i,
\end{equation}
where the left-hand $\clt_i$ is the cotangent line bundle on $\MM_\Gat$ and the right-hand $\clt_i$ is the cotangent line bundle on $\MM_\Gatc$.

Furthermore, 
\begin{enumerate}
 \item If $\Gat$ is a loop, then 
 \begin{equation}\label{eq:rhocrrtloop}
 \begin{split}
 \rho_{\Gatc}^*\crrt = \crrt + \bigl(\nc[G/G_0] &-\nc[G/G_{cut}]\\
 & + \nc[G/\gen{m_ +}]\bigr) \otimes \co_{\MM_{\Gatc}},
\end{split}
 \end{equation}
 where $G_0$ is the subgroup of $G$ fixing the connected component of a
 generic 
admissible $G$-cover $ \ce_\Gat$ with dual graph $\Gat$ containing the node which is cut to give $\sigma_+$, and
 $G_{cut}$ is the subgroup of $G$ fixing the connected component of
 $\ce_{\Gatc}$ containing the image of $\sigma_+$. 
 \item If $\Gat$ is a tree, then 
 \begin{equation}\label{eq:rhocrrttree}
\begin{split}
 \rho_{\Gatc}^*\crrt = \crrt + \bigl(\nc[G/G_0] & -\nc[G/G_{+}] -\nc[G/G_{-}]\\
& +\nc[G/\gen{m_ +}]\bigr) \otimes \co_{\MM_{\Gatc}}, 
\end{split}
 \end{equation}
 where $G_0$ is the subgroup of $G$ fixing the connected component of a generic $ \ce_\Gat$ containing the node that is cut to give $\sigma_+$, $G_{+}$ is the subgroup fixing the connected component of $\ce_{\Gatc}$ containing the image of $\sigma_+$, and $G_{-}$ is the subgroup fixing the connected component of $\ce_{\Gatc}$ containing the image of $\sigma_-$.
\item For any integer $a\ge0$ we have 
\begin{equation}
\rho^*_\Gatc \kat_a = \kat_a,
\end{equation}
where we have used $\kat_a$ to denote both the class 
$\pi_*(c_1(\omega_{\pi,\mathrm{log}})^{a+1})$ in $A^a(\MM_{g,n}(\bm))$ \emph{and} the class 
$\pi_{\Gatc *}(c_1(\omega_{\pi_\Gatc,\mathrm{log}})^{a+1})$
in $A^a(\MM_\Gatc)$.
\end{enumerate}
\end{enumerate}
\end{thm}

\begin{proof}
To prove the first equation (\ref{eq:crrt-tails}), we consider the following diagram with a Cartesian square:
\begin{equation}\label{diag:tailssquare}
\begin{diagram}
&\ce_{{g,n+1}}\\
&\dTo^{\upsilon} & \rdTo^{\taut} \\
\ce':=&\taut^*(\ce)	& \rTo^\alpha 	& \ce_{g,n} \\
&\dTo^{\pi'_{g,n}}	& 		& \dTo^{\pi_{g,n}}\\
&\MM_{g,n+1}(\bm,1)			&\rTo^{\tau} & \MM_{g,n}(\bm) 
\end{diagram},
\end{equation}
where $\ce'$ is the pullback of $\ce$ along $\tau$, and $\pi_{g,n+1} = \pi'_{g,n} \circ \upsilon$.
The map $\upsilon$ contracts any components that are made unstable by the
removal of the marked point $p_{n+1}$ and the points in its $G$-orbit, and
$\upsilon$ is an isomorphism away from these components.  Because the only
fibers that have dimension greater than zero are curves of genus-zero, the
first derived push-forward vanishes 
$R^1\upsilon_* \co_{\ce_{g,n+1}} = 0,$
 and the push-forward is the trivial bundle 
 $$\upsilon_* \co_{\ce_{g,n+1}} =  \co_{\ce'} = \alpha^* \co_{\ce_{g,n}}.$$
Therefore, by the Leray spectral sequence, we have $$R^1(\pi_{g,n+1})_* \co_{\ce_{g,n+1}} = R^1(\pi'_{g,n})_* \co_{\ce'} = R^1(\pi'_{g,n})_* \alpha^* \co_\ce.$$  
Since the square is Cartesian and $\tau$ is flat, this gives
$$\crrt_{g,n+1} = R^1(\pi_{g,n+1})_* \co_{\ce_{g,n+1}} = \tau^* R^1(\pi_{g,n})_* \co_{\ce} = \tau^*\crrt_{g,n},$$  
as desired.

Equations (\ref{eq:clt-tails}) and (\ref{eq:psit-tails}) follow from 
their counterparts
on $\M_{g,n}$.  Specifically, it is 
well-known \cite[Eq (2.36)]{Wi} that   
$$\tau^*(\cl_i) = \cl_i\otimes \co([\M_{\Ga_i}]),$$
where $\Ga_i = |\Gat_i|$ is the undecorated graph underlying the forgetting-tails graph $\Gat_i$ of Definition~\ref{df:Gati}.
By Proposition~\ref{prop:PsitPsi}, we have
\begin{align*}
 \tau^*(\clt_i^{\otimes |m_i| }) 	
 &= \tau^*(\cl_i) \\
 &= \cl_i\otimes \st^* \co([\M_{\Ga_i}])\\
 &= \clt_i^{\otimes |m_i| } \otimes \st^*\co([\M_{\Ga_i}]).
 \end{align*}
However, by Proposition~\ref{prop:FundClasses}, we have that 
$$\st^*\co([\M_{\Ga_i}]) = \co(|m_i|[\MM_{\Gat_i}]) =\co([\MM_{\Gat_i}])^{\otimes |m_i| },$$ so taking $|m_i|$-th roots gives the desired relation.

To prove Equation~(\ref{eq:rhocrrtloop}), we first note that in the following diagram the two squares are Cartesian:
\begin{equation}\label{diag:normalize}
\begin{diagram}
&\ce_{\Gatc}\\
&\dTo^{\nu} & \rdTo^{\mut_{\Gatc}} \\
\ce':=&\mu_{\Gatc}^*(\ce)	& \rTo^a 	& \ce_{\Gat} & \rTo^{\itt_{\Gat}} & \ce\\
&\dTo^{\pi'}	& 		& \dTo^{\pi_{\Gat}}& & \dTo^{\pi}\\
&\MM_{\Gatc}			&\rTo^{\mu_{\Gatc}} & \MM_{\Gat} & \rTo^{i_\Gat} & \MM_{g,n}.
\end{diagram}
\end{equation}
The upper-left vertical map $\nu$ is the normalization of the admissible $G$-cover at the nodes cut in $\Gatc$, and $\pi_{\Gatc} = \pi'\circ \nu$.  Since the squares are Cartesian, we have $$\rho_{\Gatc}^* \crrt =\mu_{\Gatc}^* i_\Ga^* R^1\pi_* \co_\ce= R^1\pi'_*(a^*\itt_{\Gat}^* \co_\ce) = R^1\pi'_* \co_{\ce'}.$$

On $\ce'$, we have the following short exact sequence:
$$ 0\rTo \co_{\ce'} \rTo \nu_*\co_{\ce_{\Gatc}} \rTo \co_{nodes} \rTo 0.$$
Pushing forward to $\MM_{\Gatc}$, we have the following long exact sequence
in $\Rep(G) \otimes K(\MM_{\Gatc})$: 
\begin{align*}
 0\rTo &\nc[G/G_0] \otimes \co_{\MM_{\Gatc}}  \rTo \nc[G/G_{cut}] \otimes \co_{\MM_{\Gatc}} \rTo \nc[G/\gen{m_+}] \otimes \co_{\MM_{\Gatc}} \\ & \rTo  
  R^1\pi'_*\co_{\ce'} \rTo R^1(\pi_{\Gatc})_*\co_{\ce_{\Gatc}}\rTo 0 .
 \end{align*}
This gives the desired relation in K-theory.
The proof of Equation~(\ref{eq:rhocrrttree}) is similar.

Equation~(\ref{eq:rhoclt}) follows immediately from the fact that in
Diagram~(\ref{diag:normalize}) the map $\nu$ is an isomorphism in a
neighborhood of the image of the sections $\sigma_\pm$. 
\end{proof}

\section{The Hurwitz-Hodge bundle and a relative Riemann-Hurwitz Theorem}  

In this section we prove the relative Riemann-Hurwitz formula, which
generalizes the formula of \cite{JKK2} and allows us to write the equivariant
K-theory class of the Hurwitz-Hodge bundle in a useful form. 

We begin with a discussion in Subsection~\ref{sec:GenFun} of several generating functions that we will need to describe the result.  Next, in Subsection~\ref{sec:punc} we will use these generating functions to describe classes $\cs_{m_i}$ in $K_G(\MM_{g,n})$ which are associated to each puncture.  In the relative Riemann-Hurwitz formula, these classes describe the contribution from each puncture to the overall formula for $\crrt$.  As our last preliminary step, in Subsection~\ref{sec:node} we describe classes
$\cs_{\Gatc}$ in $K_G(\MM_{g,n})$
associated to each cut graph $\Gatc \in \bGatc$.  These classes describe the contribution from each node to the overall formula for $\crrt$.  In Subsection~\ref{sec:MainThm} we bring all these pieces together to state and explain the main theorem and its consequences, and in Subsection~\ref{sec:ProofOfMain} we prove the main theorem.

\subsection{Generating functions and Bernoulli polynomials}\label{sec:GenFun}

\begin{df}
Let $r\geq 1$ be an integer. Define the rational functions
\begin{equation}
H_r(y) := 
\frac{y^r-1}{y-1} = \sum_{k=0}^{r-1} y^k 
\end{equation}
\begin{equation}\label{eq:BDef}
F_r(x,y) := \frac{ H_r(x y) -
H_r(y) }{x^r - 1} ,
\end{equation}
and 
\begin{equation}
\F_r(t,y) := F_r(e^{t/r},y)
\end{equation}
\end{df}

Since
\begin{equation*}
F_r(x,y)  = \sum_{k=0}^{r-1} \frac{x^k - 1}{x^r - 1} y^k
= \sum_{k=1}^{r-1} \frac{\sum_{l=0}^{k-1} x^l}{\sum_{l=0}^{r-1} x^l} y^k = \sum_{k=1}^{r-1} \frac{H_k(x)}{H_r(x)} y^k,
\end{equation*}
we see that $F_r(x,y)$ is a polynomial of order $r-1$ in $y$, i.e., 
\begin{equation}\label{eq:CDef}
F_r(x,y) = \sum_{k=0}^{r-1} F_{r,k}(x) y^k,
\end{equation}
where $F_{r,k}(x)$ is the rational function
\begin{equation}\label{eq:CDefTwo}
F_{r,k}(x) = \frac{H_k(x)}{H_r(x)} = \frac{x^k - 1}{x^r - 1} = 
\begin{cases}
\frac{\sum_{l=0}^{k-1} x^l}{\sum_{l=0}^{r-1} x^l} & \text{if $k=1,\ldots,r-1$} \\
0 & \text{ if $k = 0$.}
\end{cases}
\end{equation}
Note also that one can expand $F_{r,k}(x)$ about $x=1$ to obtain the
following power series
\begin{equation}\label{eq:FrkExpansion}
\begin{split}
F_{r,k}(x) &= \frac{k}{r} + \frac{k(k-r)}{2 r} (x-1) 
+ \frac{k (k-r)(-r+2k-3)}{12 r} (x-1)^2 \\
&+ \frac{k(k-r) (k-2) (k-r-2)}{24 r} (x-1)^3
+ \mathrm{O}((x-1)^4).
\end{split}
\end{equation}
Combining Equations (\ref{eq:FrkExpansion}) and (\ref{eq:CDef}), we may thus
regard $F_r(x,y)$ as an element in $\nq[[x-1]][y]$.

A simple computation gives the following proposition.
\begin{prop}\label{prop:Brxzeta}
Let $r>0$. If $y\in\nc$ satisfies $y^r = 1$, but $y\neq 1$,  then we have 
\begin{equation}\label{eq:Brxzeta}
F_r(x,y) = \frac{1}{x y -1}.
\end{equation}
\end{prop}
\begin{rem}This relation for $F_r(x,y)$ is not true for $y$ in a general commutative ring, and in particular, it
is not true in the representation ring of $G$. 
\end{rem}

We can expand
\begin{equation}
\F_r(t,y) = \sum_{k=0}^{r-1} \F_{r,k}(t) y^k,
\end{equation}
where
\begin{equation}
\F_{r,k}(t) := \frac{e^{k t/r}-1}{e^t-1} 
= \sum_{j=0}^\infty 
\delta B_{j+1}\left(\frac{k}{r}\right)
t^j
=\int_0^{\frac{k}{r}} B(t,z) dz, 
\end{equation}
\begin{equation}\label{eq:DeltaB}
\delta B_n(z) := \frac{B_n(z)-B_n(0)}{n!},
\end{equation}
and $B_n(z)$ is the $n$-th Bernoulli polynomial,
defined by 
\[
B(t,z) := \frac{t e^{t z}}{e^t-1} = \sum_{n=0}^\infty B_n(z) \frac{t^n}{n!}.
\]
We have the well-known relations
\begin{equation} \label{eq:BerRel}
B(-t,x) = B(t, 1-x)
\dsand
(-1)^n B_n(0) = B_n(1).
\end{equation}
We also have $B_n(0)=0$ for all odd $n>1$.  Combining this with the
definition of $\delta B_n$ and Equation~(\ref{eq:BerRel})  gives  
\begin{equation}
\delta B_n(1-x) + (-1)^{n+1} \delta B_n(x) =  \delta_n^1 \label{eq:deltaBerRel}
\end{equation}
for all integers $n\geq 1$ and $x\in \nc$ where $\delta_n^1$ is the Kronecker
delta function. 

The first few terms are
\[
B(t,z) = 1 + \left(z -\frac{1}{2}\right) \,t + 
  \left( \frac{1}{6} - z + z^2     \right) \,\frac{t^2}{2!} + 
\left( \frac{z}{2} -  \frac{3 z^2}{2} + z^3 \right) \,\frac{t^3}{3!} +
  {\mathrm{O}(t)}^4 
\]
and
\[
\F_{r,k}(t) = \frac{k}{r} + \left( \frac{k^2}{2\,r^2} - 
     \frac{k}{2\,r} \right) \,t + 
  \left( \frac{k^3}{3\,r^3} - \frac{k^2}{2\,r^2} + 
     \frac{k}{6\,r} \right) \,\frac{t^2}{2!} + 
  \left( \frac{k^4}{4\,r^4} - \frac{k^3}{2\,r^3} + 
     \frac{k^2}{4\,r^2} \right) \,\frac{t^3}{3!} + 
  {\mathrm{O}(t)}^4.
\]
Thus, $\F_r(t,y)$ may be regarded as an element of $\nq[[t]][y]$.
\begin{prop}
If $r\geq 1$ is an integer, then
\begin{equation}\label{eq:FyOne}
\begin{split}
&F_r(x,1) = \frac{1}{x-1} - \frac{r}{x^r-1} \\
&= \frac{r-1}{2} - \frac{r^2-1}{12} (x-1) + \frac{r^2-1}{24} (x-1)^2 +
 \mathrm{O}((x-1)^3)
\end{split}
\end{equation}
Furthermore, for all $s\geq 0$, 
\begin{equation}\label{eq:BernoulliSum}
\sum_{k=0}^{r-1} B_{s}\left(\frac{k}{r}\right) = B_{s}(0) r^{1-s}.
\end{equation}
\end{prop}
\begin{proof}
Equation (\ref{eq:FyOne}) follows by performing the summation in Equation 
(\ref{eq:CDef}) after plugging in $y=1$. 

Plugging in 
$x = e^{t/r}$ 
into Equation (\ref{eq:FyOne}) yields $\F_r(t,1)$ in
$\nq[[t]]$ equal to
\begin{equation}\label{eq:ChiOneSm}
\frac{r}{t}\left(\frac{t/r}{e^{t/r}-1} - \frac{t}{e^t-1}\right) = 
\sum_{k=0}^{r-1} \frac{e^{k t/r}-1}{e^t-1}.
\end{equation}
Expressing both sides of this equality in terms of Bernoulli polynomials yields
\begin{equation}\label{eq:ChiOneSmExpanded}
\sum_{j\geq 0} B_{j+1}(0) \frac{t^j (r^{-j}-r)}{(j+1)!}
=\sum_{j\geq
0}\sum_{k=0}^{r-1}\frac{B_{j+1}\left(\frac{k}{r}\right)-B_{j+1}(0)}{(j+1)!}
t^j.
\end{equation}
The result follows by equating coefficients of $t^j$ for $j\geq 0$.
\end{proof}
\begin{df}
For any function $f(x,y)$, define its \emph{dual function}  $f^*(x,y)$  by
\begin{equation}
f^*(x,y) := f(x^{-1},y^{-1}).
\end{equation}
\end{df}
\begin{prop}
Let $r\geq 1$ be an integer, then for all $k=1,\ldots,r-1$, we have 
\begin{equation}\label{eq:CPlusCBar}
F_{r,k}(x) + F_{r,r-k}^*(x) = 1,
\end{equation}
and
\begin{equation}\label{eq:BxInvyInv}
F_r^*(x,y) = H_r(y^{-1}) -1 - y^{-r} F_r(x,y).
\end{equation}
\end{prop}
\begin{proof}
Equation (\ref{eq:CPlusCBar}) is immediate. Equation (\ref{eq:BxInvyInv})
follows since
\begin{align*}
F_r(x^{-1},y^{-1}) &= \sum_{k=1}^{r-1} F_{r,k}(x^{-1}) y^{-k} 
= \sum_{k=1}^{r-1} (1-F_{r,r-k}(x)) y^{-k} \\
&= \sum_{k=1}^{r-1} y^{-k} - \sum_{k=1}^{r-1} F_{r,r-k}(x) y^{-k} 
= H_r(y^{-1}) - 1 - \sum_{k=1}^{r-1} F_{r,k}(x) y^{-(r-k)} \\
&= H_r(y^{-1}) - 1 - y^{-r} F_r(x,y).
\end{align*}
\end{proof}

Another useful generating function is the following.
\begin{df}
Let $r\geq 2$ be an integer. Let
\begin{equation}\label{eq:IF}
\begin{split}
\IF_r(\xp,\xm,y) &:= F_r(\xp,y) 
F_r(\xm,y^{-1})
 - H_r(y) \sum_{k=0}^{r-1}
F_{r,k}(\xp) F_{r,k}(\xm) \\
&=: \sum_{k=-(r-1)}^{r-1} \IF_{r,k}(\xp,\xm) y^k
\end{split}
\end{equation}
and for all $k=0,\ldots,r-1$, let
\begin{equation}
\IIF_{r,k}(\xp,\xm) := 
\begin{cases}
\IF_{r,k}(\xp,\xm) + \IF_{r,-r+k}(\xp,\xm) &
\text{if $k=1,\ldots,r-1$} \\  \IF_{r,0}(\xp,\xm) &\text{if $k=0$}.
\end{cases}
\end{equation}
\end{df}
We may regard $\IF_r(\xp,\xm,y)$ as an element of
$\nq[y,y^{-1}][[\xp-1,\xm-1]]$ since $\IF_{r,k}(\xp,\xm)$ may be regarded as
an element in $\nq[[\xp-1,\xm-1]]$. 

\begin{prop}
For all integers $r\geq 2$ and $k=0,\ldots,r-1$,
\begin{equation}\label{eq:IIF}
\begin{split}
&\IIF_{r,k}(\xp,\xm) =
\frac{\frac{\xp^k-1}{\xp^r-1}+\frac{\xm^{r-k}-1}{\xm^r-1}-1}{\xp \xm - 1} \\
& = \frac{k(k-r)}{2 r} - \frac{k(k-r)(r-2k-3)}{12 r}\left( (\xp-1) + (\xm-1)
\right) + \mathrm{O}((x-1)^2)
\end{split}
\end{equation}
and
\begin{equation}\label{eq:IFyOne}
\begin{split}
&\IF_r(\xp,\xm,1) = \frac{1}{(\xp-1)(\xm-1)} - r \frac{\xp^r\xm^r -1}{(\xp
\xm-1) (\xp^r-1) (\xm^r-1)}\\
&= \frac{1-r^2}{12} + \frac{r^2-1}{24}\left( (\xp-1) +
(\xm-1)\right) + 
\frac{r^4 - 20 r^2 + 19}{720} \left( (\xp-1)^2 +
(\xp-1)^2\right) \\
&- \frac{(r+1)(r-1)(r^2+11)}{720} (\xp-1)(\xm-1) 
+ \mathrm{O}((x-1)^3),
\end{split}
\end{equation}
where both equalities may be regarded as in $\nq[[\xp-1,\xm-1]]$.

Furthermore, if $R_r$ is the commutative algebra $R_r := \nq[v]/\gen{v^r-1}$, 
then we have the following equality in $R_r[[\xp-1,\xm-1]]$:
\begin{equation}\label{eq:IFv}
\IF_r(\xp,\xm,v) = \sum_{k=0}^{r-1} \IIF_{r,k}(\xp,\xm) v^k.
\end{equation}
\end{prop}
\begin{proof}
Equation (\ref{eq:IFv}) follows immediately from Equation (\ref{eq:IIF})
and the second equality in Equation (\ref{eq:IF}).

Equation (\ref{eq:IFyOne}) follows from Equation (\ref{eq:IIF}) after
performing the summation
\[
\IF_r(\xp,\xm,1) = \sum_{k=0}^{r-1} \IIF_{r,k}(\xp,\xm).
\]

We will now prove Equation (\ref{eq:IIF}). Plugging in definitions yields the
equality
\begin{equation*}
\begin{split}
&(\xp^r-1)(\xm^r-1)\IIF_r(\xp,\xm,y) = 
\sum_{j=1}^{r-1} \left(
\sum_{\substack{\kp-\km=j\\\kpm\geq 0}} (\xp^\kp-1)(\xm^\km-1) y^j\right.
\\ &\left. 
+\sum_{\substack{\kp-\km=j-r\\\kpm\geq 0}} (\xp^\kp-1)(\xm^\km-1)
y^{j-r} -\sum_{k=0}^{r-1}(\xp^k-1)(\xm^k-1) y^j \right).\\
\end{split}
\end{equation*}
Since the coefficient of $y^0$ of the right hand side of this equation is
zero, we have $\IIF_{r,0}(\xp,\xm) = 0$. Similarly, picking off the coefficient of
$y^j$ for all $j\in\{1,\dots,r-1\}$, we obtain
\begin{equation*}
\begin{split}
&(\xp^r-1)(\xm^r-1)\IIF_{r,j}(\xp,\xm) \\ &
= \sum_{\substack{\kp-\km=j\\\kpm\geq 0}} (\xp^\kp-1)(\xm^\km-1) 
+\sum_{\substack{\kp-\km=j-r\\\kpm\geq 0}} (\xp^\kp-1)(\xm^\km-1) 
-\sum_{k=0}^{r-1}(\xp^k-1)(\xm^k-1) \\
&= \sum_{\substack{\kp+\km=j+r\\\kpm\geq 1}} (\xp^\kp-1)(\xm^{r-\km}-1) 
+\sum_{\substack{\kp+\km=j\\\kpm\geq 1}} (\xp^\kp-1)(\xm^{r-\km}-1) 
-\sum_{k=1}^{r-1}(\xp^k-1)(\xm^k-1) \\
&= \sum_{\substack{\kp+\km=j+r\\\kpm\geq 1}} (\xp^\kp-1)(\xm^{r-(j-\kp)}-1) 
+\sum_{\substack{\kp+\km=j\\\kpm\geq 1}} (\xp^\kp-1)(\xm^{r-\km}-1) 
-\sum_{k=1}^{r-1}(\xp^k-1)(\xm^k-1) \\
&= \sum_{\kp=j+1}^{r-1} (\xp^\kp-1)(\xm^{\kp-j}-1) 
+\sum_{\kp=1}^{j-1} (\xp^\kp-1)(\xm^{\kp+r-j}-1) 
-\sum_{k=1}^{r-1}(\xp^k-1)(\xm^k-1) \\
&= \sum_{\kp=1}^{r-1-j} (\xp^{\kp+j}-1)(\xm^\kp-1) 
+\sum_{\kp=1}^{j-1} (\xp^\kp-1)(\xm^{\kp+r-j}-1) 
-\sum_{k=1}^{r-1}(\xp^k-1)(\xm^k-1) \\
\end{split}
\end{equation*}
Performing the resulting summations and solving for $\IIF_{r,j}(\xp,\xm)$
yields Equation (\ref{eq:IIF}).
\end{proof}

\subsection{Contribution from the punctures}\label{sec:punc}
In this section, we define a class $\cs_{m_i}$ in $K_G(\MM_{g,n})$ associated
to the $i$th puncture of the universal pointed admissible $G$-cover.  This
class plays 
a central role
in the formula for $\crrt$.

\begin{df}
Let  $\bm = (m_1,\ldots,m_n)$ in $G^n$ and pick $i=1,\ldots,n$. Let
$r_i = |m_i|$. 
Define
\begin{equation}\label{eq:cstmkDef}
\cst_{m_i,k} := F_{r_i,k}(\clt_i) = \frac{H_{k}(\clt_i)}{H_{r_i}(\clt_i)} = \frac{ \clt_i^k-1}{\clt_i^{r_i} -1}
\end{equation}
in 
$\Kmgnloc{(\bm)}$,
as given in Equation~(\ref{eq:CDefTwo}). 
And define 
\begin{equation}\label{eq:csexpand}
\cst_{m_i}:= F_{r_i}(\clt_i, \V_{m_i}) = \sum_{k=0}^{r_i-1} F_{r_i,k}(\clt_i) \V_{m_i}^k 
\end{equation}
in 
$\Kmgnloc{(\bm)}
\otimes\Rep(\gen{m_i})$.   
Similarly, let
\begin{equation}\label{eq:Suv}
\cs_{m_i} := \Ind^G_{\gen{m_i}} \cst_{m_i} = \sum_{k=0}^{r_i-1}
          \cst_{m_i,k} \Ind^G_{\gen{m_i}} \V_{m_i}^k
\end{equation}
in 
$\Kmgnloc{(\bm)}
\otimes\Rep(G).$ 
\end{df}

\begin{rem}
Note that the definition of $\cst_{m_i,k}$ makes sense because the element $H_{r_i}(\clt_i)=  \sum_{k=0}^{r_i-1} \clt_i^k $ is invertible in $\Kmgnloc{(\bm)}$.  This can be seen by expanding $H_{r_i}(\clt_i)$ as a power series in $(\clt_i - 1)$, which has rank zero, and is, therefore, nilpotent in $\Kmgnloc{(\bm)}$. 
\end{rem}

We also define the following more general sheaf to represent the
contribution of the punctures to a more general equation for a 
push-forward of the form $R\pi_*\cf$ from the universal $G$-cover to $\MM_{g,n}$.   
\begin{df}
For any $\cf\in K_G(\ce)$ on the universal admissible $G$-cover
$\ce\rTo^{\pi}\MM_{g,n}(\bm)$, and for any $i \in \{1,\dots,n\}$, let  
\begin{equation}\label{eq:cstF}
\cst_{m_i}(\cf):=\cst_{m_i} \otimes \Phi(\sigma_i^*(\cf))
\end{equation}
in 
$\Kmgnloc{(\bm)}
\otimes\Rep(\gen{m_i}),$ where $\sigma_i^*(\cf)$ is
regarded as an $\gen{m_i}$-equivariant sheaf.

Similarly, let   
\begin{equation}\label{eq:csF}
\cs_{m_i}(\cf) := \Ind^G_{\gen{m_i}} \I^{\gen{m_i}}\cst_{m_i}(\cf) = \I^G
\Ind^G_{\gen{m_i}} \cst_{m_i}(\cf) 
\end{equation}
in 
$\Kmgnloc{(\bm)}
\otimes\Rep(G).$ 
\end{df}
It is easy to see that $$\cs_{m_i}(\co) = \cs_{m_i},$$ and because of the $\I^G$ in the definition we may apply Equation~(\ref{eq:Ivanish}) to see that 
\begin{equation}\label{eq:csfInv}
\cs_{m_i}(\cf)^G=0
\end{equation}
for all $\cf$.

\begin{df}
The \emph{dual functions} are defined as follows:
\[
\cst_{m_i, k}^*
:= F_{r_i,k}^*(\clt_i) = F_{r_i,k}(\clt_i^{-1}),
\]
\[
\cst_{m_i}^*:= F_{r_i}^*(\clt_i , \V_{m_i} ) = F_{r_i}(\clt_i^{-1}, \V_{m_i}^{-1}),
\]
and
\[
\cs_{m_i}^* := \Ind^G_{\gen{m_i}} \cst_{m_i}^*.
\]
Similarly, for every $\cf\in K_G(\ce)$ define
\begin{equation}
\cst_{m_i}^*(\cf):= \cst_{m_i}^* \otimes \Phi(\sigma_i^*(\cf^*)) = \cst_{m_i}^* \otimes \Phi(\sigma_i^*(\Hom(\cf,\co)))
\end{equation} 
in 
$\Kmgnloc{(\bm)}
\otimes\Rep(\gen{m_i}),$ and
\begin{equation}\label{eq:csFdual}
\cs_{m_i}^*(\cf) := \Ind^G_{\gen{m_i}} \I^{\gen{m_i}}\cst_{m_i}^*(\cf) = \I^G \Ind^G_{\gen{m_i}} \cst_{m_i}^*(\cf)
\end{equation}
in 
$\Kmgnloc{(\bm)}
\otimes\Rep(G).$ 
\end{df}
\begin{rem}
Note that $\cs_{m_i}^*(\cf)$ is precisely the dual of
$\cs_{m_i}(\cf)$, that is 
$$\cs_{m_i}^*(\cf) = \Hom(\cs_{m_i}(\cf), \co).$$
\end{rem}

\begin{prop}\label{prp:SmiPlusDual}
For any  $\bm$ in $G^n$, and for all $i=1,\ldots, r_i-1$, where $r_i = |m_i|$,
we have 
\begin{equation}\label{eq:chGroupRing}
\nc[\gen{m_i}] = \sum_{k=0}^{r_i-1} \V_{m_i}^k = H_r(\V_{m_i})
\end{equation}
in $\Rep(\gen{m_i})$. For all $k=1,\ldots,r_i-1$, 
\begin{equation}\label{eq:QuvPlusQuvBar}
F_{r_i,k}(\clt_i) + F^*_{r_i,k}(\clt_i)=1
\end{equation}
in 
$\Kmgnloc{(\bm)}
\otimes\Rep(G)$,
and
\begin{equation}\label{eq:SPlusSBar}
\cs_{m_i} + \cs_{m_i}^* = \nc[G] - \nc[G/\gen{m_i}]
\end{equation}
in  $\Kmgnloc{(\bm)}
\otimes\Rep(G)$.
Finally, we have
\begin{equation}\label{eq:brkSmi}
\brk(\cs_{m_i}) = \sum_{k=0}^{r_i-1} \frac{k}{r_i} \Ind_{\gen{m_i}}^G 
\V_{m_i}^k 
\end{equation}
and 
\begin{equation}\label{eq:rkSmi}
\rk(\cs_{m_i}) = \frac{|G|}{2} \left(1 - \frac{1}{r_i}\right).
\end{equation}
\end{prop}

\begin{rem}
In the special case where $\clt_i = \clt_i^{-1}$,
then we have $\cs_{m_i,k}^* = \cs_{m_i^{-1},k}$
and 
$
\cs_{m_i}^* = \cs_{m_i^{-1}},
$
so that
$
\cs_{m_i} + \cs_{m_i^{-1}} = \nc[G] - \nc[G/\gen{m_i}].
$
Two situations in which this case occurs are, first, on
$\MM_{0,3}(\bm),$ where $\clt_i = \co$, and second,  when 
$m_i = m_i^{-1}$.
\end{rem}

\begin{proof}[Proof of Proposition~\ref{prp:SmiPlusDual}]
Equation (\ref{eq:chGroupRing}) follows from Equation
(\ref{eq:GroupRingCyclic}). Equation (\ref{eq:QuvPlusQuvBar}) 
follows from Equation (\ref{eq:CPlusCBar}). 

Now note that 
\begin{equation}\label{eq:StuvPlusStuvBar}
\V_{m_i}^{-r_i} F_{r_i}(\clt_i , \V_{m_i}) + F_{r_i}^*(\clt_i , \V_{m_i})
=  H_{r_i}(\V_{m_i}^{-1})- 1, 
\end{equation}
in  $\Kmgnloc{(\bm)}
\otimes\Rep(\gen{m_i})$.
This 
follows from
\begin{eqnarray*}
F_{r_i}^*(\clt_i , \V_{m_i} )
&=& H_{r_i}(\V_{m_i}^* ) - 1 - (\V_{m_i}^{-1})^{r_i} F_{r_i}(\clt_i,\V_{m_i} ) \\
&=& H_{r_i}(\V_{m_i}^* ) - 1 - \V_{m_i}^{-r_i} F_{r_i}(\clt_i,\V_{m_i}),
\end{eqnarray*}
where we have used Equation (\ref{eq:BxInvyInv}) in the first 
line.

We now apply $\Ind^G_{\gen{m_i}}$ to the resulting equality to obtain
\begin{multline}\label{eq:SuvPlusSuvBar}
\Ind_{\gen{m_i}}^G \left(F_{r_i}(\clt_i , \V_{m_i} )\right) +
\Ind_{\gen{m_i}}^G \left(F_{r_i}^*(\clt_i , \V_{m_i} )\right) \\=
\Ind^G_{\gen{m_i}}
H_{r_i}^*(\V_{m_i} ) 
- \nc[G/\gen{m_i}]
\end{multline}
in  $\Kmgnloc{(\bm)}
\otimes\Rep(G)$.

Equation (\ref{eq:SPlusSBar}) follows immediately from Equation
(\ref{eq:SuvPlusSuvBar}) after using Equation (\ref{eq:chGroupRing})
and the fact that $\nc[\gen{m_i}] = \nc[\gen{m_i}]^*$ in $\Rep(\gen{m_i})$. 
Equations~(\ref{eq:brkSmi}) and (\ref{eq:rkSmi}) are easily seen by noticing
that the Chern character of 
$F_{r_i,k}(\clt_i)$ is $\F_{r_i,k}(c_1(\clt_i))$, and so 
$$\rk(F_{r_i,k}(\clt_i)) = \Ch_0(F_{r_i,k}(c_1(\clt_i))) = \F_{r_i,k}(0) = k/r_i. 
$$ 
\end{proof}

\begin{prop}
For any $\cf\in K_G(\ce)$ on the universal $G$-cover
$\ce\rTo^{\pi}\MM_{g,n}(\bm)$, and for any $i\in \{1,\dots, n\}$ we have  
\begin{eqnarray}\label{eq:CstFplusCstFdual}
\cst_{m_i}^*(\cf)  + \cst_{m_i}(\Hom(\cf,\omega)) &=& \frac{r-1}{r}\left(\rk(\cf)
\right) \nc[\gen{m_i}]\notag\\
&=& \cst_{m_i}^*(\Hom(\cf,\omega))  + \cst_{m_i}(\cf) .
\end{eqnarray}
\end{prop}
\begin{proof}
For any $\ga = m_i^l \neq 1$, we can apply Proposition~\ref{prop:Brxzeta}
to get
\begin{eqnarray*}
\chi_\ga(\cst_{m_i}^*(\cf))  &+& \cst_{m_i}(\Hom(\cf,\omega))) =
\frac{\chi_{\ga}(\sig_i^*(\cf^*))}{\zeta_{r_i}^{-l}\clt^*_i - 1} +
\frac{\chi_{\ga}(\sig_i^*(\cf^*) \otimes
\sig_i^*(\omega_\pi)))}{\zeta_{r_i}^{l}\clt_i - 1}\\ 
&=& \frac{\chi_{\ga}(\sig_i^*(\cf^*))}{\zeta_{r_i}^{-1}\clt^*_i - 1} +
\frac{\chi_{\ga}(\sig_i^*(\cf^*))  \zeta_{r_i}^{l} \clt_i}{\zeta_{r_i}^{l}\clt_i - 1}
=0. 
\end{eqnarray*}
By Equation (\ref{eq:GroupRingChar}) and the fact that every representation is completely determined by its characters, it follows that   
$\cst_{m_i}^*(\cf)  + \cst_{m_i}(\Hom(\cf,\omega))$ is a scalar multiple of $\nc[\gen{m_i}]$.  Now, applying $\chi_1$, we get 
\begin{align*}
\chi_1\left(\cst_{m_i}^*(\cf)  \right.&+
\left.\cst_{m_i}(\Hom(\cf,\omega))\right) \\  
&=\sum_{k=0}^{r_i-1} \rk( F_{r_i,k}(\clt^{-1}) )\, \rk(\sig_i^*\cf^*)
+  \sum_{k=0}^{r_i-1}\rk( F_{r_i,k}(\clt ))\, \rk(\Hom(\sig_i^*\cf,\clt_i)) \\
&=2\sum_{k=0}^{r_i-1} \frac{k}{r}\rk(\sig_i^*\cf) \quad = \quad (r-1) \rk(\sig^*_i(\cf)) \\
&= \chi_1\left(\frac{r-1}{r} \rk(\sig^*_i(\cf)) \nc[\gen{m_i}]\right). 
\end{align*}
This shows that first equality holds. The proof of the second equality is similar.
\end{proof}
\begin{crl}\label{cor:csfDuality}
For every $\cf \in K_G(\ce)$ we have 
\begin{equation}\label{eq:csfDuality}
\cs_{m_i}(\Hom(\cf,\omega)) = - \cs^*_{m_i}(\cf) .
\end{equation}
\end{crl}
\begin{proof}
Applying $\I^{\gen{m_i}}$ and $\Ind^G_{\gen{m_i}}$ to Equation~(\ref{eq:CstFplusCstFdual}) we have
\begin{eqnarray*}
\cs^*_{m_i}(\cf) + \cs_{m_i}(\Hom(\cf,\omega)) &=& \Ind^G_{\gen{m_i}}\I^{\gen{m_i}}\left(\cst_{m_i}^*(\cf)  + \cst_{m_i}(\Hom(\cf,\omega)) \right)\\
&=& \Ind^G_{\gen{m_i}}\left(\frac{r-1}{r}\,\rk(\cf)\,
\I^{\gen{m_i}}(\nc[\gen{m_i}]) \right)\\
&=& 0.
\end{eqnarray*}
\end{proof}

\subsection{Contribution from the nodes}\label{sec:node}

In this section, we define a class $\cs_{\Gatc}$ in $K_G(\MM_{g,n})$
associated to each cut graph $\Gatc \in \bGatc$.  This class also plays an
important role in the formula for $\crrt$.  

Consider $\bm$ in $G^n$ 
and let $\Gatc \in \bGatcgn(\bm)$ be a choice of a cut graph decorated
with cut edges decorated by $\mpp$ and $\mmm = \mpp^{-1}$. 
Let $\rho_\Gatc:\MM_{\Gatc}\rTo\MM_{g,n}(\bm)$ denote the associated gluing
morphism. 
\begin{df} 
Consider $\cst_{\mpp}$ in $\Kgatcloc\otimes\Rep(\gen{\mpp})$ and
$\cst_{\mmm}$ in $\Kgatcloc\otimes\Rep(\gen{\mmm})$. 
Let
\begin{equation}\label{eq:SGac}
\cs_\Gatc := - \frac{|\mpp|}{2|G|} \rho_{\Gatc *} \Ind^G_{\gen{\mpp}} \left( \I^{\gen{\mpp}}
(\cst_{\mpp} \cst_{\mmm}) \right) 
\end{equation} in $\Kmgnloc{(\bm)}\otimes\Rep(G),$
where $\I^{\gen{\mpp}}$ is given by Equation~(\ref{eq:IsupG}) (but extended 
to $K(\MM_\Gat)$ coefficients).
Furthermore, define $$\cs_\Gat:=\sum_{\Gatc\in\bGatc(\Gat)} \cs_\Gatc.$$
\end{df}

It is easy to see that the terms $\rho_{\Gatc* }\clt_\pm$ do not depend on
all the data of $\Gatc$, but rather are determined only by the glued graph
$\Gat$. This yields the following proposition.  
\begin{prop}
The element $\cs_\Gatc$ depends only on the glued graph $\Gat$, and for any
$\Gatc\in \bGatc(\Gat)$ we have 
\begin{equation}\label{eq:SGat}
\begin{split}
\cs_\Gat = -\frac{|m_+|}{|\Aut(\Gat)||Z_G(m_+)|}\rho_{\Gatc *} 
\Ind^G_{\gen{\mpp}} \left( \I^{\gen{\mpp}}
(\cst_{\mpp} \cst_{\mmm}) \right) 
\end{split}
\end{equation}
\end{prop}

By Equation~(\ref{eq:Ivanish}) the projection onto the
$G$-invariant part 
$\cs_\Gatc^G$ of $\cs_\Gatc$ in $\Kmgnloc{(\bm)}$ satisfies
\begin{equation}\label{eq:GInvSNodes}
\cs_\Gatc^G   = \cs_\Gat^G = 0.
\end{equation}

\begin{rem}
Both the rank and $\Rep(G)$-valued rank of $\cs_\Gat$ and $\cs_\Gatc$ are
zero since the sheaf is only supported over the codimension-one substack
$\MM_\Ga$ in $\MM_{g,n}$. 
\end{rem}

The formula for $\cs_\Gatc$ can be rewritten in a different form that is often easier to work with, as follows.
\begin{prop}\label{prop:SimpleSGatc}
We have the following equality in $\Kmgnloc{(\bm)}\otimes\Rep(G)$:
\begin{equation}\label{eq:SimpleSGatc}
\begin{split}
&\cs_\Gatc = -\frac{|\mpp|}{2 |G|} \rho_{\Gatc *} \Ind^G_{\gen{\mpp}}
\IF_{|\mpp|}(\cltp,\cltm,\V_\mpp)  \\
& = \sum_{k=0}^{|\mpp|-1} -\frac{|\mpp|}{2 |G|} \Ind^G_{\gen{\mpp}}\V_\mpp^k
 \rho_{\Gatc *}\left[
\left(\frac{\cltp^k-1}{\cltp^{|\mpp|}-1}+\frac{\cltm^{{|\mpp|}-k}-1}{\cltm^{|\mpp|}-1}-1\right)/\left(\cltp \cltm - 1\right)
\right]
\\
\end{split}
\end{equation}
\end{prop}
\begin{proof}
Equation (\ref{eq:SimpleSGatc}) follows from Equations (\ref{eq:IF}),
(\ref{eq:chGroupRing}), (\ref{eq:IIF}) and (\ref{eq:IFv}).
\end{proof}

Finally, we define the contribution from the nodes to a more general push-forward.
\begin{df}
For any graph $\Gatc\in \bGatcgn(\bm)$, and for any $\cf\in K_G(\ce)$ let 
\begin{equation*}
\cs_{\Gatc}(\cf):=-\frac{|\mpp|}{4|G|}\I^G\Ind^G_{\gen{\mpp}} \left(  \rho_{\Gatc *} \left[\Phi\left(\sigma_+^*(\cf_{\Gatc}) + \sigma_-^*(\cf_{\Gatc})\right) \otimes  ( \cst_{\mpp} \cst_{\mmm} )
\right]\right)
\end{equation*} in $\Kmgnloc{(\bm)}\otimes\Rep(G),$ 
where $\cf_\Gatc$ is the 
 pullback of
$\cf$ to the universal $G$-cover $\ce_{\Gatc} \rTo \MM_{\Gatc}$, and
$\sig_{\pm}^*(\cf_{\Gatc})$ is regarded as an $\gen{m_\pm}$-equivariant
sheaf instead of as a $G$-equivariant sheaf.
\end{df}

As before, it is easy to check directly from the definitions that 
for any $\cf \in K_G(\ce)$ the $G$-invariants of $\cs_{\Gatc}(\cf)$ vanish: 
\begin{equation}
\left(\cs_{\Gatc}(\cf)\right)^G = 0 \label{eq:csfGaInv}
\end{equation}
and for $\cf=\co$ we obtain the original $\cs_{\Gatc}$:
\begin{equation}
\cs_{\Gatc}(\co) = \cs_{\Gatc} \label{eq:SgatcO}
\end{equation}

\subsection{The Main Theorem and Its Consequences}\label{sec:MainThm}

The main theorem of this paper is the following.
\begin{thm}[Main Theorem]\label{thm:FullBF} For any 
$\cf \in K_G(\ce)$ 
let
  $\cfb:=(\pih_*\cf)^G$.  The  
following equality holds in $\Kmgnloc{(\bm)}\otimes\Rep(G)$:
\begin{equation}\label{eq:FullBF}
\Phi(R\pi_* \cf) = R\pibar_* \cfb \otimes \nc[G]  -
\sum_{i=1}^n \cs_{m_i}(\cf) - \sum_{\substack{\Gatc\in \\ \bGatcgn(\bm)}} \cs_\Gatc(\cf).
\end{equation}
\end{thm}

We will give the proof in Section~\ref{sec:ProofOfMain}, but first  
we make a few remarks about the theorem and its consequences. The following corollary is one of the  most important consequences of the Main Theorem.
\begin{crl}[A $\Rep(G)$-valued relative Riemann-Hurwitz Theorem]\label{cor:BF} The following 
equality holds in $\Kmgnloc{(\bm)}\otimes\Rep(G)$:
\begin{equation}\label{eq:BF}
\Phi(\crrt) = \co\otimes \nc[G/G_0] +  \left(
\crr - \co \right) \otimes \nc[G]  +
\sum_{i=1}^n \cs_{m_i} + \sum_{\substack{\Gatc\in \\ \bGatcgn(\bm)}} \cs_\Gatc,
\end{equation}
where $\crr$ is the pullback of the dual Hodge bundle from $\M_{g,n}$ 
and $\co$ is the structure sheaf of $\MM_{g,n}(\bm)$. Over each connected
component $\N_{g,n}(\bm)$ of $\MM_{g,n}(\bm)$, $G_0$ denotes 
a
subgroup of
$G$ which preserves a connected component 
of a fiber of the universal 
$G$-cover $\ce$
over $\N_{g,n}(\bm)$, and $\nc[G/G_0]$ is the $G$-module generated by
the cosets $G/G_0$. 

Similarly, we have
\begin{equation}\label{eq:BFdual}
\Phi(\crrt^*) = \co\otimes \nc[G/G_0] +  \left(
\crr^* - \co \right) \otimes \nc[G]  +
\sum_{i=1}^n \cs^*_{m_i} - \sum_{\substack{\Gatc\in \\ \bGatcgn(\bm)}} \cs_\Gatc.
\end{equation}

\end{crl}
\begin{proof}[Proof of Corollary~\ref{cor:BF}]
Equation~(\ref{eq:BF}) follows immediately from the Main Theorem applied to the structure sheaf $\co$, after using the fact that $R^0\pi_* \co_\ce = \co_{\MM_{g,n}}\otimes \nc[G/G_0]$.

To see Equation~(\ref{eq:BFdual}), first use Serre duality to see that $\Phi(\crrt^*) = \Phi(R\pi_*\omega_\ce) + \nc[G/G_0]\otimes \co$.
Now note that,  by the residue map, the dualizing sheaf at a node is trivial
$$\sigp^* \omega_{\ce_\Gatc} = \sigm^*  \omega_{\ce_\Gatc} = \co_{\ce_\Gatc}.$$

Therefore, we have
\begin{eqnarray}
\cs_\Gatc(\omega_\ce) &=& -\frac{|\mpp|}{4|G|}\I^G\left( \rho_{\Gatc *} \left[\left(2\co \right) \otimes  \Ind^G_{\gen{\mpp}} ( \cst_{\mpp} \cst_{\mmm} )\right]\right)\notag\\
&=& \cs_{\Gatc}.
 \end{eqnarray}
Also, by Corollary~\ref{cor:csfDuality} we have 
$$\cs_{m_i}(\omega_\ce) = -\cs_{m_i}^*.$$
Furthermore, we have 
$$(\pih_*\omega_{\ce})^G = \omega_{\cc}.$$ 
Now applying Theorem~\ref{thm:FullBF} and using the previous relations gives
\begin{eqnarray*}
\Phi(R\pi_* \omega_\ce) &=& R\pibar_* \omega_\cc \otimes \nc[G]  -
\sum_{i=1}^n \cs_{m_i}(\omega_\ce) - \sum_{\substack{\Gatc\in \\ \bGatcgn(\bm)}} \cs_\Gatc(\omega_\ce)\notag\\
&=& (\crr^*-\co) \otimes \nc[G]  + \sum_{i=1}^n \cs^*_{m_i} - \sum_{\substack{\Gatc\in \\ \bGatcgn(\bm)}} \cs_\Gatc,
\end{eqnarray*}
as desired.
\end{proof}

\begin{rem}
Taking the $\Rep(G)$-valued rank of Equation (\ref{eq:BF}) yields the
equality in $\Rep(G)$
\begin{equation}
\brk({\crrt}) = \nc[G/G_0] + (g-1) \nc[G] + \sum_{i=1}^n \brk({\cs_{m_i})},
\end{equation}
since $\brk(\cs_\Ga) = 0$. This is precisely Equation (8.4) from
\cite{JKK2}. 
\end{rem}
\begin{crl}[Relative Riemann-Hurwitz Formula]
The following relation holds in 
$\Kmgnloc{(\bm)}$.
\begin{equation}\label{eq:RelRH}
\crrt= \frac{|G|}{|G_0|} \co +  |G| \left(
\crr
- \co \right) + \sum_{i=1}^n 
\csh_{m_i}
+ 
\sum_{\Gatc\in \bGatcgn} \frac12\rho_{\Gatc *}\csh_\Gatc,
\end{equation}
where
\begin{equation}\label{eq:Cshmi}
\begin{split}
\csh_{m_i} &:= \frac{|G|}{r_i}
\left(\frac{1}{\clt_{i}-1}-\frac{r_i}{\clt_{i}^{r_i}-1} \right) 
\end{split}
\end{equation}
and
\begin{multline}\label{eq:cshGatc}
\csh_\Gatc := 
\frac{1}{( \clt_{+}-1) ( \clt_{-}-1)}  - \rpp \frac{\clt_{+}^{\rpp}
 \clt_{-}^{\rpp} - 1}{(\clt_{+} \clt_{-}-1)
(\clt_{+}^{\rpp}-1)(\clt_{-}^{\rpp}-1)}
\end{multline}
The rank of $\crrt$ is given by the usual Riemann-Hurwitz Formula:  
\begin{equation}\label{eq:RH}
\rk(\crrt) = \frac{|G|}{|G_0|} + (g-1) |G| + \sum_{i=1}^n \frac{|G|}{2} \left(1 -
\frac{1}{r_i}\right). 
\end{equation}
\end{crl}
\begin{proof}
Take the character $\chi_\one$ of Equation (\ref{eq:BF}) and then 
apply Equations (\ref{eq:FyOne}) and (\ref{eq:IFyOne}) to obtain the desired
result. Taking the rank of Equation (\ref{eq:RelRH}) yields the usual
Riemann-Hurwitz Formula~(\ref{eq:RH}). 
\end{proof}

\begin{rem}
Equation (\ref{eq:RelRH}) can be rewritten as
\begin{equation}\label{eq:GenRelRH}
\crrt= N
\co +  \deg(\pi) \left(
\crr
- \co \right) + \sum_{i=1}^n 
\csh_{m_i}
+ 
\sum_{\Gatc\in \bGatcgn} \frac12\rho_{\Gatc *}\csh_\Gatc,
\end{equation}
where 
$N$
is the number of connected components in a fiber of
$\ce\rTo^{\pih}\cc\rTo^{\pibar}\MM_{g,n}$, $\deg(\pih) = |G|$ is the degree of 
$\pih$, and $r_i$ the order of $i$-th ramification. 
It is tempting to interpret Equation (\ref{eq:GenRelRH}) for a more general family of
curves $\ce\rTo^{\pih}\cc\rTo^{\pibar}T$, and not necessarily just for a family of admissible $G$-covers.  In this case we would regard $\csh_{m_i}$ as the sum
of all ramifications of the family occuring on the smooth locus of the family
of curves, and the sum over $\csh_\Gatc$ as ramifications occuring on the
nodal locus of the family of curves.  
\end{rem}

\begin{crl}\label{cor:MumfRel}
The following $\Rep(G)$-valued generalization of Mumford's identity holds in
$\Kmgnloc{}\otimes \Rep(G)$: 
\begin{equation}\label{eq:MumfordIdentity}
\begin{split}
\Phi(\crrt) + \Phi(\crrt^*) =
  2\co \otimes \nc[G/G_0] + 
(2g-2+n)\co\otimes\nc[G]
- \left(\sum_{i=1}^n\co \otimes \nc[G/ \gen{m_i}]\right).\quad
\end{split}
\end{equation}
Similarly,  in $A(\MM_{g,n})\otimes \Rep(G)$ we have 
\begin{equation}\label{eq:ChMumfordIdentity}
\bCh(\crrt + \crrt^*) = 2 \nc[G/G_0] + 
(2g -2 + n)\nc[G]
- \left(\sum_{i=1}^n\nc[G/ \gen{m_i}]\right).\quad
\end{equation}
\end{crl}
\begin{proof}[Proof of Corollary~\ref{cor:MumfRel}]
By equations~(\ref{eq:BF}) and (\ref{eq:BFdual})
we have
$$ \Phi(\crrt) + \Phi(\crrt^*) = \co\otimes 2\nc[G/G_0] +  \left(
\crr+\crr^* - 2\co \right) \otimes \nc[G]  +
\sum_{i=1}^n \left(\cs_{m_i} +\cs_{m_i}^*\right). 
$$
Pulling back the usual Mumford identity $\Ch(\crr)+\Ch(\crr^*) = 2g$ from $\M_{g,n}$ to $\MM_{g,n}(\bm)$, and using the fact that $\Kmgnloc{}\cong A^{\bullet}(\MM_{g,n})$ we have 
\begin{equation}\label{eq:mumfK}
\crr + \crr^* = 2g \co
\end{equation} 
in $\Kmgnloc{}$.  Applying Equations~(\ref{eq:mumfK}) and (\ref{eq:SPlusSBar}), we obtain Equation~(\ref{eq:MumfordIdentity}).  Applying the $\Rep(G)$-valued Chern character gives Equation~(\ref{eq:ChMumfordIdentity}).
\end{proof}

The Mumford identity for the ordinary Hodge bundle implies that the positive,
even-dimensional components of its Chern character must vanish. Similarly, our generalization of the Mumford identity for the Hurwitz-Hodge bundle
yields the following result.

\begin{crl}\label{crl:ChParity}
Let $G$ be a finite group. Let $\Rep(G)$ be its representation ring and $\eta$ the metric in $\Rep(G)$. Let $\crrt$ be the Hurwitz-Hodge bundle in $K_G(\MM_{g,n})$. 
For all $W$ in $\Rep(G)$, define 
$\crrt[W] := \eta(W,\Phi(\crrt))$  in $\Kmgnloc{}$.

For all $j\geq 1$ and $W$ in $\Rep(G)$, we have the equality
\begin{equation}\label{eq:ChParity}
\Ch_j(\crrt[W]) + (-1)^j\Ch_j(\crrt[W^*]) = 0.
\end{equation}
In particular, when $W^* = W$ then $\Ch_{2s}(\crrt[W]) = 0$ for all $s \geq 1$. 
\end{crl}
\begin{proof}[Proof of Corollary \ref{crl:ChParity}]
To avoid notational clutter, identify $\crrt$ in $K_G(\MM_{g,n})$
 with its
image $\Phi(\crrt)$ in $\Kmgnloc{}\otimes\Rep(G)$. Since
$\Irrep(G) =: \{\,\ee_\alpha\}_{\alpha=1}^{|\Gb|}$ is an orthonormal basis for
$\Rep(G)$, we have
\[
\crrt = \sum_{\alpha=1}^{|\Gb|} \crrt^\alpha \otimes \ee_\alpha,
\]
where $\crrt^\alpha = \eta(\crrt,\ee_\alpha)$ in $\Kmgnloc{}$.

For all $j\geq 1$, Equation
(\ref{eq:ChMumfordIdentity}) yields 
\[
\bCh_j(\crrt)+\bCh_j(\crrt^*) = 0,
\]
which, after plugging in the expansion of $\crrt$, becomes
\[
\sum_{\alpha=1}^{|\Gb|} \Ch_j(\crrt^\alpha) \otimes \ee_\alpha + \sum_{\alpha=1}^{|\Gb|}
(-1)^j\Ch_j(\crrt^\alpha) \otimes (\ee_\alpha)^*  = 0,
\]
where we have used that $\Ch_j(\cf^*) = (-1)^j\Ch_j(\cf)$ 
for all  $\cf$ in $\Kmgnloc{}$. However,
\[
\sum_{\alpha=1}^{|\Gb|} \Ch_j(\crrt^\alpha) \otimes (\ee_\alpha)^* =
\sum_{\alpha=1}^{|\Gb|} \Ch_j(\eta(\crrt,\ee_\alpha)) \otimes (\ee_\alpha)^* = 
\sum_{\alpha=1}^{|\Gb|} \Ch_j(\eta(\crrt,\ee_\alpha^*)) \otimes \ee_\alpha,
\]
where the last equality holds because $\Irrep(G)$ is preserved by
dualization. Plugging this into the previous equation yields
\[
\sum_{\alpha=1}^{|\Gb|} \left(\Ch_j(\eta(\crrt,\ee_\alpha)) +
(-1)^j\Ch_j(\eta(\crrt,\ee_\alpha^*))\right) \otimes \ee_\alpha  = 0.
\]
Given any $W = \sum_\beta W^\beta \ee_\beta$ in $\Rep(G)$, apply
$\eta(W,\cdot)$ to the previous equation to obtain the desired result.
\end{proof}

\subsection{Proof of the Main Theorem}\label{sec:ProofOfMain}

We now proceed with the proof of the main theorem (Theorem~\ref{thm:FullBF}).

\begin{lm}\label{lm:BAnsatz} Suppose that the following equality holds in
$\Kmgnloc{(\bm)}\otimes\Rep(G)$:
\begin{equation}\label{eq:BAnsatz}
\Phi(R\pi_* \cf) =  \cz \otimes \nc[G] +
\ct,
\end{equation}
where $\cz$ belongs to 
$\Kmgnloc{(\bm)}$
and the $G$-invariant part $\ct^G$
of $\ct$ vanishes:
\[
\ct^G = 0.
\]
Then
\begin{equation*}
\cz =R\pibar_* \cfb.
\end{equation*}
\end{lm}
\begin{proof}[Proof of Lemma~\ref{lm:BAnsatz}]
Consider the maps $\ce\rTo^{\pih}\cc\rTo^{\pibar}\MM_{g,n}(\bm),$ where $\pi:=
\pibar\circ\pih$. Since 
\begin{equation*}
\cfb := (\pih_*(\cf))^G,
\end{equation*} we see that
\[
(\Phi(R\pi_*\cf))^G = \Phi(R\pibar_*( (
\pih_*\cf )^G )) = R\pibar_*
\cfb.
\]
Taking $G$-invariants of both sides of Equation (\ref{eq:BAnsatz})
and using the 
vanishing of $\ct^G$,
we obtain the desired result.
\end{proof}
We are now ready to prove the Main Theorem (Theorem~\ref{thm:FullBF}).  
By Lemma \ref{lm:BAnsatz}, we need only show that Equation (\ref{eq:BAnsatz})
holds, where 
\begin{equation}\label{eq:TDef}
\ct = -\sum_{i=1}^n \cs_{m_i}(\cf) - \sum_{\substack{\Gatc\in \\ \bGatcgn(\bm)}} \cs_\Gatc(\cf),
\end{equation}
since $\cs_{m_i}(\cf)$ and $\cs_{\Gatc}(\cf)$ are both $G$-invariant (see Equations~(\ref{eq:csfInv}) and (\ref{eq:csfGaInv})). 

Suppose that for all $\ga\not=\one$ in $G$, 
\begin{equation}\label{eq:CharOne}
\chi_\ga(\Phi(R\pi_*\cf)) = - \sum_{i=1}^n
\chi_\ga(\cs_{m_i}(\cf)) + \sum_{\Ga\in\bGa_{g,n}} \chi_\ga(\cs_\Ga(\cf)),
\end{equation}
then, by Equation (\ref{eq:GroupRingChar}) and the fact that every representation is completely determined by its characters, it follows that the two sides of
Equation (\ref{eq:FullBF}) must agree up to a term proportional to $\nc[G]$; that is,  Equation (\ref{eq:BAnsatz}) must hold for some $\cz$ in
$\Kmgnloc{(\bm)}$.

Therefore, we need to prove that Equation (\ref{eq:CharOne}) holds. This is an
application of the Lefschetz-Riemann-Roch Theorem (Theorem~\ref{thm:DM-LRR}), which in this case says
\begin{align}
\chi_{\ga}\circ\Phi(R\pi_*\cf) &= \ell_{\ga}\left(R\pi_*(\cf)\right)\notag\\
& = R\pi^{\ga}_*\left(L_f(\cf)\right)\notag\\
& = -\!\!\!\!\!\!\!\!\sum_{\substack{\text{components $D$ of}\\ \text{the fixed locus}}}
\left.R\left(\pi\right|_{D}\right)_*\left( \frac{\ell_{\ga}(\cf)}{\chi_{\ga} \circ \Phi (\lambda_{-1}(\conor_{D/\ce})) }\right).\label{eq:LRRcontr}
\end{align}
The components of the fixed-point locus are of two types: first, $G$-translates of the images $D_i$ of the sections $\sigma_i$; and second, components of the nodal locus over $\MM_{\Gat}$ for certain choices of $\Gat\in \bGat_{g,n}$.  We will address the two cases in the next two subsections.

\subsubsection{Contribution from the punctures}
For each $j\in\{1,\dots,n\}$ the translates of $D_j$ which are fixed by $\ga$ are of the form $gD_j$, where there is some choice of an integer $l$ with $m_j^l = g\ga g^{-1}$.  For such a translate, the conormal bundle $\conor_{gD_j/\ce}$ is just the restriction of the canonical dualizing bundle $\omega_{\pi}$ to $gD_j$, and the element $\ga$ acts on this conormal bundle as $\zeta_j^l =\exp(-2 \pi l i/|m_j|) \neq 1$. 
Moreover, the map $\pi^{\ga}$ restricted to $gD_j$ is an isomorphism, and
the translated section $g\sigma_j$ is its inverse.   Thus, the contribution
to the LRR formula from this translate is  
\begin{equation}
-R\left(\left.\pi^{\ga}\right|_D\right)_*\left(\frac{\ell_{\ga}(\cf)}{\chi_{\ga} \circ \Phi (\lambda_{-1}(\conor_{gD_j/\ce})) }\right)
 =
\frac{\chi_\ga\circ\Phi(\cf|_{gD_j})}{\co - \zeta^l(\sigma_j^*\omega_{\pi})} 
=
\frac{\chi_\ga\circ \Phi(\sig_j^*\cf)}{\co - \zeta^l \clt_j}. 
\end{equation}
The number of distinct translates of $D_j$ that correspond to a specific choice of $m_j^l$ is $$\frac{|Z_G(\ga)|}{|\gen{m_j}|} = \frac {|G|}{r_j |\gb|}$$ (see \cite[Pf of Lm 8.5]{JKK2} for details).  So summing over all $j\in\{1,\dots,n\}$ and all contributions from translates of $D_j$ gives
$$\frac {|G|}{|\gb|} \sum_{j=1}^n \frac{1}{r_j}\sum_{m_j^l \in \gb}\frac{\chi_\ga(\sig_j^*\cf)}{\co - \zeta^l \clt_j}.$$

On the other hand, it is well-known (e.g., \cite[ex 3.19]{FH}) that for any subgroup $H \le G$ and any representation $W \in \Rep(H)$, the character of the induced representation is $$\chi_{\ga}(\Ind_H^G(W)) 
 =\frac{|G|}{|H|} \sum_{{\sigma \in H \cap \gb }}
  \frac{\chi_{\sigma}(W)}{|\gb|}.$$
Applying this to $\chi_\ga(\cs_{m_j}(\cf))$, and using Proposition~\ref{prop:Brxzeta}, 
we have
\begin{align*}
\chi_\ga(\cs_{m_j}(\cf)) &=  \frac{|G|}{r_j |\gb|} \sum_{m_j^l \in \gb } F_{r_j}(\clt_j, \zeta_j^{l}) \chi_\ga(\sig_j^*(\cf))\\
&=  \frac{|G|}{r_j |\gb|} \sum_{m_j^l \in \gb} \frac{\chi_\ga(\sig_j^*(\cf))}{\co-\zeta^l\clt_j}.
\end{align*}
This shows that the contribution from the punctures in the LRR formula (\ref{eq:LRRcontr}) is precisely 
$$\sum_{j=1}^n\chi_{\ga}(\cs_{m_j}(\cf)).$$

\subsubsection{Contribution from the nodes}

Let $\Gatc\in\bGatcgn(\bm)$ be a cut graph labeled with $m$ on the $+$ half of the cut edge, and labeled with $m^{-1}$ on the $-$ half, such that $m^l = \ga$ for some integer $l$.  Let $$D_{\Gatc} =D_+ \sqcup D_- \rInto^{\jt_\Gatc} \ce_\Gatc$$ denote the union of the images of the tautological sections $\sig_{\pm}$ \emph{and all translates of those sections which have monodromy $m$ on the $+$ side and monodromy $m^{-1}$ on the $-$ side}.  That is, $D_+$ is the union of all $Z_G(m)$-translates of the image of the section $\sigma_+$ and $D_-$ is the union of all $Z_G(m^{-1})$-translates of the image of $\sigma_-$.  

Of course, the subgroup of $Z_G(m)$ which fixes a section is exactly $\gen{m}$, so $D_+$ and $D_-$ are each principal $Z_G(m)/\gen{m}$-bundles over $\MM_\Gatc$. 

Let $\Gat$ be the graph obtained by gluing the $+$ and $-$ tails of $\Gatc$, and let $$D_\Gat = \mut_\Gatc(D_\Gatc) \rInto^{\jt_\Gat} \ce_\Gat$$ denote the image of $D_\Gatc$ under the gluing morphism $\mut_\Gatc$.

Every node fixed by $\ga$ lies over the locus $\MM_{\Gat}$ for some such choice of $m$, $l$, and $\Gatc$.  Indeed, if there are no automorphisms of the graph $\Gat$, then there are two choices of $\Gatc$ (corresponding to two choices of tail that could be labeled $+$) that give the same $\Gat$ and the same node.  If $\Gat$ has an automorphism, then there is only one such $\Gatc$.  

It might seem more natural to index the fixed nodes by $\Gat$ instead of $\Gatc$, but that will cause problems later, since we need to track the actual monodromy $m$ of the  node (which $\Gatc$ does), and not just its conjugacy class $\mb$ (which is all that $\Gat$ can track).

We will use the following diagram throughout the rest of this section.
\begin{equation}\label{diag:LRRnode}
\begin{diagram}
D_{\Gatc}  & \rTo^{\mut_\Gatc} & D_{\Gat} & \rInto^{\itt_\Gat\circ\jt_{\Gat}} & \ce_{g,n}\\
\dTo^{\pi_\Gatc} \uTo_{\sigma_{\pm}} & &\dTo^{\pi_{\Gat}} & & \dTo^{\pi}\\
\MM_{\Gatc}& \rTo^{\mu_{\Gatc}} & \MM_{\Gat} & \rInto^{i_{\Gat}} & \MM_{g,n} 
\end{diagram}.
\end{equation}

\begin{rem}
The reader should beware that although there are two tautological sections $\sigma_{\pm}:\MM_{\Gatc} \rTo D_{\Gatc}\subset \ce_{\Gatc}$, there is not necessarily a section of $\pi_{\Gat}:D_{\Gat} \rTo \MM_{\Gat}$.   

Also, one should beware that the left-hand square of this rectangle is not Cartesian, and the morphisms
$\mu_\Gatc$ and $\pi_\Gat$ are also not always \'etale.  
However, over the open locus
$\cM^G_{\Gat}$ the morphism $\mu_\Gatc$ forms a principal 
$\Aut(\Gat) \times Z_G(m)/\gen{m}$-bundle, 
and the morphism $\pi_\Gat$ forms a principal
$Z_G(m)/\gen{m}$-bundle.  This follows because in each case the
points in a given fiber are all translates of one another, and a translate by
$g\in G$ only has monodromy $m$ if $g$ is in the centralizer $Z_G(m)$.  Also,
a translate by $g$ is the same as the original node or mark precisely if $g
\in \gen{m}$.  Finally, the action of $\Aut(\Gat)$ on a point
$[E,\pt_1,\ldots,\pt_n,\pt_+,\pt_-]$ is to interchange $\pt_+$ and $\pt_-$,
whereas the action of $Z_G(m)/\gen{m}$ takes
$[E,\pt_1,\ldots,\pt_n,\pt_+,\pt_-]$ to
$[E,\pt_1,\ldots,\pt_n,\ga\pt_+,\ga\pt_-]$.  These actions clearly commute, and
so $\mut_\Gatc$ forms a principal  $\Aut(\Gat) \times
Z_G(m)/\gen{m}$-bundle. 
\end{rem}

The morphism $\mut_\Gatc$ is \'etale of degree $2|\Aut(\Gat)| |Z_G(m)/\gen{m}|$, and $\pi_\Gatc$ is \'etale of degree $2|Z_G(m)/\gen{m}|$.  
This shows that the K-theoretic push-forward $R\pi_{\Gatc *}=\pi_{\Gatc *}$ is just $|Z_G(m)/\gen{m}|$ times the pullback $\sigma_+^* + \sigma_-^*$ :
\begin{equation}\label{eq:pistarsigma}
R\pi_{\Gatc *}\cg =  \frac{|Z_G(m)|}{|m|}\left( \sigma_+^* \cg+ \sigma_-^* \cF\right)
\end{equation}
for any $\cg \in K(D_\Gatc)$.

Finally, it is known (See \cite[\S1.1]{FaPa})  that the pullback to
$D_\Gatc$ of the conormal bundle of the inclusion $D_\Gat \rInto^{\itt_\Gat
\circ \jt_\Gat}\ce_{g,n}$ is the sum of two line bundles  
\begin{equation}
\mut_{\Gatc}^* \conor_{D_{\Gat}/\ce} = \L_+ + \L_-
\end{equation}
which have the property that 
\begin{equation}\label{eq:ConorPullback}
\sigma_+^*(\L_+ + \L_-) = \sigma_-^*(\L_+ + \L_-) = \sigma_+^*(\omega_{\pi_\Gatc}) + \sigma_-^*(\omega_{\pi_\Gatc}) = \clt_+ + \clt_-.
\end{equation}

Using Diagram~(\ref{diag:LRRnode})  
we see that the contribution to the LRR formula from $D_{\Gat}$ is 
\begin{align}\label{eq:LRRpullback}
(i_{\Gat})_*R(\pi_{\Gat})_* &\left(\frac{\ell_{\ga}(\cf)}{\chi_{\ga} \circ \Phi (\lambda_{-1}(\conor_{D_{\Gat}/\ce})) }\right) \notag \\
&=\frac{1}{\deg(\mut_\Gatc)} (i_{\Gat})_*R(\pi_{\Gatc})_* \mut_{\Gatc *}\mut_\Gatc^*\left(\frac{\chi_{\ga}(\cf|_{D_\Gat})}{\chi_{\ga} \circ \Phi (\lambda_{-1}(\conor_{D_{\Gat}/\ce})) }\right) \notag \\
&=\frac{1}{\deg(\mut_\Gatc)} (i_{\Gat})_* (\mu_{\Gatc})_*R\pi_{\Gatc *} \left(\frac{\chi_{\ga}(\cf_{\Gatc})}{\chi_{\ga} \circ \Phi (\lambda_{-1}(\mut_\Gatc^*\conor_{D_\Gat/\ce})) }\right) \notag \\
&=\frac{|Z_G(m)/\gen{m}|}{\deg(\mut_\Gatc)} \rho_{\Gatc *}\left( \sigma_+^* \left(\frac{\chi_{\ga}(\cf_{\Gatc})}{\chi_{\ga} \circ \Phi (\lambda_{-1}(\L_+ + \L_-)) }\right) + \sigma_-^* \left(\frac{\chi_{\ga}(\cf_{\Gatc})}{\chi_{\ga} \circ \Phi (\lambda_{-1}(\L_+ + \L_-)) }\right)\right) \notag \\
&=\frac{1}{2|\Aut(\Gat)|}  \rho_{\Gatc *} \left(\frac{\chi_{\ga}(\sigp^*\cf_{\Gatc} + \sigm^*\cf_{\Gatc})}{\lambda_{-1}(\zeta^l
\clt_+ + \zeta^{-l}\clt_-) }\right) \notag \\ 
&=\frac{1}{2|\Aut(\Gat)|} \rho_{\Gatc *} \left(\left(\frac{\co}{\co-\zeta^l
\clt_+}\right)\left(\frac{\co}{\co-  \zeta^{-l}\clt_- }\right) 
\chi_{\ga}(\sigp^*\cf_{\Gatc} + \sigm^*\cf_{\Gatc})\right)
\notag.
\end{align}

Now, for a given choice of node in $\ce$  with monodromy $m$, there are exactly $2/|\Aut(\Gat)|$ choices of $\Gatc$ labeled with $m$ that correspond to the node, so summing over all possible nodes fixed by $\ga$ corresponds to summing over all $\Gatc$ and dividing by $2/|\Aut(\Gatc)|$.  Thus, for a given $\ga$,  the contribution to the LRR formula from the nodes is 
$$\sum_{m\in G} \sum_{l: m^l = \ga} \sum_{\substack{\Gatc\in\\ \bGatcgn(\bm,m,m^{-1})}} \frac{1}{4} \rho_{\Gatc *}
\left(\left(\frac{\co}{\co-\zeta^l \clt_+}\right)\left(\frac{\co}{\co-
\zeta^{-l}\clt_- }\right)\chi_{\ga}(\sigp^*\cf_{\Gatc} + \sigm^*\cf_{\Gatc})\right).$$ 
If we replace $m$ by $m'\in\mb$, or if we replace $\ga$ by $\ga' \in \gb$ the bundles $\rho_{\Gatc *}\clt_{\pm}$ are unchanged, and the terms of the form $\zeta_m^l$ are unchanged, so we can sum over $m$ and $l$ such that $m^l\in \gb$ and divide by the order of the conjugacy class $\gb$.  Thus the contribution to the LRR formula can be rewritten as
\begin{align}
&\frac{1}{|\gb|}
\sum_{m^l \in \gb}
\sum_{\substack{\Gatc\in\\ \bGatcgn(\bm,m,m^{-1})}} \frac{1}{4} \rho_{\Gatc *}
\left(\left(\frac{\co}{\co-\zeta^l \clt_+}\right)\left(\frac{\co}{\co-
\zeta^{-l}\clt_- }\right)\chi_{\ga}(\sigp^*\cf_{\Gatc} + \sigm^*\cf_{\Gatc})\right)\notag \\
&= \frac{1}{4}\sum_{\substack{\Gatc\in\\ \bGatcgn(\bm)}}\frac{1}{|\gb|}
\sum_{l:\mpp^l \in \gb}
 \rho_{\Gatc *}
\left(\left(\frac{\co}{\co-\zeta^l \clt_+}\right)\left(\frac{\co}{\co-
\zeta^{-l}\clt_- }\right)\chi_{\ga}(\sigp^*\cf_{\Gatc} + \sigm^*\cf_{\Gatc})\right)\notag \\
&= \sum_{\substack{\Gatc\in\\ \bGatcgn(\bm)}} \frac{|\mpp|}{4|G|} \chi_{\ga} \I^{G}\Ind_{\gen{\mpp}}^G  \rho_*
\left(\cst_{m_+}\cst_{m_-}(\sigp^*\cf_{\Gatc} + \sigm^*\cf_{\Gatc})\right)\notag\\ 
&=\sum_{\substack{\Gatc\in\\ \bGatcgn(\bm)}} \chi_{\ga} \cs_{\Gatc}(\cf).
\end{align}
This completes the proof of Theorem~\ref{thm:FullBF}.

\subsection{Group automorphisms and the Hurwitz-Hodge bundle} 

In this section, we study the action of the automorphism group of $G$
on the moduli space of $G$-covers and the Hurwitz-Hodge bundle.

Let $G$ be a group. Let $\Aut(G)$ denote the automorphism group of $G$. Given
any element $\ga$ in $G$, the map $\Ad_\ga:G\rTo G$, which takes $m\mapsto
\ga m \ga^{-1}$ for all $m$ in $G$, is an \emph{inner automorphism of
  $G$}. The group $\In(G)$ of of all inner automorphisms of $G$ is a normal
subgroup of $\Aut(G)$. The \emph{outer automorphism group $\Out(G)$ of $G$}
is the quotient group $\Aut(G)/\In(G)$. 

We will now describe an action of $\Aut(G)$ on $\MM_{g,n}$ and $\M_{g,n}(\BG)$.
Let $(E,\varrho;\pt_1,\ldots,\pt_n)$ denote a pointed admissible $G$-cover
with monodromies $\bm = (m_1,\ldots,m_n)$, where the $G$-action on $E$ is denoted
by $\varrho:G\rTo\Aut(E)$. Given any $\theta$ in $\Aut(G)$,
$(E,\varrho\circ\theta^{-1};\pt_1,\ldots,\pt_n)$ is also a pointed admissible
$G$-cover, but with monodromies $\theta(\bm) := (\theta(m_1),\ldots,\theta(m_n))$.
Furthermore, if  $$f:(E,\varrho;\pt_1,\ldots,\pt_n)\rTo
(E',\varrho';\pt_1',\ldots,\pt_n')$$  is a morphism which is $G$-equivariant with respect to
the $G$ actions $\varrho$ and $\varrho'$,  then the same map
$f:(E,\varrho\circ\theta^{-1};\pt_1,\ldots,\pt_n)\rTo
(E',\varrho\circ\theta^{-1};\pt_1',\ldots,\pt_n')$ is $G$-equivariant with
respect to the $\varrho\circ\theta^{-1}$ and $\varrho'\circ\theta^{-1}$ $G$ actions. 
Hence, $\Aut(G)$ acts on the category of pointed admissible
$G$-covers. 

Since the same discussion applies to families of (pointed) admissible $G$-covers,
we obtain
an action $L:\Aut(G)\rTo \Aut(\MM_{g,n})$ of $\Aut(G)$ on
$\MM_{g,n}$,  via
\[
L(\theta)(E;\varrho;\pt_1,\ldots,\pt_n) :=
(E;\varrho\circ\theta^{-1};\pt_1,\ldots,\pt_n). 
\]
Furthermore, for all $\bm$ in $G^n$, we have
$L(\theta):\MM_{g,n}(\bm)\rTo\MM_{g,n}(\theta(\bm))$. 

The action of $\Aut(G)$ on $\M_{g,n}(\BG)$, also denoted by
$L(\theta):\M_{g,n}(\BG)\rTo\M_{g,n}(\BG)$, is defined in the same way. Since
$\theta$ respects the conjugacy classes of $G$, it induces an action on
$\Gb$. Therefore, $L(\theta)$ takes
$\M_{g,n}(\BG;\bmb)\rTo\M_{g,n}(\BG;\theta(\bmb))$ for all $\bmb$ in $\Gb^n$,
where $\theta(\bmb)$ is defined by acting with $\theta$ componentwise.

A similar construction endows the category of $G$-modules with an action of
$\Aut(G)$. Therefore, $\Rep(G)$ is an $\Aut(G)$-module, 
where the map
$L(\theta):\Rep(G)\rTo\Rep(G)$ 
preserves the multiplication, the pairing, and
dualization on $\Rep(G)$. 

This action of $\Aut(G)$ on $\Rep(G)$ factors
through the action of $\Out(G)$, since if $\varrho:G\rTo\Aut(W)$ is a
$G$-module and $\theta = \Ad_\ga$ is an inner automorphism for some $\ga$ in
$G$, then $\varrho\circ\Ad_{\ga^{-1}}:G\rTo\Aut(W)$ is another $G$-module which is
isomorphic to $\varrho$ under the isomorphism $\varrho(\ga)$.

The action of $\Aut(G)$ on $\MM_{g,n}$ induces an action of $\Aut(G)$ on
$K_G(\MM_{g,n})$ and $K(\MM_{g,n})$ and the map 
$\Phi:K_G(\MM_{g,n})\rTo K(\MM_{g,n})_{\loc}\otimes\Rep(G)$ 
is $\Aut(G)$-equivariant. Since the monodromies change
by conjugation under the action of an inner automorphism, $\In(G)$ need not
act trivially upon $K_G(\MM_{g,n})$ or $K(\MM_{g,n})$.  However, if $\theta =
\Ad_\ga$ is an inner automorphism, then $\theta(\bmb) = \bmb$ and $L(\Ad_\ga):
\M_{g,n}(\BG)\rTo\M_{g,n}(\BG)$ induces the identity map on both
$K_G(\M_{g,n}(\BG))$ and $K(\M_{g,n}(\BG))$, since twisting the group action
$\varrho$ by an inner automorphism $\Ad_\ga$ yields a new group action
$\varrho\circ\Ad_{\ga^{-1}} = \varrho(\ga^{-1})\circ\varrho\circ\varrho(\ga)$  that
is isomorphic to $\varrho$ via the isomorphism $\varrho(\ga)$. Therefore, the
action of $\Aut(G)$ on $K_G(\M_{g,n}(\BG))$ and $K(\M_{g,n}(\BG))$ factors
through an action of the outer automorphism group $\Out(G)$.  

\begin{prop}
Let 
$L:\Aut(G)\rTo \Aut(K_G(\MM_{g,n}))$ 
be the action of $\Aut(G)$ induced from
its action on $\MM_{g,n}$. 
Let 
$\widetilde{L}:\Aut(G)\rTo \Aut(K_G(\ce))$
 be the action of $\Aut(G)$ induced
from its action on the universal $G$-cover $\ce\rTo^{\pi}\MM_{g,n}$. 

For all $\theta$ in $\Aut(G)$
and $\cf$ in $K_G(\ce)$,
 the following properties hold:
\begin{enumerate}
\item We have the equality in $K_G(\MM_{g,n})$
\begin{equation}\label{eq:AutInvGeneral}
L(\theta)R\pi_*\cf = R\pi_*\widetilde{L}(\theta)\cf.
\end{equation}
In particular, the class of the Hurwitz-Hodge bundle $\crrt$ in
$K_G(\MM_{g,n})$ satisfies 
\begin{equation}\label{eq:AutInvcrrt}
L(\theta)\crrt = \crrt. 
\end{equation}
Furthermore, Equations (\ref{eq:AutInvGeneral}) and (\ref{eq:AutInvcrrt})
continue to hold if $\MM_{g,n}$ is everywhere replaced by $\M_{g,n}(\BG)$
above.
\item For all $i=1,\ldots,n$ and $\bm$ in $G^n$,
we have the equality in $K(\MM_{g,n}(\theta(\bm)))_{\loc}\otimes\Rep(G)$,
\begin{equation}\label{eq:LSm}
L(\theta) \cs_{m_i}(\cf) = \cs_{\theta(m_i)}(\widetilde{L}(\theta)\cf).
\end{equation}
\item 
For all $\Gatc$ in $\bGatcgn(\bm)$, 
we have the equality in $K(\MM_{g,n}(\theta(\bm)))_{\loc}\otimes\Rep(G)$,
\begin{equation}\label{eq:LSGat}
L(\theta)\cs_\Gatc(\cf) =       \cs_{\theta(\Gatc)}(\widetilde{L}(\theta)\cf),
\end{equation}
where $\theta(\Gatc)$ in  $\bGatcgn(\theta(\bm))$
replaces  all decorations of $\Gatc$ 
by  their images under $\theta$.
\item 
Using the notation from Theorem \ref{thm:FullBF} and Corollary \ref{cor:BF},
the following equalities holds in $\Kmgnloc{(\theta(\bm))}\otimes\Rep(G)$:
\begin{equation}\label{eq:LFullBF}
L(\theta)\Phi(R\pi_* \cf) = R\pibar_* \cfb \otimes \nc[G]  -
\sum_{i=1}^n \cs_{\theta(m_i)}(\widetilde{L}(\theta)\cf) -
\sum_{\substack{\Gatc\in \\ \bGatcgn(\theta(\bm))}}
\cs_{\theta(\Gatc)}(\widetilde{L}(\theta)\cf),
\end{equation}
and
\begin{equation}\label{eq:LBF}
L(\theta)\Phi(\crrt) = \co\otimes \nc[G/G_0] +  \left(
\crr - \co \right) \otimes \nc[G]  +
\sum_{i=1}^n \cs_{\theta(m_i)} + \sum_{\substack{\Gatc\in
    \\ \bGatcgn(\theta(\bm))}} \cs_\Gatc.
\end{equation}
\end{enumerate}
\end{prop}
\begin{proof}
Consider $\theta$ in $\Aut(G)$. We have the commutative diagram
\begin{equation}\label{eq:LtAndL}
\begin{diagram}
\ce(\bm) & \rTo^{\widetilde{L}(\theta)} & \ce(\theta(\bm))\\
\dTo^{\pi_\bm}	& & \dTo^{\pi_{\theta(\bm)}}\\
\MM_{g,n}(\bm) & \rTo^{L(\theta)} & \MM_{g,n}(\theta(\bm)),
\end{diagram}
\end{equation}
where $\pi_\bm$ and $\pi_{\theta(\bm)}$ are the universal $G$-covers and
$\widetilde{L}(\theta)$ is the canonical lift of 
the isomorphism
$L(\theta)$. Furthermore,
$\widetilde{L}(\theta)$ takes the $G$-action on $\ce(\bm)$ to the $G$-action
on $\ce(\theta(\bm))$, \ie\ if $\varrho_\bm:G\rTo\Aut(\ce(\bm))$ and
$\varrho_{\theta(\bm)}:G\rTo\Aut(\ce(\theta(\bm)))$ are the group actions,
then 
\[
\varrho_{\theta(\bm)}(\ga) =
\widetilde{L}(\theta)\circ\varrho_\bm(\ga)\circ\widetilde{L}^{-1}(\theta)
\]
for all $\ga$ in $G$. In other words, $\widetilde{L}(\theta)$ is a
$G$-equivariant 
isomorphism which induces the isomorphism
$\widetilde{L}(\theta):K_G(\ce(\bm))\rTo K_G(\ce(\theta(\bm)))$. 
Equation (\ref{eq:AutInvGeneral}) follows from the fact that $L(\theta)$ and
$\widetilde{L}(\theta)$ are $G$-equivariant isomorphisms for all $\theta$ in
$\Aut(G)$.  Equation (\ref{eq:AutInvcrrt}) arises when $\cf$ is chosen to be
the structure sheaf $\co$. The same arguments hold for $\M_{g,n}(\BG)$. This
proves the first claim.

The second claim follows from the fact that the action of $\Aut(G)$ on
$\Rep(G)$ takes 
\[
\Ind^G_{m_i} \V_{m_i}^k \mapsto \Ind^G_{\theta(m_i)} \V_{\theta(m_i)}^k,
\]
while the action of $L(\theta)$ on $\MM_{g,n}(\bm)\rTo\MM_{g,n}(\theta(\bm))$ 
preserves $\clt_i$ for all $i\in\{1,\ldots,n\}$, 
since the definition of $\clt_i$ in $K(\MM_{g,n})$ is independent of the $G$-action. 
This proves the
second claim.  The proof of the third claim is similar.

The proof of the last statement follows from the fact that the commutative
diagram (\ref{eq:LtAndL}) has rows that are isomorphisms and
$\widetilde{L}(\theta)$ is $G$-equivariant. Therefore, in the proof of 
Equation (\ref{eq:BF}), each of the terms is mapped to its counterpart by
$L(\theta)$ in the Lefschetz-Riemann-Roch Theorem.
\end{proof}

\section{Two Chern characters of the Hurwitz-Hodge bundle}

In this section, we introduce certain tautological classes on $\MM_{g,n}$ and
calculate the Chern characters, $\Ch$ and $\bCh$, of the Hurwitz-Hodge
bundle $\crrt$. There are actually at least two distinct ways to compute the
Chern character of the Hurwitz-Hodge bundle.  The first method is to use
Grothendieck-Riemann-Roch and adapt the arguments of 
\cite[\S5]{Mum} to obtain $\Ch(\crrt)$. However, this will not yield any
information about the representations of $G$.  Our second method of computing
the Chern character is to apply 
a more refined
Chern Character $\bCh$ to Corollary~\ref{cor:BF}
which permits
us to track the representation theory.   
We will show that the former result can be obtained from the latter.

\subsection{Computation of the Chern character using GRR}

We continue to use the notation $r_i:=|m_i|$ and  $\rpp:=|\mpp| = |\mmm|$.  

Applying the Grothendieck-Riemann-Roch theorem directly to the definition of the Hurwitz-Hodge bundle yields the following theorem.
\begin{thm}\label{thm:GRRChernCh}
The Chern character $\Ch(\crrt)$ of the dual $\crrt$  of the Hurwitz-Hodge
bundle on $\MM_{g,n}$ satisfies the following equality in $A^\bullet(\MM_{g,n})$:
\begin{multline}
\Ch(\crrt) = |G/G_0| + \sum_{\ell=1}^\infty
\frac{B_\ell(0)}{\ell!}\left[-|G|\kappa_{\ell-1}  
+ \sum_{i=1}^n \frac{|G|}{r_i^\ell}\psi_i^{\ell-1} 
\right. \\ 
\left. 
 -\frac12 \sum_{\Gatc\in\bGatc(\bm)}\frac{1}{\rpp^{\ell-2}} \rho_{\Gatc *}
 \left(\sum_{q=0}^{\ell-2}(-1)^q \psip^q \psim^{\ell-2-q} \right)
\right], 
 \end{multline}
as in Section~\ref{sec:MainThm}, over each connected component $\N_{g,n}(\bm)$ 
of $\MM_{g,n}(\bm)$, we denote by $G_0$ the subgroup of
$G$ which preserves a connected component of a fiber of the universal
$G$-curve over $\N_{g,n}(\bm)$, and the sum over $q$ is understood to be zero
when $\ell = 1$.
\end{thm}
\begin{proof}
Recall the diagram~(\ref{diag:LRRnode}) used in the proof of the main theorem.
\begin{equation*}
\begin{diagram}
D_{\Gatc}  & \rTo^{\mut_\Gatc} & D_{\Gat} & \rInto^{\itt_\Gat\circ\jt_{\Gat}} & \ce_{g,n}\\
\dTo^{\pi_\Gatc} \uTo_{\sigma_{\pm}} & &\dTo^{\pi_{\Gat}} & & \dTo^{\pi}\\
\MM_{\Gatc}& \rTo^{\mu_{\Gatc}} & \MM_{\Gat} & \rInto^{i_{\Gat}} & \MM_{g,n} 
\end{diagram}.
\end{equation*}
The singular locus, $\mathrm{Sing}$, of $\ce$ consists of the union of the
images of all the loci $D_\Gatc$  in $\ce$, but each of these loci
$\rhot_{\Gatc}(D_\Gatc) = \itt_\Gat\circ\jt_{\Gat}\circ
{\mut_\Gatc}(D_\Gatc)$ 
appears
twice if there are no automorphisms
of $\Gat$, since there is a choice of which side of the cut edge to label
with $+$.   

Let $\conor_{\mathrm{Sing}}$ be the conormal bundle of the singular locus ${Sing}$ in $\ce$, and let $P$ be Mumford's polynomial \cite[Lem 5.1]{Mum}:
$$P(A_1+A_2,A_1\cdot A_2) = \sum_{\ell=1}^\infty \frac{(-1)^\ell B_\ell(0)}{\ell!}
\left(\frac{A_1^{\ell-1} +  A_2^{\ell-1}}{A_1+ A_2}\right) $$
using the convention that
$$\left(\frac{A^{\ell-1} +  B^{\ell-1}}{A+B}\right) := 0 \qquad \text{ when }
\ell = 1.$$  
Notice that if $s \geq 1$ then
\begin{equation}\label{eq:PSeries}
\left(\frac{A^s +  B^s}{A+B}\right) = \sum_{q=0}^{s-1} (-1)^ q A^q B^{s-1-q}.
\end{equation}

The Grothendieck-Riemann-Roch theorem states that 
$$\Ch(R\pi_* \omega_\pi) = \pi_*(\Ch(\omega_\pi)\Td^{\vee}(\Omega_\pi)).$$
Combining the argument of \cite[\S5]{Mum} with an argument similar to that given in Equation~(\ref{eq:LRRpullback}), we have
\begin{align*}
&\Ch(R\pi_* \omega_\pi) 
= \sum_{\ell=1}^\infty \frac{(-1)^\ell B_\ell(0)}{\ell!}\kat'_{\ell-1} 
+ \pi_*(\Td^{\vee}(\co_{\mathrm{Sing}}) -1)\notag \\
&=\sum_{\ell=1}^\infty \frac{(-1)^\ell B_\ell(0)}{\ell!}\kat'_{\ell-1} 
+ \pi_*P(c_1(\conor_{\mathrm{Sing}}),c_2(\conor_{\mathrm{Sing}}))  \notag \\
&=\sum_{\ell=1}^\infty \frac{(-1)^\ell B_\ell(0)}{\ell!}\kat'_{\ell-1} 
+ \frac{|\Aut(\Gat)|}{2} \sum_{\Gatc\in\bGatcgn(\bm)} i_{\Gat *}\pi_{\Gat *}P(c_1(\conor_{D_{\Gat/\ce}}),c_2(\conor_{D_{\Gat/\ce}}))\notag \\
&=\sum_{\ell=1}^\infty \frac{(-1)^\ell B_\ell(0)}{\ell!}\kat'_{\ell-1} 
+ \frac{|\Aut(\Gat)|}{2\deg(\mut_\Gatc)} \sum_{\Gatc\in\bGatcgn(\bm)} i_{\Gat *}\pi_{\Gat *}\mut_{\Gatc *}\mut_{\Gatc}^* P(c_1(\conor_{D_{\Gat/\ce}}),c_2(\conor_{D_{\Gat/\ce}}))\notag \\
&=\sum_{\ell=1}^\infty \frac{(-1)^\ell B_\ell(0)}{\ell!}\kat'_{\ell-1} 
+ \frac{|\Aut(\Gat)|}{2\deg(\mut_\Gatc)} \sum_{\Gatc\in\bGatcgn(\bm)}   i_{\Gat *}\pi_{\Gat *}\mut_{\Gatc *} P(c_1(\mut_{\Gatc}^*\conor_{D_{\Gat/\ce}}),c_2(\mut_{\Gatc}^*\conor_{D_{\Gat/\ce}}))\notag \\
&=\sum_{\ell=1}^\infty \frac{(-1)^\ell B_\ell(0)}{\ell!}\kat'_{\ell-1} 
+ \frac{|\Aut(\Gat)|}{2\deg(\mut_\Gatc)} \sum_{\Gatc\in\bGatcgn(\bm)}   \rho_{\Gatc *} \pi_{\Gatc *}P(c_1(\mut_{\Gatc}^*\conor_{D_{\Gat/\ce}}),c_2(\mut_{\Gatc}^*\conor_{D_{\Gat/\ce}}))\notag \\
&=\sum_{\ell=1}^\infty \frac{(-1)^\ell B_\ell(0)}{\ell!}\kat'_{\ell-1} 
+ \frac{1}{4} \sum_{\Gatc\in\bGatcgn(\bm)}  \rho_{\Gatc *} \left(\sigma_+^*(P(c_1(\mut_{\Gatc}^*\conor_{D_{\Gat/\ce}}),c_2(\mut_{\Gatc}^*\conor_{D_{\Gat/\ce}}))\right. \notag \\
& \left. \qquad + \sigma_-^*(P(c_1(\mut_{\Gatc}^*\conor_{D_{\Gat/\ce}}),c_2(\mut_{\Gatc}^*\conor_{D_{\Gat/\ce}}))\right)\notag \\
&=\sum_{\ell=1}^\infty \frac{(-1)^\ell B_\ell(0)}{\ell!}\kat'_{\ell-1} 
+ \frac{1}{2} \sum_{\Gatc\in\bGatcgn(\bm)}   \rho_{\Gatc *} \left(P(\psitp+\psitm,\psitp\psitm)\right) \notag \\
&=\sum_{\ell=1}^\infty \frac{(-1)^\ell B_\ell(0)}{\ell!}\left[\kat'_{\ell-1} 
+ \frac{1}{2} \sum_{\Gatc\in\bGatc(\bm)} \rho_{\Gatc *} \left(\frac{\psitp^{\ell-1} +  \psitm^{\ell-1}}{\psitp+ \psitm}\right)\right].\notag
\end{align*}
Now using Equation~(\ref{eq:MumKappaAC}) to relate $\kat_a$  and Mumford's $\kat'_a$, we have 
\begin{multline}
\Ch(R\pi_* \omega_\pi) =  \sum_{\ell=1}^\infty \frac{(-1)^\ell B_\ell(0)}{\ell!}\left[\kat_{\ell-1} 
 - \sum_{i=1}^n \frac{|G|}{r_i}\psit_i^{\ell-1}
\right. \\ 
\left. 
 +\frac12 \sum_{\Gatc\in\bGatc(\bm)} \rho_{\Gatc *} \left(\frac{\psitp^{\ell-1} +  \psitm^{\ell-1}}{\psitp+ \psitm}\right)\right].
 \end{multline}
Applying the relations between $\psit_i$ and $\psi_i$ given in Equation~(\ref{eq:psit}) and the relations between $\kat$ and $\kappa$ given in Equation~(\ref{eq:kappat}) we get,
\begin{eqnarray*}
\Ch(R\pi_* \omega_\pi) &=&  \sum_{\ell=1}^\infty \frac{(-1)^\ell B_\ell(0)}{\ell!}\left[|G|\kappa_{\ell-1} 
 - \sum_{i=1}^n \frac{|G|}{r_i^\ell}\psi_i^{\ell-1}
\right. \\
 &+& \qquad \left.\frac12 \sum_{\Gatc\in\bGatc(\bm)}\left(\frac{1}{\rpp}\right)^{\ell-2} \rho_{\Gatc *} \left(\frac{\psip^{\ell-1} +  \psim^{\ell-1}}{\psip+ \psim}\right)\right].
\end{eqnarray*}
Finally, we have $R^1\pi_*(\omega_\pi) = \alpha\co$, where $\alpha$ is the number of connected components of a fiber of $\ce$ over a general point of a given connected component $\N_{g,n}(\bm)$ of $\MM_{g,n}$.  Since $\alpha = |G/G_0|$, this  finishes the proof.
\end{proof}

\subsection{Computation of the the Chern character using Corollary~\ref{cor:BF}}

In this section, we will compute the $\Rep(G)$-valued Chern character 
$\bCh(\crrt)$ 
of the Hurwitz-Hodge bundle using Corollary~\ref{cor:BF}.

\begin{df}
Define the formal power series
\begin{equation}\label{df:DD}
\DD(u) := \frac{1}{B(u,0)} = \frac{e^u-1}{u}
\end{equation}
in $\nq[[u]]$. \end{df}
\begin{rem}
The Todd class $\Td(L)$ of a line bundle $L$ is
$B(-c_1(L),0) = \DD^{-1}(-c_1(L))$.
\end{rem}

\begin{lm}
For all $\Gatc\in \bGatcgn(\bm,\mpp,\mmm)$, let $\cF$ be any element of $K_G(\MM_\Gatc)$.  We have
\begin{equation}\label{eq:bChGRR}
\bCh\left(\rho_{\Gatc *}\cF\right) = \rho_{\Gatc *} \left[\bCh\left(
\cF\right)\DD\left(\frac{\psip+\psim}{|\mpp|}\right)
\right].
\end{equation}
\end{lm}
\begin{proof}
Throughout this proof we will refer to Diagram~(\ref{diag:LRRnode}).
Let $T_f = \Omega_f^*$ denote the relative tangent bundle of a morphism $f$ and
$\nor_i = -T_i$ denote the normal bundle of a regular embedding $i$.  We have
\begin{eqnarray*}
\bCh(\rho_{\Gatc *}\cF) 
&=& \frac1{\deg(\pi_\Gatc)} \bCh(\rho_{\Gatc *} \pi_{\Gatc *}\pi_{\Gatc}^* \cF) \\
&=& \frac1{\deg(\pi_\Gatc)} \bCh(\pi_* \itt_{\Gat *} \jt_{\Gat *} \mut_{\Gatc *}\pi_{\Gatc}^* \cF) \\
&=& \frac1{\deg(\pi_\Gatc)} \pi_*\left( \bCh( \itt_{\Gat *} \jt_{\Gat *}\mut_{\Gatc *} \pi_{\Gatc}^* \cF)\Tdch(\Omega_\pi)\right)\\
&=& \frac1{\deg(\pi_\Gatc)} \pi_*\left(\itt_{\Gat *}\jt_{\Gat *}\left( \bCh(\mut_{\Gatc *} \pi_{\Gatc}^* \cF)\Td(-\nor_{\itt_{\Gat}\jt_{\Gat}})\right)\Tdch(\Omega_\pi)\right)\\
&=& \frac1{\deg(\pi_\Gatc)} \pi_*\left(\itt_{\Gat *}\jt_{\Gat *}\left(\mut_{\Gatc *} \bCh( \pi_{\Gatc}^* \cF)\Td(-\nor_{\itt_{\Gat}\jt_{\Gat}})\right)\Tdch(\Omega_\pi)\right)\\
&=& \frac1{\deg(\pi_\Gatc)} \pi_*\itt_{\Gat *}\jt_{\Gat *}\mut_{\Gatc *} \left(\bCh( \pi_{\Gatc}^* \cF)\mut_{\Gatc}^*\Td(-\nor_{\itt_{\Gat}\jt_{\Gat}})\mut_{\Gatc}^*\jt_{\Gat}^*\itt_{\Gat}^*\Tdch(\Omega_\pi)\right)\\
&=& \frac1{\deg(\pi_\Gatc)} \rho_{\Gatc *}\pi_{\Gatc *} \left(\bCh( \pi_{\Gatc}^* \cF)\mut_{\Gatc}^*\Td(-\nor_{\itt_{\Gat}\jt_{\Gat}})\mut_{\Gatc}^*\jt_{\Gat}^*\itt_{\Gat}^*\Tdch(\Omega_\pi)\right),
\end{eqnarray*}
where the third and fourth equalities use the Grothendieck-Riemann-Roch Theorem, the
fifth uses the fact that $\mut_\Gatc$ is \'etale, and the sixth uses the projection formula.
Using Equations~(\ref{eq:pistarsigma}) and (\ref{eq:ConorPullback}), we now have
\begin{eqnarray*}
\bCh(\rho_{\Gatc *}\cF) 
&=& \frac{|Z_G(\mpp)|}{\rpp\deg(\pi_\Gatc)} \rho_{\Gatc *}(\sig_+^*+\sig_-^*) \left(\bCh( \pi_{\Gatc}^* \cF)\mut_{\Gatc}^*\Td(-\nor_{\itt_{\Gat}\jt_{\Gat}})\mut_{\Gatc}^*\jt_{\Gat}^*\itt_{\Gat}^*\Tdch(\Omega_\pi)\right)\\
&=& \frac{1}{2} \rho_{\Gatc *} \left[\left(\bCh(  \cF)\sig_+^*\mut_{\Gatc}^*\Td(-\nor_{\itt_{\Gat}\jt_{\Gat}})\sig_+\mut_{\Gatc}^*\jt_{\Gat}^*\itt_{\Gat}^*\Tdch(\Omega_\pi)\right)\right.\\
&& \left. \qquad + \left(\bCh( \cF)\sig_-^*\mut_{\Gatc}^*\Td(-\nor_{\itt_{\Gat}\jt_{\Gat}})\sig_-\mut_{\Gatc}^*\jt_{\Gat}^*\itt_{\Gat}^*\Tdch(\Omega_\pi)\right)\right]\\
&=& \frac{1}{2} \rho_{\Gatc *} \left[\left(\bCh(  \cF)\Tdch(-\cltp-\cltm))\sig_+^*\mut_{\Gatc}^*\jt_{\Gat}^*\itt_{\Gat}^*\Tdch(\Omega_\pi)\right)\right.\\
&& \left. \qquad + \left(\bCh( \cF)\Tdch(-\cltp-\cltm))\sig_-^*\mut_{\Gatc}^*\jt_{\Gat}^*\itt_{\Gat}^*\Tdch(\Omega_\pi)\right)\right].
\end{eqnarray*}
By a simple argument given in \cite[\S5]{Mum}, the term $\Tdch(\Omega_\pi)$ can be written as
$$\Tdch(\Omega_\pi) = \Tdch(\omega_\pi) \Tdch(-\co_{\mathrm{Sing}}),$$
where $\co_{\mathrm{sing}}:= \sum_{\Gat} \itt_{\Gat *} \jt_{\Gat *} \co_{D_\Gat}$.  But by the residue map, we also have $\jt^*\itt^*\omega_\pi \cong \co_{D_\Gat}$, from which we deduce
\begin{eqnarray*}
\sigpm^*\mut_{\Gatc}^*\jt_{\Gat}^*\itt_{\Gat}^*\Tdch(\Omega_\pi) &=& \Tdch(-\sigpm^*\mut_{\Gatc}^*\jt_{\Gat}^*\itt_{\Gat}^*\itt_{\Gat *}\jt_{\Gat *}\co_{D_\Gat})\\
&=& \Tdch(-\sigpm^*\mut_{\Gatc}^*\lambda_{-1}(\nor_{\itt_{\Gat}\jt_{\Gat}}))\\
&=& \Tdch(-\lambda_{-1}(\sigpm^*\mut_{\Gatc}^*\nor^*_{\itt_{\Gat}\jt_{\Gat}}))\\
&=& \Tdch(-\lambda_{-1}(\cltp + \cltm))\\
&=& \Tdch\left(-(1-\cltp) (1-\cltm)\right)\\
&=& \Tdch\left(-1+\cltp +\cltm - \cltp\otimes\cltm\right)\\
&=& \frac{\DD(\psitp+\psitm)}{\DD(\psitp)\DD(\psitm)}.
\end{eqnarray*}
Plugging this back into our earlier calculation gives the desired result.
\end{proof}
\begin{thm}\label{thm:bChtS}
The Chern character,  $\bCh(\crrt)$, of the dual Hurwitz-Hodge bundle in
$A^\bullet(\MM_{g,n}(\bm))\otimes\Rep(G)$ is
\begin{equation}\label{eq:bChR}
\bCh(\crrt) = \bone\otimes\nc[G/G_0] + (\Ch(\crr)-\bone)\otimes\nc[G] +
\sum_{i=1}^n \bCh(\cs_{m_i}) + \sum_{\Gatc\in \bGatcgn(\bm)} \bCh(\cs_\Gatc),
\end{equation}
where $\bCh(\cs_{m_i})$ and $\bCh(\cs_\Gatc)$ are given explicitly as follows:

For all $i=1,\dots,n$ we have
\begin{eqnarray}\label{eq:bChtS}
\bCh(\cs_{m_i}) &=& 
\Ind^G_{\gen{m_i}}\left(\frac{1}{e^{\psi_i/r_i} \V_{m_i}-1} - \nc[\gen{m_i}]
\frac{1}{e^{\psi_i}-1} \right) \\
&=&\sum_{k=0}^{r_i-1} \frac{\exp(k \psi_i/r_i)-1}{\exp(\psi_i)-1} 
\Ind^G_{\gen{m_i}} \V_{m_i}^k \\
&=& \sum_{k=0}^{r_i-1} \sum_{j=0}^\infty \frac{B_{j+1}(k/r_i) - B_{j+1}(0)}{(j+1)!}
\psi_i^j \Ind^G_{\gen{m_i}} \V_{m_i}^k
\end{eqnarray}
in 
$A^\bullet(\MM_{g,n}(\bm))\otimes\Rep(G)$ for all $i=1,\ldots, n$. 

Similarly, for all $\Gatc$ in $\bGatcgn(\bm,\mpp,\mmm)$ we have
\begin{eqnarray}\label{eq:bChtSGat}
\bCh(\cs_\Gatc) 
&=& -\frac{\rpp}{2|G|}  \rho_{\Gatc
*}\left(\Ind^G_{\gen{\mpp}} \I^{\gen{\mpp}}(\bCh(\cst_{\mpp})
\bCh(\cst_{\mmm})\DD((\psip+\psim)/\rpp))\right)\notag 
\end{eqnarray}
in $A^\bullet(\MM_{g,n}(\bm))\otimes\Rep(G),$ 
and
\begin{equation}\label{eq:bChttStp} 
\bCh(\cst_{m_\pm}) = \sum_{k=0}^{\rpp-1} 
\frac{\exp(k\psi_{\pm}/\rpp)-1}{\exp(\psi_\pm)-1}  \V_{m_\pm}^k
= \sum_{k=0}^{\rpp-1} \sum_{j=0}^\infty \frac{B_{j+1}(k/\rpp)-B_{j+1}(0)}
{(j+1)!}\psi_\pm^j \V_{m_\pm}^k.
\end{equation}

In addition, $\bCh(\cs_\Gatc)$ can be rewritten as
\begin{equation}\label{eq:SimplebChcs}
\bCh(\cs_{\Gatc}) = -\frac{\rpp^2}{2|G|} \sum_{k=0}^{\rpp-1}\sum_{j\geq 1}
\Ind^G_{\gen{\mpp}}\V_\mpp^k \delta B_{j+1}(k/\rpp)
\sum_{\substack{\jjp+\jjm = j-1 \\ \jjpm\geq 0}} (-1)^{\jjm} \rho_{\Gatc
*}\left(\psip^\jjp\psim^\jjm\right) 
\end{equation}
\end{thm}

\begin{proof}
Equation (\ref{eq:bChtS}) follows from the definition of $\cs_{m_i}$
given in Equation (\ref{eq:Suv}), from the relation between $\psi_i$ and
$\psit_i$ given in Proposition (\ref{prop:PsitPsi}).

Equation (\ref{eq:bChtSGat}) follows from Equation
(\ref{eq:bChGRR}) and the fact
that the Chern character is a ring 
homomorphism.

We will now prove Equation (\ref{eq:SimplebChcs}). Equation (\ref{eq:bChGRR})
implies that 
\[
\bCh\left(\rho_{\Gatc *}\left(\I^{\gen{\mpp}}\cst_{\mpp}\cst_{\mmm}\right)\right) 
= \rho_{\Gatc *} \left( 
\IF_{\rpp}(e^{\psip/\rpp},e^{\psim/\rpp},\V_\mpp)
\Delta\left((\psip+\psim)/\rpp\right)
\right).
\]
We have
\begin{equation*}
\begin{split}
&\IF_{\rpp}(e^{\psip/\rpp},e^{\psim/\rpp},\V_\mpp)
\Delta\left((\psip+\psim)/\rpp\right) \\
&= \sum_{k=1}^{\rpp-1} \V_\mpp^k
\frac{\rpp}{\psip+\psim}\left(-1 + \frac{e^{k\psip/\rpp}-1}{e^{\psip}-1} + 
\frac{e^{(\rpp-k)\psim/\rpp}-1}{e^{\psim}-1}\right) \\
&= \sum_{k=1}^{\rpp-1} \V_\mpp^k
\frac{\rpp}{\psip+\psim}\left(-1  
+\sum_{j\geq 0}\left[\delta B_{j+1}(k/\rpp)\psip^j 
+ \delta B_{j+1}(1-k/\rpp)\psim^j \right]
\right) \\
&= \sum_{k=1}^{\rpp-1} \V_\mpp^k
\frac{\rpp}{\psip+\psim}\left(-1 + \delta B_1(k/\rpp) + \delta B_1(1-k/\rpp)
+\sum_{j\geq 1}\delta B_{j+1}(k/\rpp)\left(\psip^j 
+ (-1)^{j+1}\psim^j \right)
\right) \\
&= \sum_{k=1}^{\rpp-1} \V_\mpp^k 
\sum_{j\geq 1}\rpp \delta B_{j+1}(k/\rpp) \frac{\psip^j 
+ (-1)^{j+1}\psim^j}{\psip+\psim}
\\
\end{split}
\end{equation*}
where we have plugged in the definitions and canceled the numerator of
$\Delta$ with the denominator of $\IF_{\rpp}$ in the first equality and we
have used Equation (\ref{eq:DeltaB}) in the second. In the third and fourth
equalities, we have used  Equation~(\ref{eq:deltaBerRel}).
But for all $j\geq 1$, 
\begin{equation}
\frac{\psip^j + (-1)^{j+1}\psim^j}{\psip+\psim} = \sum_{\substack{\jjp+\jjm=
j-1\\ \jjpm\geq 0}} \psip^\jjp\psim^\jjm (-1)^\jjm.
\end{equation}
Plugging this in and then applying $-\frac{\rpp}{2 |G|}\Ind^G_{\gen{\mpp}}$ 
yields Equation (\ref{eq:SimplebChcs}).
\end{proof} 

\subsection{Relating the Chern characters of the Hurwitz-Hodge bundle}
It is not \emph{a priori} obvious that the two ways of computing the Chern
character are consistent.  We will now show that one is a special case of the
other, that is to say, that the Chern character $\Ch:\Kmgnloc{}\otimes\Rep(G)\rTo A^\bullet(\MM_{g,n})$
can be obtained from the $\Rep(G)$-valued Chern character $\bCh$ via
\begin{equation}
\Ch = \chi_1\circ\bCh.
\end{equation}

\begin{prop}
Applying $\chi_1$ to Equation (\ref{eq:bChR}) yields Theorem
\ref{thm:GRRChernCh}. We also have the following equalities in
$A^\bullet(\MM_{g,n}(\bm))$
\begin{equation}\label{eq:Chcs}
\Ch(\cs_{m_i}) = \sum_{j\geq 1} |G| (r_i^{-j}-1) \psi_i^{j-1} \frac{B_j(0)}{j!}
\end{equation}
\begin{equation}\label{eq:Chcsn}  
\Ch(\cs_\Gatc) = -\frac{1}{2} \rho_{\Gatc *}\left(\sum_{j=1}^\infty 
\frac{B_j(0)}{j!} (\rpp^{2-j} - \rpp^2) \sum_{q=0}^{j-2} (-1)^q \psip^q
\psim^{j-2-q} \right),
\end{equation}
for all $\Gatc$ in $\bGatcgn(\bm,\mpp,\mmm)$, 
and
\begin{equation}\label{eq:Chhod}
\Ch (\crr) = g + \sum_{j=1}^\infty \frac{B_{j}(0)}{j!}\left[-\kappa_{j-1}
+ \sum_{i=1}^n \psi_i^{j-1} - \sum_{\Gatc\in \bGatcgn(\bm)}
  \frac{\rpp^2}{2 |G|}\rho_{\Gatc *}\left(\sum_{q=0}^{j-2} (-1)^q
  \psip^q\psim^{j-2-q}\right)\right],
\end{equation}
where it is understood that the sum over $q$ vanishes when $j=1$ in the last
two equations.
\end{prop}
\begin{proof}
We first observe that for all $m$ in $G$ and $k=0,\ldots,|m|-1$, 
\begin{equation}\label{eq:RankOfInduced}
\chi_1(\Ind^G_{\gen{m}} \V_m^k ) = \frac{|G|}{|m|}.
\end{equation}

Comparing Equations (\ref{eq:ChiOneSm}), (\ref{eq:ChiOneSmExpanded}), and
Equation (\ref{eq:bChtS}) yields Equation (\ref{eq:Chcs}).

To prove Equation (\ref{eq:Chcsn}), we apply $\chi_1$ to Equation
(\ref{eq:SimplebChcs}) to obtain
\begin{equation*}
\begin{split}
&\Ch(\cs_\Gatc) = \chi_1(\bCh(\cs_\Gatc)) \\
&= \frac{|G|}{\rpp} \left(-\frac{\rpp^2}{2|G|}\right)
\sum_{k=0}^{\rpp-1}\sum_{j\geq 1} \delta B_{j+1}(k/\rpp)
\sum_{\substack{\jjp+\jjm = j-1 \\ \jjpm\geq 0}} (-1)^{\jjm} \rho_{\Gatc
*}\left(\psip^\jjp\psim^\jjm\right)  \\
&= -\frac{\rpp}{2} \sum_{j\geq 1} \left[\sum_{k=0}^{\rpp-1}
\delta B_{j+1}(k/\rpp)\right] \sum_{\substack{\jjp+\jjm =
j-1 \\ \jjpm\geq 0}} (-1)^{\jjm} \rho_{\Gatc
*}\left(\psip^\jjp\psim^\jjm\right)  \\ 
&= -\frac{\rpp}{2} \sum_{j\geq 1} \frac{B_{j+1}(0)}{(j+1)!} (\rpp^{-j}-\rpp)
 \sum_{\substack{\jjp+\jjm = j-1 \\ \jjpm\geq 0}} (-1)^{\jjm} \rho_{\Gatc
*}\left(\psip^\jjp\psim^\jjm\right)  \\ 
&=  \sum_{j\geq 1} -\frac{\rpp^{1-j}-\rpp^2}{2}\frac{B_{j+1}(0)}{(j+1)!} 
 \sum_{\substack{\jjp+\jjm = j-1 \\ \jjpm\geq 0}} (-1)^{j-1-\jjp} \rho_{\Gatc
*}\left(\psip^\jjp\psim^\jjm\right)  \\ 
&= \sum_{j\geq 1} -\frac{\rpp^{1-j}-\rpp^2}{2} \frac{B_{j+1}(0)}{(j+1)!} 
 \sum_{\substack{\jjp+\jjm = j-1 \\ \jjpm\geq 0}} (-1)^{\jjp} \rho_{\Gatc
*}\left(\psip^\jjp\psim^\jjm\right)  \\ 
&= \sum_{j\geq 2} -\frac{\rpp^{2-j}-\rpp^2}{2} \frac{B_{j}(0)}{j!} 
 \sum_{\substack{\jjp+\jjm = j-2 \\ \jjpm\geq 0}} (-1)^{\jjp} \rho_{\Gatc
*}\left(\psip^\jjp\psim^\jjm\right)  \\ 
\end{split}
\end{equation*}
where we have used Equation (\ref{eq:BernoulliSum}) in the fourth equality,
and have used the fact that $B_{j+1}(0) = 0$ unless $j\geq 1$ is odd in the
sixth equality. This proves Equation (\ref{eq:Chcsn}). 

Finally, on $\M_{g,n}$, we have the usual formula due to Mumford
\begin{equation}
\Ch(\crr) = g + \sum_{\ell=2}^\infty
\frac{B_{\ell}(0)}{\ell!}\left[-\kappa_{\ell-1} + \sum_{i=1}^n
\psi_i^{\ell-1} - \frac{1}{2} \sum_\Gac \rho_{\Gac*}\left(\sum_{q=0}^{\ell-2} (-1)^q
  \psip^q\psim^{\ell-2-q}\right)\right]. 
\end{equation}
Applying $\st^*$, we obtain Equation (\ref{eq:Chhod}) after using 
Equation~(\ref{eq:StRhoFullPushPull}), which in this case gives 
\begin{equation}\label{eq:TylerCheckThis}
\st^*\left(\frac{1}{2}\sum_\Gac \rho_{\Gac *} (\psip^a\psim^b)\right)  =
\sum_\Gatc \frac{\rpp^2}{2 |G|} \rho_{\Gatc *}(\psip^a\psim^b).
\end{equation}

Now, plug in Equations (\ref{eq:Chcs}), (\ref{eq:Chcsn}), (\ref{eq:Chhod})
into Equation (\ref{eq:bChtS}) and use that the rank of $\nc[G/G_0]$ is
$\frac{|G|}{|G_0|}$ and we obtain Theorem \ref{thm:GRRChernCh}.
\end{proof}

\begin{crl}
We have the identities
\begin{equation}\label{eq:lat}
\begin{split}
\bCh_1(\crrt) &= -\frac{1}{12}\nc[G] \otimes\kappa_1  +
\sum_{i=1}^n \left(\frac{1}{12}\nc[G] +\sum_{k=0}^{|m_i|-1}
\Ind^G_{\gen{m_i}} \V_{m_i}^k \frac{k (k-|m_i|)}{2 |m_i|^2}\right)
\otimes \psi_i \\ 
&- \sum_{\Gatc\in \bGatcgn(\bm)}\left(\frac{|\mpp|^2}{24|G|} \nc[G]
+\sum_{k=0}^{|\mpp|-1} \frac{k (k-|\mpp|)}{4 |G|}
\Ind^G_{\gen{\mpp}}\V_\mpp^k\right) \otimes \rho_{\Gatc *}(\one)\\ 
\end{split}
\end{equation}
and
\begin{equation}\label{eq:GACKappa}
\Ch_1(\crrt) = -\frac{|G|}{12}\kappa_1  + \sum_{i=1}^n \frac{|G|}{12 |m_i|^2}
\psi_i - \sum_{\Gatc\in\bGatcgn(\bm)}\frac{1}{24} \rho_{\Gatc*}(\one).
\end{equation}
\end{crl}

\begin{proof}

Since $B_2(x)-B_2(0) = x^2-x$, Equation (\ref{eq:SimplebChcs}) yields
\begin{equation*}
\sum_{\Gatc} \bCh_1(\cs_\Gatc) = \sum_{\Gatc} -\frac{1}{4|G|}
\sum_{k=0}^{\rpp-1} k
(k-|\mpp|)
\Ind_{\gen{m}}^G \V_\mpp^k \otimes \rho_{\Gatc
*}(\one).
\end{equation*}
Now plug in 
\begin{equation}
\Ch_1(\crr) = \frac{1}{12}(-\kappa_1+\sum_{i=1}^n\psi_i - \sum_{\Gatc}
\frac{|\mpp|^2}{2|G|} \rho_{\Gatc*}(\one))
\end{equation}
from part of Equation (\ref{eq:Chhod}) to obtain the desired result. 

The second equation can be obtained from the first by applying $\chi_1$.
\end{proof}
\begin{rem}
In the special case that $G=\{1\}$, Equation~(\ref{eq:GACKappa}) reduces to the well-known relation \cite[pg.~102]{MumLEns} 
\begin{equation}\label{eq:lambdakappa}
12 \lambda_1 = \kappa_1 - \sum_{i=1}^n \psi_i + \frac12\sum_{\Gac\in\bGacgn} \rho_{\Gac *} (\one).
\end{equation}
 To see this, recall that $\lambda_1 :=
-c_1(\crrt)$ and that Mumford uses the class we call $\kappa'_1$ instead of
$\kappa_1$, which is why the equation here differs from his by the term
$\sum_{i=1}^n \psi_i$.  
\end{rem}

\bibliographystyle{amsplain}

\providecommand{\bysame}{\leavevmode\hboxto3em{\hrulefill}\thinspace}

\end{document}